\documentclass{article}
\usepackage{amsmath, amssymb}
\usepackage{index}
\usepackage{t1enc}
\usepackage[english]{babel}
\newtheorem{theorem}{Theorem}[section]
\newtheorem{lemma}[theorem]{Lemma}
\newtheorem{corollary}[theorem]{Corollary}
\newtheorem{proposition}[theorem]{Proposition}
\newtheorem{definition}[theorem]{Definition}

\newtheorem{remark}[theorem]{Remark}

\makeindex

\newindex{sub}{sdx}{snd}{Subject Index}
\newindex{aut}{adx}{and}{Author Index}
\newindex{not}{ndx}{nnd}{Notation Index}

\renewcommand{\abstractname}{Abstract.}
\renewcommand\abstract{\hfil\break\topsep=0pt\partopsep=0pt\parsep=0pt\itemsep=0pt\relax
\trivlist\item[\hskip\labelsep
{\bfseries\abstractname}]\if!\abstractname!\hskip-\labelsep\fi}
\newcommand{\email}[1]{{(e-mail: #1)}}

\def\title#1{\hfil\break\hfil\break
\hfil\break\par\addvspace\baselineskip\noindent
\ignorespaces{\LARGE\bf#1}\hfil\break}
\def\author#1{\par\addvspace\baselineskip\noindent
\ignorespaces{\large\bf#1}}

\begin{document}

\title{Pricing, Hedging and Optimally Designing Derivatives\vspace{3mm}\\
via Minimization of Risk Measures \footnote{The authors would like
to thank their respective husband, Laurent Veilex and Faycal El
Karoui, for their constant support and incredible patience.}}
\begin{center}{\large\today{}}
\end{center}
\author{Pauline Barrieu \footnote{L.S.E., Statistics Department, Houghton Street, London,
WC2A 2AE, United Kingdom. \email{p.m.barrieu@lse.ac.uk}\\
Research supported in part by the EPSRC, under Grant
GR-T23879/01.} and Nicole El Karoui \footnote{C.M.A.P., Ecole
Polytechnique, 91128 Palaiseau C\'{e}dex, France.
\email{elkaroui@cmapx.polytechnique.fr}}}
\vspace{2cm}\\
The question of pricing and hedging a given contingent claim has a
unique solution in a complete market framework. When some
incompleteness is introduced, the problem becomes however more
difficult. Several approaches have been adopted in the literature
to provide a satisfactory answer to this problem, for a particular
choice criterion. Among them, \index[aut]{Hodges} Hodges and
\index[aut]{Neuberger} Neuberger \cite{HodgesNeuberger} proposed
in 1989 a method based on utility maximization. The price of the
contingent claim is then obtained as the smallest (resp. largest)
amount leading the agent indifferent between selling (resp.
buying) the claim and doing nothing. The price obtained is the
indifference seller's (resp. buyer's) price. Since then, many
authors have used this approach, the exponential utility function
being most often used (see for instance, \index[aut]{El Karoui} El
Karoui and \index[aut]{Rouge} Rouge \cite{ElKarouiRouge},
\index[aut]{Becherer} Becherer \cite{Becherer}, Delbaen et al.
\cite{Delbaenetal}\index[aut]{Delbaen}\index[aut]{Grandits}
\index[aut]{Rheinlander}\index[aut]{Samperi}\index[aut]{Schweizer}\index[aut]{Stricker},
\index[aut]{Musiela} Musiela and \index[aut]{Zariphopoulou}
Zariphopoulou
 \cite{MZ1} or \index[aut]{Mania} Mania and \index[aut]{Schweizer} Schweizer \cite{ManiaSchweizer}...).\vspace{1mm}\\
In this chapter, we also adopt this exponential utility point of
view to start with in order to find the optimal hedge and price of
a contingent claim based on a non-tradable risk. But soon, we
notice that the right framework to work with is not that of the
exponential utility itself but that of the certainty equivalent
which is a convex functional satisfying some nice properties among
which that of cash translation invariance. Hence, the results
obtained in this particular framework can be immediately extended
to functionals satisfying the same properties, in other words to
convex risk measures as introduced by \index[aut]{F\"{o}llmer}
F\"{o}llmer and \index[aut]{Schied} Schied
\cite{FoellmerSchied02a} and \cite{FoellmerSchied02b} or by
\index[aut]{Frittelli} Frittelli and \index[aut]{Gianin} Gianin
\cite{FrittelliGianin02}. Starting with a utility maximization
problem, we end up with an equivalent risk measure minimization in
order to
price and hedge this contingent claim.\vspace{1mm}\\
Moreover, this hedging problem can be seen as a particular case of
a more general situation of risk transfer between different
agents, one of them consisting of the financial market. Therefore,
we consider in this chapter the general question of optimal
transfer of a non-tradable risk and specify the results obtained
in the particular situation
of an optimal hedging problem.\vspace{1mm}\\
Both static and dynamic approaches are considered in this chapter,
in order to provide constructive answers to this optimal risk
transfer problem. Quite recently, many authors have studied
dynamic version of static risk measures (see for instance, among
many other references, \index[aut]{Cvitanic} Cvitanic and
\index[aut]{Karatzas} Karatzas \cite{CvitanicKaratzas99},
\index[aut]{Scandolo} Scandolo \cite{Scandolo}, \index[aut]{Weber}
Weber \cite{Weber}, \index[aut]{Artzner} \index[aut]{Delbaen}
\index[aut]{Eber} \index[aut]{Heath}Artzner et al.
\cite{ArtznerDelbaenEberHeath2}, \index[aut]{Cheridito} Cheridito,
\index[aut]{Delbaen} Delbaen and \index[aut]{Kupper} Kupper
\cite{CheriditoDelbaenKupper1}, \index[aut]{Frittelli} Frittelli
and \index[aut]{Gianin} Gianin \cite{FrittelliGianin04},
\index[aut]{Gianin} Gianin \cite{Gianin}, \index[aut]{Riedel}
Riedel \cite{Riedel} or \index[aut]{peng} Peng \cite{Peng03}).
When considering a dynamic framework, our main purpose is to find
a trade-off between static and very abstract risk measures as we
are more interested in tractability issues and interpretations of
the dynamic risk measures we obtain rather than the ultimate
general results. Therefore, after introducing a general axiomatic
approach to dynamic risk measures, we relate the dynamic version
of convex risk measures to BSDEs. For the sake of a better
understanding, a whole section in the second part is dedicated to
some key results and properties of BSDEs, which are essential to
this definition of dynamic convex risk measures.
\begin{center}
\section*{Part I: Static Framework}
\end{center}
In this chapter, we focus on the question of optimal hedging of a
given risky position in an incomplete market framework. However,
instead of adopting a standard point of view, we look at it in
terms of an optimal risk transfer between different economic
agents, one of them being
possibly a financial market.\\
The risk that we consider here is not (directly) traded on any
financial market. We may think for instance of a weather risk, a
catastrophic risk (natural catastrophe, terrorist attack...) but
also of any global insurance risk that may be securitized, such as
the
longevity of mortality risk...\\
First adopting a static point of view, we proceed in several
steps. In a first section, we relate the notion of indifference
pricing rule to that of transaction feasibility, capital
requirement, hedging and naturally introduce convex risk measures.
Then, after having introduced some key operations on convex risk
measures, in particular the dilatation and the inf-convolution, we
study the problem of optimal risk transfer between two agents. We
see how the risk transfer problem can be reduced to an
inf-convolution problem of convex functionals. We solve it
explicitly in the dilated framework
and give some necessary and sufficient conditions in the general framework.\vspace{2mm}\\
\section{Indifference Pricing, Capital Requirement and Convex Risk Measures\label{section risk measure et pricing}}
As previously mentioned in the introduction, since 1989 and the
seminal paper by \index[aut]{Hodges} Hodges and
\index[aut]{Neuberger} Neuberger \cite{HodgesNeuberger}
indifference pricing based on a utility criterion has been a
popular (academic) method to value claims in an incomplete market.
Taking the buyer point of view, the \index[sub]{indifference
price} indifference price corresponds to the maximal amount $\pi$,
the agent having a utility function $u$ is ready to pay for a
claim $X$. In other words, $\pi$ is determined as the amount the
agent pays such that her expected utility remains unchanged when
doing the transaction:
$$\mathbb{E}[u(X-\pi)]=u(0).$$
This price is not a transaction price. It gives an upper bound
(for the buyer) to the price of this claim so that a transaction
will take place. $\pi$ also corresponds to the certainty
equivalent of the claim payoff $X$. Certain properties this
indifference price should have are rather obvious: first it should
be an increasing function of $X$ but also a convex function in
order to take into account the diversification aspect of
considering a portfolio of different claims rather than the sum of
different individual portfolios. Another property which is rather
interesting is the \index[sub]{translation invariance} cash
translation invariance property. More precisely, it seems natural
to consider the situation where translating the payoff of the
claim $X$ by a non-risky amount $m$ simply leads to a translation
of the price $\pi$ by the same amount. It is the case, as we will
see for the exponential utility in the following subsection.
\subsection{The Exponential Utility Framework}
First, let us notice that exponential utility functions have been
widely used in the financial literature. Several facts may justify
their relative importance compared to other utility functions but,
in particular, the absence of constraint on the sign of the future
considered cash flows and its relationship with probability
measures make them very convenient to use.
\subsubsection{Indifference Pricing Rule
\label{sssect:indifference pricing rule}} \index[sub]{indifference
pricing rule}
In this introductory subsection, we simply consider an agent,
having an exponential utility function $U(x) =-\gamma \exp \big( -\frac{1}{%
\gamma }x\big)$, where $\gamma$ is her \index[sub]{risk tolerance
coefficient}risk tolerance coefficient. She evolves in an
uncertain universe modelled by a standard probability space $(
\Omega ,\Im ,\mathbb{P}) $ with time horizon $T$. The wealth $W$
of the agent at this future date $T$ is uncertain, since $W$ can
be seen as a particular position on a given portfolio or as the
book of the agent. To reduce her risk, she can decide whether or
not to buy a contingent claim with a payoff $X$ at time $T$. For
the sake of simplicity, we neglect interest rate between $0$ and
$T$ and assume that both random variables $W$ and $X$ are bounded.
\\
In order to decide whether or not she will buy this claim, she
will find the maximum price she is ready to pay for it, her
indifference price $\pi(X)$ for the claim $X$ given by the
constraint $\mathbb{E}_{\mathbb{P}}\big[U(W+X-\pi(X))\big]=
\mathbb{E}_{\mathbb{P}}\big[ U(W)\big]$. Then,
\begin{eqnarray*}
\mathbb{E}_{\mathbb{P}}\big[\exp \big( -\frac{1}{%
\gamma }(W+X-\pi(X))\big)\big]&=&\mathbb{E}_{\mathbb{P}}\big[\exp \big( -\frac{1}{%
\gamma }W\big)\big]\\
\Leftrightarrow \quad
\pi(X|W)&=&e_{\gamma}\big(W\big)-e_{\gamma}\big(W+X\big)
\end{eqnarray*}
where \index[not]{$e_{\gamma}$} $e_{\gamma}$ is the opposite of
the certainty equivalent, defined for any bounded random variable
$\Psi$
\begin{equation}
\label{definition e_gamma} e_{\gamma }\big( \Psi \big)
\triangleq \gamma \ln \mathbb{E}_{\mathbb{P}}\Big[ \exp \big( -\frac{1}{%
\gamma }\Psi \big) \Big].
\end{equation}
The indifference pricing rule $\pi(X|W)$ has the desired property
of increasing monotonicity, convexity and translation invariance:
$\pi(X+m|W)=\pi(X|W)+m$. Moreover, the functional $e_{\gamma
}(X)=-\pi(X|W=0)$ has similar properties; it is decreasing, convex
and translation invariant in the following sense: $e_{\gamma }(
\Psi +m) =e_{\gamma }( \Psi) -m$.
\subsubsection{Some Remarks on the "Price" $\pi(X)$
\label{subsection transaction prices}}
$\pi(X)$ does not correspond to a transaction price but simply
gives an indication of the transaction price range since it
corresponds to the maximal amount the agent is ready to pay for
the claim $X$ and bear the associated risk given her initial
exposure. This dependency seems quite intuitive: for instance, the
considered agent can be seen as a trader who wants to buy the
particular derivative $X$ without knowing its price. She
determines it by considering the contract
{\it relatively} to her existing book.\\
For this reason and for the sake of a better understanding, we
will temporarily denote it by $\pi^{b}(X|W)$, the upper-script "b"
standing for "buyer". This heavy notation underlines the close
relationship between the pricing rule and the actual exposure of
the agent. The considered framework is symmetric since there is no
particular requirement on the sign of the different quantities.
Hence, it is possible to define by simple analogy the
\index[sub]{indifference seller's price} {\it indifference
seller's price} of the claim $X$. Let us denote it by
$\pi^{s}(X|W)$, the upper-script "s" standing for "seller". Both
seller's and buyer's indifference pricing rules are closely
related as
$$\pi^{s}(X|W)=-\pi^{b}(-X|W).$$
Therefore, the seller's price of $X$ is simply the opposite of the
buyer's price of $-X$.\\
Such an axiomatic approach of the pricing rule is not new. This
was first introduced in insurance under the name of
\index[sub]{convex premium principle} {\it convex premium
principle} (see for instance the seminal paper of
\index[aut]{Deprez} Deprez and \index[aut]{Gerber} Gerber
\cite{DeprezGerber} in 1985) and then developed in continuous time
finance (see for instance \index[aut]{El Karoui} El
Karoui and \index[aut]{Quenez} Quenez \cite{ElKarouiQuenez95}).\\
When adopting an exponential utility criterion to solve a pricing
problem, the right framework to work with seems to be that of the
functional $e_{\gamma }$ and not directly that of utility. This
functional, called entropic risk measure, holds some key
properties of convexity, monotonicity and cash translation
invariance. It is therefore possible to generalize the utility
criterion to focus more on the notion of price keeping in mind
these wished properties. The convex risk measure provides such a
criterion as we will see in the following.
\subsection{Convex Risk Measures: Definition and Basic Properties
\label{risk measure basic properties}}
Convex risk measures can have two possible interpretations
depending on the representation which is used: they can be
considered either as a pricing rule or as a capital requirement
rule. We will successively present both of them in the following,
introducing each time the vocabulary associated with this
particular approach.
\subsubsection{Risk Measure as an Indifference Price}
We first recall the definition and some key properties of the
convex risk measures introduced by \index[aut]{F\"{o}llmer}
F\"{o}llmer and \index[aut]{Schied} Schied
\cite{FoellmerSchied02a} and \cite{FoellmerSchied02b}. The
notations, definitions and main properties may be found in this
last reference \cite{FoellmerSchied02b}. In particular, we assume
that uncertainty is described through a measurable space $(\Omega,
\Im)$, and that risky positions belong to the linear space of
bounded functions (including constant functions), denoted by
\index[not]{$\mathcal{X}$} $\mathcal{X}$.
\begin{definition}
The functional \index[not]{$\rho$} $\rho :\mathcal{X}\rightarrow
\mathbb{R}$ is a \index[sub]{convex risk measure} {\rm (monetary)
convex risk measure} if, for any $\Phi$ and $\Psi$ in
$\mathcal{X}$, it satisfies
the following properties:\\
$a)$ Convexity: $\forall \lambda \in [ 0,1] \quad \rho\big(
\lambda \Phi+( 1-\lambda ) \Psi\big) \leq \lambda \rho( \Phi) +(
1-\lambda
) \rho( \Psi) $;\\
$b)$ Monotonicity: $\Phi\leq \Psi\Rightarrow \rho( \Phi) \geq \rho
( \Psi)
$;\\
$c)$ Translation invariance: $\forall m\in \mathbb{R}\quad \rho(
\Phi+m) =\rho( \Phi) -m$. \\
A convex risk measure $\rho$ is \index[sub]{coherent risk measure} coherent if it satisfies also:\\
$d)$ Homogeneity : $\forall \lambda \in \mathbb{R}^+ \quad
\rho\big( \lambda\,\Phi\big)=\lambda\rho\big(\Phi\big)$.
\end{definition}
 Note that the convexity property is essential: this translates the
natural fact that diversification should not increase risk. In
particular, any convex combination of ''admissible'' risks should
be ''admissible''. One of the major drawbacks of the famous risk
measure VAR (Value at Risk) is its failure to meet this criterion.
This may lead to arbitrage opportunities inside the financial
institution using it as risk measure as observed by
\index[aut]{Artzner} Artzner, \index[aut]{Delbaen}Delbaen,
\index[aut]{Eber}Eber and \index[aut]{Heath}Heath in their seminal
paper
\cite{ArtznerDelbaenEberHeath}.\\
Intuitively, given the translation invariance, $\rho(X) $ may be
interpreted as the amount the agent has to hold to completely
cancel the risk associated with her risky position $X$ since
\begin{equation}
\label{equation rho-rho} \rho(X+\rho(X))=\rho(X)-\rho(X)=0.
\end{equation}
$\rho(X)$ can be also considered as the opposite of the "buyer's
indifference price" of this position, since when paying the amount
$-\rho(X)$, the new exposure $X-(-\rho(X))$ does not carry any
risk with positive measure, i.e. the agent is somehow indifferent
using this criterion between doing nothing and
having this "hedged" exposure.\\
The convex risk measures appear therefore as a natural extension
of utility functions as they can be seen directly as an
indifference pricing rule.
\subsubsection{Dual Representation \label{dual representation}}
In order to link more closely both notions of pricing rule and
risk measure, the duality between the Banach space $\mathcal{X}$
endowed with the supremum norm $\|.\|$ and its dual space
$\mathcal{X}'$, identified with the set $\mathbf{M}^{\rm ba}$ of
finitely additive set functions with finite total variation on $(
\Omega ,\Im) $, can be used as it leads to a \index[sub]{dual
representation}dual representation. The properties of monotonicity
and cash invariance allow to restrict the domain of the dual
functional to the set $\mathbf{M}_{1,f}$ of all finitely additive
measures (Theorem $4.12$ in
\index[aut]{F\"{o}llmer}\index[aut]{Schied}\cite{FoellmerSchied02b}).
The following theorem gives an "explicit" formula for the risk
measure (and as a consequence for the price) in terms of expected
values:
\begin{theorem}\label{Theorem dual representation}
Let $\mathbf{M}_{1,f}$ be the set of all finitely additive
measures on $(\Omega,\Im)$, and \index[not]{$\alpha$} $\alpha
\big( \mathbf{Q}\big) $ the minimal penalty function
\index[sub]{penalty function} taking values in $\mathbb{R}\cup
\big\{ +\infty \big\} $:
\begin{eqnarray}
& & \forall \mathbf{Q}\in \mathbf{M}_{1,f}\quad \alpha \big( %
\mathbf{Q}\big) =\sup_{\Psi \in \mathcal{X}}\big\{ \mathbb{E}_{\mathbf{Q}%
}[-\Psi] -\rho ( \Psi) \big\} \qquad \big( \geq -\rho( 0) \big).
\label{penalty 2}
\\
& &\text{\bf Dom}(\alpha)=\{\mathbf{Q}\in \mathbf{M}_{1,f}|\>\>\alpha \big( %
\mathbf{Q}\big)<+\infty\}
\end{eqnarray}
The Fenchel duality relation holds :
\begin{eqnarray}
& & \forall \Psi \in \mathcal{X\quad }\rho( \Psi) =\sup_{\mathbf{Q}%
\in \mathbf{M}_{1,f}}\big\{ \mathbb{E}_{\mathbf{Q}}[-\Psi] -\alpha
\big( \mathbf{Q}\big) \big\}   \label{risk measure with penalty}
\end{eqnarray}
Moreover, for any $\Psi \in \mathcal{X}$ there exists an optimal
additive measure $\mathbf{Q}_\Psi \in \mathbf{M}_{1,f}$ such that
\begin{equation*}
\rho( \Psi) =\mathbb{E}_{\mathbf{Q}_\Psi}[-\Psi] -\alpha \big(
\mathbf{Q}_\Psi\big)= \max_{\mathbf{Q} \in \mathbf{M}_{1,f}}\big\{
\mathbb{E}_{\mathbf{Q}}[-\Psi] -\alpha \big( \mathbf{Q}\big)
\big\}.
\end{equation*}
\end{theorem}
Henceforth, $\alpha(\mathbf{Q})$ is the minimal penalty function,
denoted by $\alpha_{\min}(\mathbf{Q})$ in \index[aut]{F\"{o}llmer}\index[aut]{Schied} \cite{FoellmerSchied02b}.\\
The dual representation of $\rho$ given in Equation (\ref{risk
measure with penalty}) emphasizes the interpretation in terms of a
worst case related to the agent's (or regulator's) beliefs.
\vspace{-3mm}
\paragraph{Convex Analysis Point of view}
We start with Remark 4.17 and the Appendices 6 and 7 in
\index[aut]{F\"{o}llmer}\index[aut]{Schied}
\cite{FoellmerSchied02b}. The penalty function $\alpha$ defined in
\eqref{risk measure with penalty} corresponds to the
Fenchel-Legendre transform on the Banach space $\mathcal{X}$ of
the convex risk-measure $\rho$. The dual space $\mathcal{X}'$ can
be identified with the set $\mathbf{M}^{\rm ba}$ of finitely
additive set functions with finite total variation. Then the
subset $\mathbf{M}_{1,f}$ of "finite probability measure" is
weak*-compact in $\mathcal{X}'=\mathbf{M}^{\rm ba}$ and the
functional $\mathbf{Q}\rightarrow \alpha(\mathbf{Q)}$ is
weak*-lower semi-continuous (or weak*-closed) as supremum of
affine functionals. This terminology from convex analysis is based
upon the observation that lower semi-continuity and the closure of
the level sets $\{\phi \leq c\}$ are equivalent properties.
Moreover $\rho$ is lower semi-continuous (lsc) with respect to the
weak topology $\sigma(\mathcal{X},\mathcal{X}')$ since any set
$\{\rho\leq c\}$ is convex and strongly closed given that $\rho$
is strongly Lipschitz-continuous. Then, general duality theorem
for conjugate functional yields to
$$\rho(\Psi)=\sup_{r\in \mathbf{M}^{\rm ba}}(r(\Psi)-\rho^*(r)),
\qquad \rho^*(r)=\sup_{\Psi\in \mathcal{X}}(r(\Psi)-\rho(\Psi))$$
with the convention
$r_{\mathbf{Q}}(\Psi)=\mathbb{E}_{\mathbf{Q}}[-\Psi]$ for
$\mathbf{Q}\in \mathbf{M}_{1,f}$. We then use the properties of
monotonicity and cash invariance of $\rho$ to prove that when
$\rho^*(r)<+\infty$, $-r\in \mathbf{M}_{1,f}$. Moreover by
weak*-compacity of $\mathbf{M}_{1,f}$, the upper semi-continuous
functional $\mathbb{E}_{\mathbf{Q}}[-\Psi]-\alpha(\mathbf{Q})$
attains its maximum on $\mathbf{M}_{1,f}$.\\
In the second part of this chapter, we will intensively used
convex analysis point of view when studying dynamic convex risk
measures. \vspace{-3mm}
\paragraph{Duality and Probability Measures}
We are especially interested in the risk measures that admit a
representation \eqref{risk measure with penalty} in terms of
$\sigma$-additive probability measures $\mathbb{Q}$. In
this paper, for the sake of simplicity and clarity, we use the notation $%
\mathbf{Q}\in \mathbf{M}_{1,f}$ when dealing with additive
measures and $\mathbb{Q}\in \mathcal{M}_1$ when considering
probability measures. So, we are looking for the following
representation on $\cal X$
\begin{equation}
\label{eq:duality proba}
\rho( \Psi) =\sup_{\mathbb{Q}%
\in {\mathcal M}_1}\big\{ \mathbb{E}_{\mathbb{Q}}[-\Psi] -\alpha
\big( \mathbb{Q}\big) \big\}.
\end{equation}
We can no longer expected that the supremum is attained without
additional assumptions. Such representation on $\mathcal{M}_1$ is
closely related to some continuity properties of the convex
functional $\rho$ (Lemma 4.20 and Proposition 4.21 in
\index[aut]{F\"{o}llmer}\index[aut]{Schied}
\cite{FoellmerSchied02b}).
\begin{proposition}
\label{th:duality probability} $i)$ Any convex risk measure $\rho$
defined on $\cal X$ and satisfying \eqref{eq:duality proba} is
continuous from above, in the sense that
$$ \Psi _{n}\searrow \Psi \quad
\Longrightarrow\quad \rho \big( \Psi _{n}\big) \nearrow \rho(
\Psi).$$ $ii)$ The converse is not true in general, but holds
under continuity from below assumption:
 $$ \Psi _{n}\nearrow \Psi \quad
\Longrightarrow \quad \rho \big( \Psi _{n}\big) \searrow \rho(
\Psi).$$ Then any additive measure ${\mathbf Q}$  such that
$\alpha({\mathbf Q})<+\infty$ is $\sigma$-additive and
\eqref{eq:duality proba} holds true. Moreover, from $i)$, $\rho$
is also continuous by above.
\end{proposition}
\subsubsection
{Risk Measures on $\mathbb{L}_{\infty}(\mathbb{P})$}
The representation theory on $\mathbb{L}_{\infty}(\mathbb{P})$ was
developed in particular by \index[aut]{Delbaen}Delbaen
\cite{Delbaen00} and extended by \index[aut]{Frittelli}Frittelli
and \index[aut]{Gianin}Gianin \cite{FrittelliGianin02} and
\cite{FrittelliGianin04} and
\index[aut]{F\"{o}llmer}\index[aut]{Schied}
\cite{FoellmerSchied02b}. When a probability measure $\mathbb{P}$
is given, it is natural indeed to define risk measures $\rho$ on
$\mathbb{L}_{\infty}(\mathbb{P})$ instead of on $\cal X$
satisfying the compatibility condition:
\begin{equation}
\label{eq:riskmeasure equalite} \rho(\Psi)=\rho(\Phi)\quad
\text{if}\quad\Psi=\Phi\quad \mathbb{P}-a.s.
\end{equation}
Let us introduce some new notations:$ \>
\mathbf{M}_{1,ac}(\mathbb{P})$ is the set of finitely additive
measures absolutely continuous w. r. to $\mathbb{P}$ and $
\mathcal{M}_{1,ac}(\mathbb{P})$ is the set of probability measures
absolutely continuous w. r. to $\mathbb{P}$.\\ We also define
natural extension of continuity from below in the space
$\mathbb{L}_{\infty}(\mathbb{P})$: ($\> \Psi _{n}\searrow \Psi
\quad \mathbb{P}-a.s.\Longrightarrow\quad \rho \big( \Psi
_{n}\big) \nearrow \rho( \Psi)$), or continuity from above in the
space $\mathbb{L}_{\infty}(\mathbb{P})$: ($\> \Psi _{n}\nearrow
\Psi \quad \mathbb{P}-a.s.\Longrightarrow\quad \rho \big( \Psi
_{n}\big) \searrow \rho( \Psi)$).\\ These additional results on
conjugacy relations are given
in\index[aut]{F\"{o}llmer}\index[aut]{Schied}
\cite{FoellmerSchied02b} Theorem 4.31 and in \index[aut]{Delbaen}
Delbaen \cite{Delbaen00} Corollary 4.35). Sometimes, as in
\index[aut]{Delbaen} \cite{Delbaen00}, the continuity from above
is called the Fatou property.
\begin{theorem} \label{Theorem rep duale Mac}
Let $\mathbb{P}$ be a given probability measure.\vspace{-1mm}
\begin{enumerate}
\item Any convex risk measure $\rho$ on $\mathcal{X}$ satisfying \eqref{eq:riskmeasure
equalite} may be considered as a risk measure on $\mathbb{L}_{\infty}(\mathbb{P})$. A dual
representation holds true in terms of absolutely continuous additive measures $\mathbf{Q}\in
\mathbf{M}_{1,ac}(\mathbb{P})$.\vspace{-2mm}
\item $\rho$ admits a dual representation on $
\mathcal{M}_{1,ac}(\mathbb{P})$:
$$\alpha(\mathbb{Q})= \sup_{\Psi \in
\mathbb{L}_{\infty}(\mathbb{P})} \big\{\mathbb{E}_{\mathbb{Q}}[-\Psi]-\rho(\Psi) \big\},
\qquad \rho( \Psi) =\sup_{\mathbb{Q}%
\in \mathcal{M}_{1,ac}(\mathbb{P})}\big\{ \mathbb{E}_{\mathbb{Q}}[-\Psi] -\alpha \big(
\mathbb{Q}\big) \big\}$$
if and only if one of the equivalent properties holds:\\
$a) \quad$ $\rho$ is continuous from above (Fatou property);\\
$b) \quad$ $\rho$ is closed for the weak*-topology
$\sigma(\mathbb{L}_{\infty},\mathbb{L}_1)$;\\
$c) \quad$ the acceptance set $\{\rho\leq 0\}$ is weak*-closed in
$\mathbb{L}_{\infty}(\mathbb{P})$.\vspace{-2mm}
\item Assume that $\rho$ is a coherent
(homogeneous) risk measure, satisfying the Fatou property. Then,
\begin{equation}
\label{eq:coherent duality} \rho(\Psi)=\sup_{\mathbb{Q}\in
\mathcal{M}_{1,ac}(\mathbb{P})}\big\{\mathbb{E}_{\mathbb{Q}}[-\Psi]\>|\>\alpha(\mathbb{Q})\>=0\big\}
\end{equation}
The supremum in \eqref{eq:coherent duality} is a {\bf maximum} iff
one of the following equivalent properties holds:\\
$a) \quad$ $\rho$ is continuous from below;\\
$b) \quad$ the convex set $\mathcal{Q}=\{\mathbb{Q}\in
\mathcal{M}_{1,ac}\big|\>\alpha(\mathbb{Q})\>=\>0\}$ is weakly
compact in $\mathbb{L}^1(\mathbb{P})$.
\end{enumerate}
\end{theorem}
According to the Dunford-Pettis theorem, the weakly relatively
compact sets of $\mathbb{L}^1(\mathbb{P})$ are sets of uniformly
integrable variables and La Vall\'{e}e-Poussin gives a criterion
to check this property. Therefore, the subset $\mathcal{A}$ of
$\mathbb{L}^1(\mathbb{P})$ is weakly relatively compact iff it is
closed and uniformly integrable. Moreover, according to the La
Vall\'{e}e-Poussin criterion, an increasing convex continuous
function $\Phi : \mathbb{R}_{+}\rightarrow \mathbb{R}$, (also
called Young's function) such that:
$$\lim_{x\rightarrow \infty}  \frac{\Phi(x)}{x}=+\infty \qquad
{\rm and}  \qquad \sup_{\mathbb{Q}\in \mathcal{A}}
\mathbb{E}_{\mathbb{P}}\left[\Phi(\frac{d\mathbb{Q}}{d\mathbb{P}})\right]<+\infty.$$
\subsection{Comments on Measures of Risk and Examples}
\subsubsection{About Value at Risk}
Risk measures, just as utility functions, go beyond the simple
problem of pricing. Both are inherently a choice or decision {\it
criterion}. More precisely, when assessing the risk related to a
given position in order to define the amount of capital
requirement, a first natural approach is based on the distribution
of the risky position itself. In this framework, the most
classical measure of risk is simply the \emph{variance} (or the
mean-variance analysis). However, it does not take into account
the whole distribution's features (as asymmetry or skewness) and
especially it does not focus on the ''real'' financial risk which
is the downside risk. Therefore different methods have been
developed to focus on the risk of losses: the most widely used (as
it is recommended to bankers by many
financial institutions) is the so-called \emph{Value at Risk} \index[sub]{value at risk} (denoted by $%
VAR$), based on quantiles of the lower tail of the distribution.
More precisely, the $VAR$  associated with the position $X$ at a
level $\varepsilon $ is defined as
\begin{equation*}
VAR_{\varepsilon }\big( X\big) =\inf \big\{ k:\mathbb{P}( X+k<0)
\leq \varepsilon \big\}.
\end{equation*}
The $VAR$ corresponds to the minimal amount to be added to a given
position to make it acceptable. Such a criterion satisfies the key
properties of decreasing monotonicity, translation invariance
since $\forall m\in \mathbb{R}$, $VAR_{\varepsilon }\big(
X+m\big) =VAR_{\varepsilon }\big( X\big) -m$ and finally,
the $VAR$ is positive homogeneous as $\forall \lambda \geq 0$, $%
VAR_{\varepsilon }\big( \lambda X\big) =\lambda VAR_{\varepsilon
}\big(
X\big) $.\\
This last property reflects the linear impact of the size of the
position on the risk measure. However, as noticed by Artzner et
al.\index[aut]{Artzner}\index[aut]{Delbaen}\index[aut]{Eber}\index[aut]{Heath}
\cite{ArtznerDelbaenEberHeath} this criterion fails to meet a
natural consistency requirement: it is not a convex risk measure
while the convexity property translates the natural fact that
diversification should not increase risk. In particular, any
convex combination of ''admissible'' risks should be
''admissible''. The absence of convexity of the $VAR$ may lead to
arbitrage opportunities inside the financial institution using
such criterion as risk measure. Based on this logic, Artzner et
al.\index[aut]{Artzner}\index[aut]{Delbaen}\index[aut]{Eber}\index[aut]{Heath}
\cite{ArtznerDelbaenEberHeath} have adopted a more general
approach to risk measurement. Their paper is essential as it has
initiated a systematic axiomatic approach to risk measurement. A
{\it coherent measure of risk} should be convex and satisfy the
three key properties of the $VAR$ \vspace{-9mm}
\paragraph{Conditional Value at Risk}
 For instance, a coherent version of the Value at Risk is the
so-called {\it Conditional Value at Risk} \index[sub]{conditional
value at risk} as observed by \index[aut]{Rockafellar} Rockafellar
and \index[aut]{Uryasev} Uryasev \cite{RockafellarUryasev00}. This
risk measure is denoted by $CVAR_{\varepsilon}$ and defined as
\begin{equation*}
CVAR_{\lambda}(X)
=\inf_K\mathbb{E}\Big[\frac{1}{\lambda}(X-K)^{-}-K\Big].
\end{equation*}
This coincide with the {\it Expected Shortfall}
\index[sub]{expected shortfall} under some assumptions for the
$X$-distribution (for more details, see Corollary 5.3 in
\index[aut]{Acerbi} Acerbi and \index[aut]{Tasche} Tasche
\cite{AcerbiTasche}). In this case, the $CVAR$ can be written as
\begin{equation*}
CVAR_{\lambda}(X) =\mathbb{E}[-X|X+VAR_{\lambda}(X)<0].
\end{equation*}
Moreover, the $CVAR$ also coincides with another coherent version
of the $VAR$, called {\it Average Value at Risk}
\index[sub]{average value at risk} and denoted by $AVAR$. This
risk measure is defined as:
$$AVAR_{\lambda}(\Psi)=\frac{1}{\lambda} \int_0^{\lambda} VAR_{\epsilon}(\Psi) d{\epsilon}.$$
For more details, please refer for instance to
\index[aut]{F\"{o}llmer} F\"{o}llmer and
\index[aut]{Schied} Schied \cite{FoellmerSchied02b} (Proposition 4.37).\vspace{2mm}\\
More recently, the axiom of positive homogeneity has been
questioned. Indeed, such a condition does not seem to be
compatible with the notion of liquidity risk existing on the
market as it implies that the size of the risky position has
simply a linear impact on the risk measure. To tackle this
shortcoming, \index[aut]{F\"{o}llmer}F\"{o}llmer and
\index[aut]{Schied}Schied consider, in \cite{FoellmerSchied02a}
and \cite{FoellmerSchied02b}, {\it convex risk measures} as
previously defined.
\subsubsection{Risk Measures and Utility Functions}
\paragraph{Entropic Risk Measure}
The most famous convex risk measure on
$\mathbb{L}_{\infty}(\mathbb{P})$ is certainly the entropic risk
measure \index[sub]{entropic risk measure} defined as the
functional $e_{\gamma}$ in the previous section when considering
an exponential utility framework. The dual formulation of this
continuous from below functional justifies the name of {\it
entropic risk measure} since:
\begin{equation*}
\forall \Psi \in \mathbb{L}_{\infty}(\mathbb{P})\quad \quad
e_{\gamma }( \Psi ) =\gamma \ln \mathbb{E}_{\mathbb{P}}\Big[ \exp
\big( -\frac{1}{\gamma }\Psi
\big) \Big] =\sup_{\mathbb{Q}\in \mathcal{M}_{1}}\big\{ \mathbb{E}_{%
\mathbb{Q}}[-\Psi] -\gamma h\big( \mathbb{Q}|\mathbb{P}\big)
\big\}
\end{equation*}
where $h( \mathbb{Q}|\mathbb{P}) $ is the relative entropy of $%
\mathbb{Q}$ with respect to the prior probability measure
$\mathbb{P}$, defined by
\begin{equation*}
h\big( \mathbb{Q}|\mathbb{P}\big) =\mathbb{E}_{\mathbb{P}}\Big[\frac{d\mathbb{Q}}{d\mathbb{P}}\ln \frac{d%
\mathbb{Q}}{d\mathbb{P}}\Big] \quad \text{if }\mathbb{Q\ll P}
\qquad \mbox{and} \quad +\infty\quad \ \text{otherwise}.
\end{equation*}
Since $e_{\gamma}$ is continuous from below in
$\mathbb{L}_{\infty}(\mathbb{P})$, by the previous theorem
\begin{equation*}
\forall \Psi \in \mathbb{L}_{\infty}(\mathbb{P})\quad \quad
e_{\gamma }( \Psi ) =\gamma \ln \mathbb{E}_{\mathbb{P}}\Big[ \exp
\big( -\frac{1}{\gamma }\Psi
\big) \Big] =\mathbf{ max}_{\mathbb{Q}\in \mathcal{M}_{1,ac}}\big\{\mathbb{E}_{%
\mathbb{Q}}[-\Psi] -\gamma h\big( \mathbb{Q}|\mathbb{P}\big)
\big\}
\end{equation*}
As previously mentioned in Paragraph \ref{sssect:indifference
pricing rule}, this particular convex risk measure is closely
related to the exponential utility function and to the associated
indifference price. However, the relationships between risk
measures and utility functions can be extended. \vspace{-3mm}
\paragraph{Risk Measures and Utility Functions}
More generally, risk measures and utility functions have close
relationships based on the hedging and super-replication problem.
It is however
possible to obtain a more general connection between them using the notion of shortfall risk.\\
More precisely, any agent having a utility function $U$ assesses
her risk by taking the expected utility of the considered position
$\Psi \in \mathbb{L}_{\infty}(\mathbb{P})$:
$\mathbb{E}_{\mathbb{P}}[U(\Psi)]$. If she focuses on her "real"
risk, which is the downside risk, it is natural to consider
instead the {\it loss function} \index[sub]{loss function}
$\mathcal{L}$ defined by $\mathcal{L}(x)=-U(-x)$
(\index[aut]{F\"{o}llmer}\index[aut]{Schied}\cite{FoellmerSchied02b}
Section 4.9). As a consequence, $\mathcal{L}$ is a convex and
increasing function and maximizing the expected utility is
equivalent to minimize the expected loss (also called
the {\it shortfall risk} \index[sub]{shortfall risk}), $\mathbb{E}_{\mathbb{P}}[\mathcal{L}(-\Psi)]$.\\
It is then natural to introduce the following risk measure as the
opposite of the indifference price:
$$\rho(X)=\inf\{m \in
\mathbb{R}\big|\>\mathbb{E}_{\mathbb{P}}[\mathcal{L}(-\Psi-m)]
\leq l(0)\}.$$ Moreover, there is an explicit formula for the
associated penalty function given in terms of the
\index[sub]{Fenchel-Legendre transform} Fenchel-Legendre transform
$\mathcal{L}^{*}(y)=\sup\{-xy-l\mathcal{L}(x)\}$ of the convex
function $\mathcal{L}$ (\cite{FoellmerSchied02b} Theorem 4.106):
$$\alpha(\mathbb{Q})=\inf_{\lambda
>0}\Big\{\frac{1}{\lambda}\Big(\mathcal{L}(0)+\mathbb{E}_{\mathbb{P}}\big[\mathcal{L}^{*}\big(\lambda
\frac{d\mathbb{Q}}{d\mathbb{P}}\big)\big]\Big)\Big\}.$$
\subsection{Risk Measures and Hedging \label{section risk measures and hedging}}
In this subsection, we come back to the possible interpretation of
the risk measure $\rho(X)$ in terms of capital
requirement\index[sub]{capital requirement}. This leads also to a
natural relationship between risk measure and hedging. We then
extend it to a wider perspective of super-hedging.
\subsubsection{Risk Measure and Capital Requirement}
Looking back at Equation (\ref{equation rho-rho}), the risk
measure $\rho(X)$ gives an assessment of the minimal {\it capital
requirement} to be added to the position as to make it acceptable
in the sense that the new position ($X$ and the added capital)
does not carry any risk with non-negative measure any more. More
formally, it is natural to introduce the \emph{acceptance set}
\index[sub]{acceptance set}$\mathcal{A}_{\rho }$ related to
$\rho$\index[not]{$\mathcal{A}_{\rho}$} defined as the set of all
acceptable positions in the sense that they do not require any
additional capital:
\begin{equation}
\mathcal{A}_{\rho }=\big\{ \Psi \in \mathcal{X},\quad \rho( \Psi)
\leq 0\big\}.  \label{acceptance set}
\end{equation}
Given that the epigraph of the convex risk measure $\rho$ is ${\rm
epi}(\rho)=\{(\Psi,m)\in {\mathcal{X}}\times
\mathbb{R}\big|\>\rho(\Psi)\leq m\}= \{(\Psi,m)\in
{\mathcal{X}}\times \mathbb{R}\big|\>\rho(\Psi+m)\leq 0\}$, the
characterization of $\rho $ in terms of $\mathcal{A}_{\rho }$ is
easily obtained
\begin{equation*}
\rho(X) =\inf \big\{ m\in \mathbb{R};m+X \in \mathcal{A}%
_{\rho }\big\}. \label{rho and acceptance set}
\end{equation*}
This last formulation makes very clear the link between risk
measure and capital requirement.\\
From the definition of both the convex risk measure $\rho $ and the acceptance set
$\mathcal{A}_{\rho }$ and the dual representation of the risk measure $\rho$, it is possible
to obtain another characterization of the associated penalty function $\alpha$ as:
\begin{equation}
\alpha
\big( \mathbf{Q}\big) =\sup_{\Psi \in \mathcal{A}_{\rho }}\mathbb{E}_{%
\mathbf{Q}}[-\Psi],\quad\text{if}\quad \mathbf{Q}\mathbb{\in }\mathbf{M}_{1,f},\qquad
=+\infty,\quad \text{if not}. \label{penalty 1 (acceptance set)}
\end{equation}
$\alpha(\mathbf{Q})$ is the support function of $-\mathcal{A}_{\rho}$, denoted by
$\Sigma^{\mathcal{A}_{\rho}}(\mathbf{Q})$. When $\mathcal{A}_{\rho} $ is a cone, i.e. $\rho$
is a coherent (positive homogeneous) risk measure, then $\alpha(\mathbf{Q})$
only takes the values $0$ and $+\infty$.\\
By definition, the set $\mathcal{A}_{\rho}$ is "too large" in the
following sense: even if we can write $m+X \in \mathcal{A}_{\rho}$
as $m+X = \xi \in \mathcal{A}_{\rho}$, we cannot have an explicit
formulation for $\xi$ and in particular cannot compare $m+X$ with
$0$. Therefore, it seems natural to consider a (convex) class of
variables $\mathcal{H}$ such that $m+X \geq H \in \mathcal{H}$.
$\mathcal{H}$ appears as a natural (convex) set from which a risk
measure can be generated.
\subsubsection{Risk Measures Generated by a Convex Set\label{risk measure generated by convex set}}
\paragraph{Risk Measures Generated by a Convex in $\mathcal{X}$}
In this section, we study the generation of a convex risk measure
from a general convex set.
\begin{definition}
\label{definition risk measure generated by convex set} Given a
non-empty convex subset $\mathcal{H}$ of $\mathcal{X}$ such that
$\inf\{m\in \mathbb{R}\>\big|\>\>\exists \xi \in
\mathcal{H},\text{m }\geq \xi\}>-\infty$, the functional
\index[not]{$\nu^{\mathcal{H}}$} $\nu ^{\mathcal{H}}$ on
$\mathcal{X}$
\begin{equation}
\label{eq: risk measure generated by a set} \nu ^{\mathcal{H}}(
\Psi) =\inf \big\{ m\in \mathbb{R};\>\>\exists \xi \in
\mathcal{H},\text{m }+\Psi \geq \xi \big\}
\end{equation}%
is a convex risk measure. Its minimal penalty function $\alpha^{\mathcal{H}}$ is given by:
$\alpha^{\mathcal{H}}(\mathbf{Q})=\sup_{H \in \mathcal{H}} \mathbb{E}_\mathbf{Q}[-H]$.
\end{definition}
The main properties of this risk measure are listed or proved below:
\begin{enumerate}
\item The acceptance set of $\nu ^{\mathcal{H}}$ contains the
convex  subsets $ \mathcal{H}$ and
$\mathcal{A}_{\mathcal{H}}=\big\{ \Psi \in \mathcal{X},\exists \xi
\in \mathcal{H},\>\> \Psi \geq \xi \big\}$. Moreover,
$\mathcal{A}_{\nu ^{\mathcal{H}}}=\mathcal{A}_{\mathcal{H}}$ if
the last subset is closed in the following sense: For $\xi
\in\mathcal{A}_{\mathcal{H}}$ and $\Psi \in \mathcal{X}$, the set
$\{\lambda\in[0,1]>\big|\>\lambda \xi+(1-\lambda)\Psi\in
\mathcal{A}_{\mathcal{H}}\}$ is closed in $[0,1]$ (see Proposition
4.6 in
\index[aut]{F\"{o}llmer}\index[aut]{Schied}\cite{FoellmerSchied02b}).
\item The penalty function $\alpha^{\mathcal{H}}$ associated with
$\nu ^{\mathcal{H}}$ is the support function of
$-\mathcal{A}_{\nu^{\mathcal{H}}}$ defined by
$\alpha^{\mathcal{H}}\big( \mathbf{Q}\big)
=\Sigma^{\mathcal{A}_{\nu^{\mathcal{H}}}}\big( \mathbf{Q}\big) =
\sup_{X\in \mathcal{A}_{\nu^{\mathcal{H}}}}\mathbb{E}_{\mathbf{Q}%
}[-X]$. Let us show that $\alpha^{\mathcal{H}}$ is also nothing
else but $\Sigma^{\mathcal{H}}$:\\
For any $X\in\mathcal{A}_{\nu^{\mathcal{H}}}$ there exist
$\epsilon>0$ and $\xi \in \mathcal{H}$ such that $-X \leq
-\xi+\epsilon.$ Taking the "expectation" with respect to the
additive measure $\mathbf{Q}\in \mathbf{M}_{1,f}$, it follows that
$\mathbb{E}_{\mathbf{Q}}[-X]\leq
\mathbb{E}_{\mathbf{Q}}[-\xi]+\varepsilon\leq
\Sigma^{\mathcal{H}}\big( \mathbf{Q}\big)+\varepsilon$ where
$\Sigma^{\mathcal{H}}\big( \mathbf{Q}\big) =
\sup_{H\in \mathcal{H}}\mathbb{E}_{\mathbf{Q}%
}[-H]$. Taking the supremum with respect to
$X\in\mathcal{A}_{\nu^{\mathcal{H}}}$ on the left hand side, we
deduce that $\Sigma^{\mathcal{A}_{\nu^{\mathcal{H}}}}\leq
\Sigma^{\mathcal{H}}$; the desired result follows from the
observation that $\mathcal{H}$ is included in
$\mathcal{A}_{\nu^{\mathcal{H}}}$. \item When $\mathcal{H}$ is a
cone, the corresponding risk measure is coherent (homogeneous).
The penalty function $\alpha ^{\mathcal{H}}$ is the indicator
function (in the sense of the convex analysis) of the orthogonal
cone $\mathbf{M}_{\mathcal{H}}$:
\index[not]{$l^{\mathbf{M}_{\mathcal{H}}}$}
$l^{\mathbf{M}_{\mathcal{H}}}\big( \mathbf{Q}\big) =0 \quad
\text{if }\quad \mathbf{Q}\mathbb{\in
}\mathcal{M}_{\mathcal{H}}\>,\>\> +\infty \> \text{otherwise}$,
where
\begin{equation*}
\mathbf{M}_{\mathcal{H}}=\big\{ \mathbf{Q\in
}\mathbf{M}_{1,f};\forall \xi \in \mathcal{H},\
\mathbb{E}_{\mathbf{Q}}[-\xi] \leq 0\big\}.
\end{equation*}
The dual formulation of $\nu^{\mathcal{H}}$ is simply given  for
$\Psi \in \mathcal{X}$ by: $\quad \nu ^{\mathcal{H}}( \Psi )
=\sup_{\mathbf{Q}\in
\mathbf{M}_{\mathcal{H}}}\mathbb{E}_{\mathbf{Q}}[-\Psi].$
\end{enumerate}
It is natural to associate the convex indicator $l^{\mathcal{H}}$
on $\mathcal{X}$ with the set $\mathcal{H}$, $l^{\mathcal{H}}(X)
=0 \> \text{if }X \in \mathcal{H} \>;\>\> +\infty \>
\text{otherwise}$. This convex functional is not translation
invariant, and therefore it is not a convex risk measure.
Nevertheless, $l^{\mathcal{H}}$ and $\nu^{\mathcal{H}}$ are
closely related as follows:
\begin{corollary}\label{re:worst risk measure and hedging}
Let $l^{\mathcal{H}}$ be the convex indicator on $\mathcal{X}$ of
the convex set $\mathcal{H}$.\\
The risk measure $\nu^{\mathcal{H}}$, defined in Equation
(\ref{eq: risk measure generated by a set}), is the largest convex
risk measure dominated by $l^{\mathcal{H}}$ and it can be
expressed as:
$$\nu
^{\mathcal{H}}(\Psi)=\inf_{\xi \in \mathcal{X}}\{\rho_{\rm
worst}(\Psi-\xi)+ l^{\mathcal{H}}(\xi)\}$$ where $\rho_{\rm
worst}(\Psi)=\sup_{\omega\in\Omega}\{-\Psi(\omega)\}$ is the worst
case risk measure.
\end{corollary}
{\bf Proof:} Let $\mathcal{L}=\{m \in \mathbb{R},\>\> \exists \xi
\in \mathcal{H},\>\> m \geq \xi \}$. This set is a half-line with
lower bound $\inf_{\xi \in \mathcal{H}}\sup_{\omega}\xi(\omega)$.\\
Moreover, for any $m_0 \notin \mathcal{L}$, $m_0 \leq
\inf_{\xi \in \mathcal{H}}\sup_{\omega}\xi(\omega)$. Therefore,
$\nu^{\mathcal{H}}(0) =
\inf_{\xi \in \mathcal{H}}\sup_{\omega}\xi(\omega)=\inf_{\xi}\rho_{\rm
worst}(-\xi)$.\\ The same arguments hold for
$\nu^{\mathcal{H}}(\Psi)$.
$\square$\vspace{2mm}\\
Therefore, $\nu ^{\mathcal{H}}$ may be interpreted as the worst
case risk measure $\rho_{\rm worst}$ reduced by the use of
(hedging) variables in $\mathcal{H}$. This point of view would be
generalized in Corollary \ref{corollary modification de rho par H}
in terms of the inf-convolution $\nu^{\mathcal{H}}=\rho_{\rm
worst}\square l^{\mathcal{H}}.$\vspace{-3mm}
\paragraph{Risk Measures Generated by a Convex Set in $\mathbb{L}_{\infty}(\mathbb{P})$}
\label{par:riskmeasureconvexset}
Assume now $\mathcal{H}$ to be a convex subset of
$\mathbb{L}_{\infty}(\mathbb{P})$. The functional
{$\nu^{\mathcal{H}}$} on $\mathbb{L}_{\infty}(\mathbb{P})$ is
still defined by the same formula \eqref{eq: risk measure
generated by a set}, in which the inequality has to be understood
in $\mathbb{L}_{\infty}(\mathbb{P})$, i.e. $\mathbb{P}-a.s.$, with
a penalty function only defined on $\mathbf{M}_{1,ac}(\mathbb{P})$
and given by $\alpha^{\mathcal{H}}(\mathbf{Q})=\sup_{H \in
\mathcal{H}}\mathbb{E}_{\mathbf{Q}}[-H]$.\\
The problem is then to give condition(s) on the set $\mathcal{H}$
to ensure that the dual representation holds on
$\mathcal{M}_{1,ac}(\mathbb{P})$ and not only on
$\mathbf{M}_{1,ac}(\mathbb{P})$. By Theorem \ref{Theorem rep duale
Mac}, this problem is equivalent to the continuity from above of
the risk measure $\nu^{\mathcal{H}}$ or equivalently to the
weak*-closure of its acceptance set $\mathcal{A}_{\mathcal{H}}$.
Properties of this kind are difficult to check and in the
following, we will simply give some examples where this property
holds.
\subsection{Static Hedging and Calibration \label{subsection hedging description}}
In this subsection, we consider some examples motivated by
financial risk hedging
problems.\vspace{-1mm}
\subsubsection{Hedging with a Family of Cash Flows}
We start with a very simple model where it is only possible to
hedge statically over a given period using a finite family of {\em
bounded} cash flows $\{C_1,C_2,...,C_d\}$, the (forward) price of
which is known at time $0$ and denoted by
$\{\pi_1,\pi_2,...,\pi_d\}$. All cash flows are assumed to be
non-negative and non-redundant. Constants may be included and then
considered as assets.\\
We assume that the different prices are {\bf coherent} in the sense that
$$\exists \>\mathbb{Q}_0  \sim\mathbb{P},s.t.\quad \forall i,\>\>
\mathbb{E}_{\mathbb{Q}_0}[C_i]=\pi_i$$ Such an assumption implies
in particular that any inequality on the cash flows is preserved
on the prices. The quantities of interest are often the gain
values of
the basic strategies, $G_i=C_i-\pi_i$.\\
We can naturally introduce the non-empty set $\mathcal{Q}_e$ of equivalent
{\em "martingale measures"} as
$$\mathcal{Q}_e=\{\>\mathbb{Q}|\>\>\mathbb{Q}  \sim\mathbb{P},s.t.\quad \forall i,\>\>
\mathbb{E}_\mathbb{Q}[G_i]=0\>\}$$
The different instruments we consider are very liquid; by selling
or buying some quantities $\theta_i$ of such instruments, we
define the family $\mathbf{\Theta}$ of gains associated with
trading strategies $\theta$:
$$\mathbf{\Theta}=\Big\{\>G(\theta)=\sum_{i=1}^d\theta_i\,G_i,\>\theta\in \mathbb{R}^d,\quad \text{with initial
value}\quad\sum_{i=1}^d\theta_i\,\pi_i\Big\}$$ This framework is
very similar to Chapter 1 in \index[aut]{F\"{o}llmer} F\"{o}llmer
and \index[aut]{Schied} Schied \cite{FoellmerSchied02b} where it
is shown that the assumption of coherent prices is equivalent to
the absence of arbitrage opportunity in the market defined as
$$\text{\bf(AAO)}\qquad G(\theta)\geq 0 \quad \mathbb{P}\>\>a.s. \quad
\Rightarrow G(\theta) = 0 \quad \mathbb{P}\>\>a.s.$$ These
strategies can be used to hedge a risky position $Y$. In the
classical financial literature, a superhedging strategy is a par
$(m,\theta)$ such that $m+ G(\theta)\geq Y,\>a.s.$ This leads to
the notion of superhedging (super-seller) price $\pi^{\rm
sell}_{\uparrow}(Y)=\inf\{\>m \>|\>\exists \>G(\theta)\>
s.t.\>\>m+ G(\theta)\geq Y\}.$ In terms of risk measure, we are
concerned with the static superhedging price of $-Y$. So, by
setting $\mathcal{H}=-\mathbf{\Theta}$, we define the risk measure
$\nu^{\mathcal{H}}$ as
$$\nu^{\mathcal{H}}(X)=\pi^{\rm
sell}_{\uparrow}(-X)=\inf\{m \in \mathbb{R},\>\> \exists \theta
\in \mathbb{R}^d: m+X +G(\theta)\geq 0\}$$ Let us observe that the
no arbitrage assumption implies that
$\mathbb{E}_{\mathbb{Q}_0}[G(\theta)]=0$. Hence,
$\nu^{\mathcal{H}}(0)\geq
\mathbb{E}_{\mathbb{Q}_0}[-X]>-\infty$.
Moreover, the dual representation of the risk measure
$\nu^{\mathcal{H}}$ in terms of probability measures is closely
related to the absence of arbitrage opportunity as underlined in
the following proposition (Chapter 4 in \cite{FoellmerSchied02b}):
\begin{proposition}
$i)$ If the market is arbitrage-free, i.e. {\bf (AAO)} holds true,
the convex risk measure $\nu^{\mathcal{H}}$ can be
represented in terms of the set of equivalent "martingale"
measures $\mathcal{Q}_e$ as
\begin{equation}
\label{eq:superprice} \nu^{\mathcal{H}}(\Psi)=\sup_{\mathbb{Q}\in
\mathcal{Q}_e}\mathbb{E}_\mathbb{Q}(-X),
\quad\text{where}\quad\mathcal{Q}_e=\{\mathbb{Q} \sim
\mathbb{P},\>\mathbb{E}_{\mathbb{Q}}(G_i)=0, \>\forall i=1...d\}.
\end{equation}
By Theorem \ref{Theorem dual representation}, this
$\mathbb{L}_{\infty}(\mathbb{P})$-risk measure is continuous from above.\\
$ii)$ Moreover, the market is arbitrage-free if
$\nu^{\mathcal{H}}$ is {\bf sensitive} in the sense that
$\nu^{\mathcal{H}}(\Psi) > \nu^{\mathcal{H}}(0)$ for all $\Psi$
such that $\mathbb{P}(X<0)>0$ and $\mathbb{P}(X \leq 0)=1$.
\end{proposition}
\subsubsection{Calibration Point of View and Bid-Ask Constraint}
This point of view is often used on financial markets when cash flows depend on some basic
assets $(S_1,S_2,...,S_n)$, whose
characteristics will be given in the next paragraph.\\
We can consider for instance $(C_i)$ as payoffs of derivative
instruments, sufficiently liquid to be used as {\it calibration
tools and static hedging strategies}. So far, all agents having
access to the market agree on the derivative prices, and do not
have any
restriction on the quantity they can buy or sell.\\
We now take into account some restrictions on the trading. We
first introduce a bid-ask spread on the (forward) price of the
different cash flows. We denote by $\pi^{ask}_i(C_i)$ the market
buying price and by $\pi^{bid}_i(C_i)$ the market selling price.
The price coherence is now written as
$$\>\exists \>\mathbb{Q}_0 \sim \mathbb{P},\quad \forall i,\>\>
\pi^{ask}_i(C_i) \leq \mathbb{E}_{\mathbb{Q}_0}[C_i] \leq
\pi^{bid}_i(C_i)$$ To define the gains family, we need to make a
distinction between cash flows when buying and cash flows when
selling. To do that, we double the number of basic gains, by
associating, with any given cash-flow $C_i$, both gains
$G^{bid}_i=C_i-\pi^{bid}_i$ and $G^{ask}_i=\pi^{ask}_i(C_i)-C_i$.
Henceforth, we do not make distinction of the notation and we
still denote any gain by $G_i$. The price coherence is then
expressed as
$$\>\exists \>\mathbb{Q}_0 \sim \mathbb{P},\quad \forall\> i=1....\mathbf{2d},\>\>
\mathbb{E}_{\mathbb{Q}_0}[G_i] \leq 0.$$ The set of such probability
measures, called {\em super-martingale measures}, is denoted by
$\mathcal{Q}_e^s$. Note that the coherence of the prices implies
that the set $\mathcal{Q}_e^s$ is non empty.\\
Using this convention, a strategy is defined by a $2d$-dimensional
vector $\theta$, the components of which are all non-negative.
More generally, we can introduce more trading restriction on the
size of the transaction by constraining $\theta$ to belong to a
convex set $\mathcal{K}\subseteq \mathbb{R}^{2d}_+$ such that
$0\in \mathcal{K}$. Note that we can also take into account some
limits to the resources of the investor, in such way the initial
price $\langle\theta,\pi\rangle$ has an an upper bound. In any
case, we still denote the set of admissible strategies by
$\mathcal{K}$ and the family of associated gains by:
$\mathbf{\Theta}=\Big\{\>G(\theta)=\sum_{i=1}^{2d}\theta_iG_i,\>\theta\in
\mathcal{K}\Big\}$.\\
In this constrained framework, the relationship between price
coherence and ({\bf AAO}) on
$\mathbf{\Theta}$ has been studied in details in \index[aut]{Bion-Nadal}Bion-Nadal
\cite{Bion-Nadal2} but also in Chapter 1 of \cite{FoellmerSchied02b}.\\
More precisely, as above, the price coherence implies that the risk measure
$\nu^{\mathcal{H}}$ related to $\mathcal{H}=-\mathbf{\Theta}$ is
not identically
$-\infty$.\\
A natural question is to extend the duality relationship
\eqref{eq:superprice} using the subset of {\em super-martingale
measures}. Using Paragraph \ref{par:riskmeasureconvexset}, this
question is equivalent to show that the minimal penalty function
is infinite outside of the set of absolutely continuous
probability measures and that
$\nu^{\mathcal{H}}$ is continuous from above.\\
When studying the risk measure $\nu^{\mathcal{H}}$ (Definition \ref{definition risk measure generated by convex set}
and its properties), we have proved that:
$$\forall \mathbb{Q} \in \mathcal{M}_{1,ac}(\mathbb{P}), \quad
\alpha^{\mathcal{H}}(\mathbb{Q})=\sup_{\xi \in
\mathcal{H}}\mathbb{E}_{\mathbb{Q}}[-\xi] =\sup_{\theta \in
\mathcal{K}}\mathbb{E}_{\mathbb{Q}}[G(\theta)].$$ In particular,
since $0 \in \mathcal{K}$, if $\mathbb{Q}\in \mathcal{Q}_e^s$,
then $\alpha(\mathbb{Q})=0$. Moreover, if $\Theta$ is a cone, then
$\alpha^{\mathcal{H}}$ is the indicator function of
$\mathcal{Q}^s_e$.\vspace{1mm}\\
It remains to study the continuity from above of
$\nu^{\mathcal{H}}$ and especially to relate it with the absence
of arbitrage opportunity in the market. We summarize below the
results \index[aut]{F\"{o}llmer}F\"{o}llmer and
\index[aut]{Schied}Schied obtained in Theorem 4.95 and Corollary
9.30 \cite{FoellmerSchied02b}.
\begin{proposition}
Let the set $\mathcal{K}$ be a closed subset of $\mathbb{R}^d$. Then, the market is arbitrage-free if and only if
the risk measure $\nu^{\mathcal{H}}$ is sensitive. In this case, $\nu^{\mathcal{H}}$ is continuous from above
and admits the dual representation:
$$\nu^{\mathcal{H}}(\Psi)=\sup_{\mathbb{Q} \in \mathcal{M}_{1,ac}}
\left\{\mathbb{E}_\mathbb{Q}[-\Psi]-\alpha^{\mathcal{H}}(\mathbb{Q}) \right\}.$$
\end{proposition}
\subsubsection{Dynamic Hedging \label{Dynamic hedging - static}}
A natural extension of the previous framework is the multi-period
setting or more generally the continuous-time setting. We briefly
present some results in the latter case. Note that we will come
back to these questions, in the second part of this chapter, under
a slightly different form, assuming that basic
asset prices are It\^o's processes.\\
We now consider a time horizon $T$, a filtration $(\mathcal{F}_t;t
\in [0,T])$ on the probability space $(\Omega, \mathcal{F},
\mathbb{P})$ and a financial market with $n$ basic assets, whose
(non-negative) vector price process $S$ follows a special locally
bounded semi-martingale under $\mathbb{P}$. To avoid arbitrage, we
assume that:
$$\text{\bf
(AAO)}\qquad\text{\em There exists a probability measure}\quad
\mathbb{Q}_0 \sim \mathbb{P}\quad \text{\em such that}\quad S
\quad\text{is a }\>\mathbb{Q}_0-\text{\em local-martingale}.$$ Let
$\mathcal{Q}_{\rm ac}$ be the family of absolutely continuous
martingale measures: $\>\>\mathcal{Q}_{\rm
ac}=\{\>\mathbb{Q}\>|\>\mathbb{Q}\ll\mathbb{P},\quad S\>\>\text{is
a }\>\mathbb{Q}\>\>\text{local-martingale}\}$. {\bf (AAO)} ensures
that the set $\mathcal{Q}_{\rm ac}$ is non empty. Then, as in \index[aut]{Delbaen}Delbaen \cite{Delbaen00},
$\mathcal{Q}_{ac}$ is a closed convex subset of $\mathbb{L}^1(\mathbb{P})$.\\
Let us now introduce dynamic strategies as predictable processes
$\theta$ and their gain processes
$G_t(\theta)=\int_0^t\langle\theta_u,dS_u\rangle=(\theta.S)_t$. We
only consider bounded gain processes and define:
$$\>\Theta^S_T=\{G_T(\theta)=(\theta.S)_T\>|\>\theta.S\>\>\text{is
bounded}\}$$  \index[aut]{Delbaen}Delbaen and
\index[aut]{Schachemayer}Schachermayer have established in
\cite{DelbaenSchachermayer}, as in the static case, the following
duality relationship,
$$\sup\{\>\mathbb{E}_\mathbb{Q}[-X]\>|\> \mathbb{Q}\in \mathcal{Q}_{\rm ac}\}=
\inf\{\>m\>|\>\exists \>\>G_T(\theta)\in \Theta^S_T\>\>s.t.\quad
m+X+G_T(\theta)\geq 0\}$$ Putting $\mathcal{H}=-\Theta^S_T$, this
equality shows that $\nu^{\mathcal{H}}$ is a {\em coherent convex
risk measure continuous from above}. \vspace{-3mm}
\paragraph{Constrained portfolios}
When constraints are introduced on the strategies, everything
becomes more complex. Therefore, we refer to the course held by
\index[aut]{Schied}Schied \cite{Schied_LN} for more details.\\
We assume that hedging positions live in the following convex set:
$$\>\Theta^S_T=\{G_T(\theta)=(\theta.S)_T\>|\>\theta.S\>\>\text{is
bounded by below},\theta \in \mathcal{K}\}$$ The set of
constraints is closed in the following sense: the set $\{\int
\theta dS\> |\> \theta \in \mathcal{K}\}$ is closed in the
semi-martingale or \'{E}mery topology. The optional decomposition
theorem of \index[aut]{F\"{o}llmer}F\"{o}llmer and
\index[aut]{Kramkov}Kramkov \cite{FollmerKramkov} implies the
following dual representation for the risk measure
$\nu^{\mathcal{H}}$:
$$\nu^{\mathcal{H}}(\Psi)=\sup_{\mathbb{Q} \in \mathcal{M}_{1,ac}}
\big\{\mathbb{E}_{\mathbb{Q}}[-\Psi]
-\mathbb{E}_{\mathbb{Q}}[A_T^{\mathbb{Q}}]\big\}$$ where
$A_.^{\mathbb{Q}}$ is the optional process  defined by
$A_0^{\mathbb{Q}}=0$ and $dA_t^{\mathbb{Q}}={\rm ess}\sup_{\xi \in
\mathcal{K}} \mathbb{E}_{\mathbb{Q}}[\theta_tdS_t | \mathcal{F}_t].$\\
The penalty function $\alpha^{\mathcal{H}}$ of the risk measure
$\nu^{\mathcal{H}}$ can be described as
$\mathbb{E}_{\mathbb{Q}}[A_T^{\mathbb{Q}}]$ provided that
$\mathbb{Q}$ satisfies the three following conditions:
\begin{itemize}
\item $\mathbb{Q}$ is equivalent to $\mathbb{P}$; \item Every
process $\theta.S$ with $\theta \in \mathcal{K}$ is a special
semi-martingale under $\mathbb{Q}$; \item $\mathbb{Q}$ admits the
upper variation process $A^{\mathbb{Q}}$ for the set $\{\theta.S
\>|\> \theta \in \mathcal{K}\}$.
\end{itemize}
We can set $\alpha^{\mathcal{H}}(\mathbb{Q})=+\infty$ when one of
these conditions does not hold.
\begin{remark}
Note that there is a fundamental difference between static hedging
with a family of cash flows and dynamic hedging. In the first
case, the initial wealth is a market data: it corresponds to the
(forward) price of the considered cash flows. The underlying logic
is based upon calibration as the probability measures we consider
have to be consistent with the observed market prices of the
hedging instruments. In the dynamic framework, the initial wealth
is a given data. The agent invests it in a self-financing
admissible portfolio which may be rebalanced in continuous time.\\
The problem of dynamic hedging with calibration constraints is a
classical problem for practitioners. This will be addressed in
details after the introduction of the inf-convolution operator.
Some authors have been looking at this question (see for instance
\index[aut]{Bion-Nadal}Bion-Nadal \cite{Bion-Nadal} or
\index[aut]{Cont}Cont \cite{Cont}).
\end{remark}
\section{Dilatation of Convex Risk Measures, Subdifferential and Conservative Price}
\subsection{Dilatation: $\gamma$-Tolerant Risk Measures \label{subsection dilated static risk measure}}
For non-coherent convex risk measures, the impact of the size of
the position is not linear. It seems therefore natural to consider
the relationship between "risk tolerance" and the perception of
the size of the position. To do so, we start from a given root
convex risk measure $\rho$. The risk tolerance coefficient is
introduced as a parameter describing how agents penalize compared
with this root risk measure. More precisely, denoting by $\gamma$
the risk tolerance, we define
\index[not]{$\rho_{\gamma}$}$\rho_{\gamma}$ as:
\begin{equation}
\label{eq:dilated} \rho _{\gamma }( \Psi ) =\gamma \rho \big(
\frac{1}{\gamma }\Psi \big).
\end{equation}
$\rho_{\gamma}$ satisfies a \index[sub]{tolerance}{\it tolerance
property} or a \index[sub]{dilatation} {\it dilatation property}
with respect to the size of the position, therefore it is called
the {\it $\gamma$-tolerant risk
measure}\index[sub]{$\gamma$-tolerant risk measure} associated
with $\rho$ (also called the {\it dilated risk
measure}\index[sub]{dilated risk measure} associated with $\rho$
as in \index[aut]{Barrieu}Barrieu and \index[aut]{El Karoui}El
Karoui \cite{Barrieu-ElKaroui6}). A typical example is the
entropic risk measure where $e_{\gamma}$ is simply the
$\gamma$-dilated of $e_1$. These dilated risk measures satisfy the
following nice property:
\begin{proposition}
\label{Theorem dilated risk measure monotonicity} Let $\big( \rho
_{\gamma },\gamma
>0\big) $ be the family of $\gamma$-tolerant risk measures issued of $\rho$.
Then,\\
$(i)$ The map $\gamma \rightarrow (\rho_\gamma-\gamma \rho(0))$ is non-increasing,\\
$(ii)$ For any $\gamma,\gamma^{'}>0$, $(\rho_{\gamma})_{\gamma ^{'}}=\rho
_{\gamma\,\gamma^{'}}$.\\
$(iii)$ The perspective functional defined on $]0,\infty[\times \mathcal{X}$ by
$$p_{\rho}(\gamma,X)=\gamma\rho(\frac{X}{\gamma})=\rho_{\gamma}(X)$$
is a homogeneous convex functional, cash-invariant with respect to
$X$ (i.e. a coherent risk measure in $X$).
\end{proposition}
\textbf{Proof:} $(i)$ We can take $\rho(0)=0$ without loss of
generality of the arguments. By applying the convexity inequality
to $\frac{X}{\gamma}$ and $0$ with the coefficients
$\frac{\gamma}{\gamma+h}$ and $\frac{h}{\gamma+h}$ $(h>0)$, we
have, since $\rho(0)=0$:
\begin{equation*}
\rho(\frac{X}{\gamma+h})\leq
\frac{\gamma}{\gamma+h}\rho(\frac{X}{\gamma})+
\frac{h}{\gamma+h}\rho(0)\leq\frac{\gamma}{\gamma+h}\rho(\frac{X}{\gamma}).
\end{equation*}
$(ii)$ is an immediate consequence of the definition and
characterization of tolerant risk measures.\\
$(iii)$ The
perspective functional is clearly homogeneous. To show the
convexity, let $\beta_1\in[0,1]$ and $\beta_2=1-\beta_1$ two real
coefficients, and $(\gamma_1,X_1)$, $(\gamma_2,X_2)$ two points in
the definition space of $p_{\rho}$. Then, by the convexity of
$\rho$,
\begin{eqnarray*}
p_{\rho}\big(\beta_1(\gamma_1,X_1)+\beta_2(\gamma_2,X_2)\big)&=&²(\beta_1\gamma_1+\beta_2\gamma_2)\>\rho\Big(\frac{\beta_1
X_1+\beta_2 X_2}
{\beta_1\gamma_1+\beta_2\gamma_2}\Big)\\
&\leq& (\beta_1\gamma_1+\beta_2\gamma_2)\>\Big[\frac{\beta_1
\gamma_1}{\beta_1\gamma_1+\beta_2\gamma_2}\rho(\frac{X_1}{\gamma_1})+\frac{\beta_2
\gamma_2}{\beta_1\gamma_1+\beta_2\gamma_2}\rho(\frac{X_2}{\gamma_2})\Big]
\\
&\leq&
\beta_1\>\rho_{\gamma_1}(X_1)\>+\>\beta_2\>\rho_{\gamma_2}(X_2).
\end{eqnarray*}
The other properties are obvious.$\quad\square$
\vspace{2mm}\\
So, we naturally are looking for the asymptotic behavior of the
perspective risk measure when the risk tolerance either tends to
$+\infty$ or tends to $0$.
\subsection{Marginal Risk Measures and Subdifferential}
\subsubsection{Marginal Risk Measure}
Let us first observe that $\rho $ is a \textit{coherent} risk
measure if and only if $\rho _{\gamma }\equiv \rho $. We then
consider the behavior of the family of $\gamma$-tolerant risk
measures when the tolerance becomes infinite.
\begin{proposition}
\label{prop:marginriskmeasure} Suppose that $\rho( 0) =0$, or
equivalently
$\alpha(\mathbf{Q})\geq 0 \quad \forall \mathbf{Q}\in \mathbf{M}_{1,f} $.\\
$a)$ The \index[sub]{marginal risk measure} marginal risk measure
\index[not]{$\rho_{\infty}$} $\rho _{\infty }$, defined as the non-increasing limit of
$\rho_{\gamma}$ when $\gamma$ tends to infinity, is a coherent risk measure with penalty
function $\alpha_{\infty}=\lim_{\gamma \rightarrow +\infty}(\gamma \alpha)$ that is:
$$\begin{array}{llll}
\alpha_{\infty}(\mathbf{Q})&:=\sup_{\Psi}\big\{
\mathbb{E}_{\mathbf{Q}}[-\Psi]-\rho_{\infty}(\Psi)\big\}=\>\>0 \quad{\rm if}\>\>
\alpha(\mathbf{Q})=0\>\>, \quad +\infty \>\>{\rm if \>not},\\
 \rho _{\infty }( \Psi ) &=
\sup_{\mathbf{Q\in }\mathbf{M}_{1,f} } \big\{
\mathbb{E}_{\mathbf{Q}}[-\Psi]\>\big|\>\> \alpha \big(
\mathbf{Q}\big) =0\big\}.
\end{array}
$$
$b)$ Assume now that $\rho$ is a
$\mathbb{L}_{\infty}(\mathbb{P})$-risk measure such that
$\rho(0)=0$.\\
If $\rho$ is continuous from below, the $\rho_{\infty}$ is continuous from below and admits
a representation in terms of absolutely continuous probability measures as:
$$\rho_{\infty}(\Psi)= \mathbf{\max}_{\mathbb{Q\in }\mathcal{M}_{1,ac} }
\big\{ \mathbb{E}_{\mathbb{Q}}[-\Psi]\>\big|\>\> \alpha \big( \mathbb{Q}\big) =0\big\},$$
and the set $\big\{\mathbb{Q} \in \mathcal{M}_{1,ac}\>\big|\>\> \alpha \big( \mathbb{Q}\big)
=0\big\}$ is non empty, and weakly compact in $\mathbb{L}^{1}(\mathbb{P})$.
\end{proposition}
\noindent \textbf{Proof:} $a)$ Thanks to Theorem \ref{Theorem
dilated risk measure monotonicity}, for any $\Psi\in \mathcal{X}$
$\rho _{\gamma }(\Psi)\searrow \rho _{\infty }(\Psi)$ when $\gamma
\rightarrow +\infty$. Given the fact that $-m \geq
\rho_{\gamma}(\Psi) \geq -M$ when $m \leq \Psi \leq M$, we also
have $-m \geq \rho_{\infty}(\Psi) \geq -M$ and $\rho_{\infty}$ is
finite.\\
Convexity, monotonicity and cash translation invariance properties
are preserved when taking the limit. Therefore, $\rho_{\infty}$ is
a convex risk measure with $\rho_{\infty}(0)=0$.\\
Moreover, given that $(\rho_{\delta})_{\gamma}=\rho_{\delta
\gamma}=(\rho_{\gamma})_{\delta}$, we have that
$(\rho_{\delta})_{\infty}=\rho_{\infty}=(\rho_{\infty})_{\delta}$
and $\rho_{\infty}$ is a coherent risk measure.\vspace{1mm}\\
Since $\alpha \geq 0$, the minimal penalty function is:
$$
\begin{array}{lll} \alpha_{\infty}(\mathbf{Q})&=\sup_{\xi} \big\{
\mathbb{E}_{\mathbf{Q}}[-\xi]-\rho_{\infty}(\xi)\big\}\\
&=\sup_{\xi}\sup_{\gamma>0}\big\{
\mathbb{E}_{\mathbf{Q}}[-\xi]-\gamma\rho(\frac{\xi}{\gamma})\big\}\\
&=\sup_{\gamma >0}\big\{\gamma \alpha(\mathbf{Q})\big\}= \>\>0\>\>
\>{\rm if}\>\> \alpha(\mathbf{Q})=0\>\>, \quad +\infty \>\>{\rm if
\>not}.
\end{array}
$$
Moreover, $\alpha_{\infty}$ is not identically equal to $+\infty$
since the set $\big\{\mathbf{Q} \in \mathbf{M}_{1,f}\>\big|\>\>
\alpha \big( \mathbf{Q}\big) =0\big\}$ is not empty given that
$\rho(0)=0=\max\big\{-\alpha(\mathbf{Q})
\big\}=-\alpha(\mathbf{Q}_0)$ for some additive measure
$\mathbf{Q}_0 \in \mathbf{M}_{1,f}$ , from Theorem
\ref{Theorem dual representation}.\\
Assume now that $\rho$ is continuous from below and consider a
non-decreasing sequence $(\xi_n \in \mathcal{X})$ with limit $\xi
\in \mathcal{X}$. By monotonicity,
$$\rho_{\infty}(\xi)=\inf_{\gamma}\rho_{\gamma}(\xi)
=\inf_{\gamma}\inf_{\xi_n}\rho_{\gamma}(\xi_n)=\inf_{\xi_n}\inf_{\gamma}\rho_{\gamma}(\xi_n)
=\inf_{\xi_n}\rho_{\infty}(\xi_n).$$
Then, $\rho_{\infty}$ is also continuous from below.\vspace{2mm}\\
$b)$ When $\rho$ is a $\mathbb{L}_{\infty}(\mathbb{P})$-risk
measure, continuous from below, $\rho$ is also continuous from
above and the dual representation holds in terms of absolutely
continuous probability measures. Using the same argument as above,
we can prove that $\rho_{\infty}$ is a coherent
$\mathbb{L}_{\infty}(\mathbb{P})$-risk measure, continuous from
below with minimal penalty function:
$$\begin{array}{ll}
\alpha_{\infty}(\mathbb{Q})&= 0 \>\> {\rm if}\>\>
\alpha(\mathbb{Q})=0 \quad {\rm and} \quad \mathbb{Q} \in
\mathcal{M}_{1,ac}\\
&= +\infty \>\> {\rm otherwise}.
\end{array}$$
Moreover, thanks to Theorem \ref{Theorem rep duale Mac}, the set $\big\{\mathbb{Q} \in
\mathcal{M}_{1,ac}\>\big|\>\> \alpha \big( \mathbb{Q}\big) =0\big\}$ is non empty and weakly
compact in $\mathbb{L}^{1}(\mathbb{P})$. $\square$\vspace{2mm}\\
To have some intuition about the interpretation in terms of marginal risk measure, it is
better to refer to the risk aversion coefficient $\epsilon=1/\gamma$. $\rho _{\infty
}(\Psi)$ appears as the limit of $\frac{1}{\epsilon}\big(\rho(\epsilon \Psi)-\rho(0)\big)$,
i.e. the right-derivative at $0$ in the direction of $\Psi$ of the risk measure $\rho$, or
equivalently, the marginal risk measure. For instance $e_{\infty}(\Psi) =
\mathbb{E}_{\mathbb{P}}(-\Psi).$\vspace{1mm}\\
In some cases, and in particular when the set
$\mathcal{Q}_{\infty}^{\alpha }$ has a single element, the pricing
rule $\rho_{\infty}(-\Psi)$ is a linear pricing rule and can be
seen as an extension of the notion of marginal utility pricing and
of the Davis price (see \index[aut]{Davis}Davis \cite{Davis97} or
\index[aut]{Karatzas} Karatzas and \index[aut]{Kou}Kou \cite{KK}).
\subsubsection{Subdifferential and its Support Function}
\paragraph{Subdifferential}
Let us first recall the definition of the subdifferential of a
convex functional.
\begin{definition}
Let $\phi$ be a convex functional on $\mathcal{X}$. The
subdifferential of $\phi$ at $X$ is the set $$\partial
\phi(X)=\big\{\mathbf{q} \in \mathcal{X}'\>|\>\>\forall X\in
\mathcal{X}, \>\> \phi(X+Y) \geq \phi(X)+\mathbf{q}(-Y)\big\}$$
\end{definition}
The subdifferential of a convex risk measure $\rho$ with penalty
function $\alpha(q)=\sup_{Y}\{q(-Y)-\rho(Y)\}$ is included in
${\bf Dom}(\alpha)$ since when $ q\in \partial \rho(\xi)$, then
$\alpha(q)-(q(-\xi)-\rho(\xi))\leq 0$. So, we always refer to
finitely additive measure $\mathbf{Q}$ when working with risk
measure subdifferential. In fact, we have the well-known
characterization of the subdifferential:\\
{\em  $\mathbf{q} \in
\partial \rho(\xi)$ if and only if $\mathbf{q}\in
\mathbf{M}_{1,f}$ is  optimal for the maximization program
$\mathbb{E}_{\mathbf{Q}}[-\xi]-\alpha(\mathbf{Q}) \>\>
\longrightarrow \>\>\max_{\mathbf{Q}\in
\mathbf{M}_{1,f}}$}.\vspace{2mm}\\
We can also relate it with the notion of marginal risk measure,
when the root risk measure is now centered around a given element
$\xi \in \mathcal{X}$, i.e. $\rho_\xi(X)=\rho(X+\xi)-\rho(\xi)$,
by defining:
$$\rho_{\infty,\xi}(\Psi)
\equiv \lim_{\gamma \rightarrow +\infty} \gamma
\Big(\rho\big(\xi+\frac{\Psi}{\gamma}\big)-\rho(\xi)\Big).$$ Using
Proposition \ref{prop:marginriskmeasure}, since the $\rho_\xi$
penalty function is $\alpha_{\xi}(\mathbf{Q}) \equiv
\alpha(\mathbf{Q})-\mathbb{E}_{\mathbf{Q}}[-\xi]+\rho(\xi)$,
$\rho_{\infty,\xi}$ is coherent and
$$\rho_{\infty,\xi}(\Psi)=\sup_{\mathbf{Q\in }\mathbf{M}_{1,f} }
\big\{ \mathbb{E}_{\mathbf{Q}}[-\Psi]\>\big| \>\> \rho(\xi)=
\mathbb{E}_{\mathbf{Q}}[-\xi]-\alpha(\mathbf{Q})\big\}.$$
\begin{proposition}
The coherent risk measure $\rho_{\infty,\xi}(\Psi) \equiv
\lim_{\gamma \rightarrow +\infty} \gamma
\Big(\rho\big(\xi+\frac{\Psi}{\gamma}\big)-\rho(\xi)\Big)$ is the
support function of the subdifferential $\partial \rho(\xi)$ of
the convex risk measure $\rho$ at $\xi$:
$$\rho_{\infty,\xi}(\Psi)=\sup_{\mathbf{Q}\in \mathbf{M}_{1,f} } \big\{
\mathbb{E}_{\mathbf{Q}}[-\Psi]\>\big| \>\>\rho(\xi)=
\mathbb{E}_{\mathbf{Q}}[-\xi]-\alpha(\mathbf{Q})\big\}=\sup_{\mathbf{Q}\in
\partial \rho(\xi)}\mathbb{E}_{\mathbf{Q}}[-\Psi]$$
\end{proposition}
{\bf Proof:} From the definition of the subdifferential,
$$\begin{array}{lll}
\partial \rho(\xi)&=\big\{\mathbf{q} \in \mathcal{X}'\>|\>\>\forall \Psi \in \mathcal{X}, \>\>
\rho(\xi+\Psi) \geq \rho(\xi)+\mathbf{q}(-\Psi)\big\}\\
&=\big\{\mathbf{q} \in \mathcal{X}'\>|\>\>\forall \Psi \in
\mathcal{X},\>\>\rho_{\infty,\xi}(\Psi) \geq
\mathbf{q}(-\Psi)\big\}\\
&=\partial \rho_{\infty,\xi}(0).
\end{array}$$
But $q\in \partial \rho_{\infty,\xi}(0)$ iff
$\alpha_{\infty,\xi}(\mathbf{Q}_{\mathbf{q}})=0$. So the proof is
complete. $\quad \square$
\paragraph{The $\mathbb{L}_{\infty}(\mathbb{P})$ case:}
When working with $\mathbb{L}_{\infty}(\mathbb{P})$-risk measures, following
\index[aut]{Delbaen}Delbaen \cite{Delbaen00} (Section 8), the natural definition of the
subdifferential is the following:
$$\partial \rho(\xi)=\big\{f \in \mathbb{L}^1(\mathbb{P})\>|\>\>\forall \Psi \in \mathbb{L}_{\infty}(\mathbb{P})
, \>\> \rho(\xi+\Psi) \geq
\rho(\xi)+\mathbb{E}_{\mathbb{P}}[f(-\Psi)]\big\}$$ Using the same
arguments as above, we can prove that every $f \in \partial
\rho(\xi)$ is non-negative with a $\mathbb{P}$-expectation equal
to $1$. Since $\partial \rho(\xi)$ is also the subdifferential of
$\rho_{\infty,\xi}(0)$, the properties of $\partial \rho(\xi)$ may
be deduced from those of the coherent risk measure
$\rho_{\infty,\xi}$, for which we have already shown that if
$\rho$ is continuous from below and $\rho(0)=0$ then for any
$\xi$, the effective domain of $\alpha_{\infty,\xi}$ is non empty.
Then, under this assumption, $\partial \rho(\xi)$ is non empty and
we have the same characterization of the subdifferential as:
$$\mathbb{Q} \in \partial \rho(\xi) \Longleftrightarrow
\rho(\xi)=\mathbb{E}_{\mathbb{Q}}[-\xi]-\alpha(\mathbb{Q}).$$
We
now summarize these results in the following proposition:
\begin{proposition}
Let $\rho$ be a $\mathbb{L}_{\infty}(\mathbb{P})$-risk measure,
continuous from below. Then, for any $\xi \in
\mathbb{L}_{\infty}(\mathbb{P})$, $\rho_{\infty,\xi}$ is the
support function of the non empty subdifferential $\partial
\rho(\xi)$, i.e.:
$$\rho_{\infty,\xi}(\Psi)=\sup\big\{\mathbb{E}_{\mathbb{Q}}[-\Psi]\>;\>
\mathbb{Q} \in \partial \rho(\xi)\big\}.$$ and the supremum is attained by some
$\mathbb{Q}\in \partial \rho(\xi)$.
\end{proposition}
\subsection{Conservative Risk Measures and Super-Price \label{Conservative risk measures}}
We now focus on the properties of the $\gamma$-tolerant risk
measures when the risk tolerance coefficient tends to $0$ or
equivalently when the risk aversion coefficient goes to $+\infty$.
The \index[sub]{conservative risk measure} conservative risk
measures that are then obtained can be reinterpreted in terms of
super-pricing rules\index[sub]{super price}. Using vocabulary from
convex analysis, these risk measures are related to recession (or
asymptotic) functions.
\begin{proposition}
$(a)$ When $\gamma$ tends to $0$, the family of $\gamma$-tolerant
risk measures $(\rho_{\gamma})$ admits a limit
\index[not]{$\rho_{0^+}$} $\rho_{0^+}$, which is a coherent risk
measure. This conservative risk measure $\rho_{0^+}$ is simply the
''super-price'' of $-\Psi $:
\begin{equation*}
\rho _{0^+}( \Psi ) =\lim_{\gamma\downarrow 0}\nearrow (\rho
_{\gamma}( \Psi )-\gamma \rho(0))=\sup_{\mathbf{Q\in }\mathbf{M%
}_{1,f} } \big\{\mathbb{E}_{\mathbf{Q}}[-\Psi]\big|\>\alpha \big(
\mathbf{Q}\big)  <\infty \big\}.
\end{equation*}
Its minimal penalty function is
$$\alpha_{0^+}(\mathbf{Q})=0 \>\>
if\>\> \alpha(\mathbf{Q})<+\infty \quad {\rm and} \quad =+\infty
\>\> {\rm if \>\>not}.$$ $(b)$ If $\rho$ is continuous from above
on $\mathbb{L}^{\infty}(\mathbb{P})$, then $\rho_{0^+}$ is
continuous from above and
\begin{equation*}
\rho _{0^+}( \Psi ) =\sup_{\mathbb{Q\in }\mathcal{M%
}_{1,ac} } \big\{\mathbb{E}_{\mathbb{Q}}[-\Psi]\big|\>\alpha \big(
\mathbb{Q}\big)  <\infty \big\}.
\end{equation*}
\end{proposition}
\textbf{Proof:} Let us first observe that $\rho_{\gamma}(\xi) = \gamma \big(
\rho(\frac{\xi}{\gamma})-\rho(0)\big) + \gamma \rho(0)$ is the sum of two terms. The first
term is monotonic while the second one goes
to $0$.\\
The functional $\rho _{0^+}$ is coherent (same proof as for $\rho
_{\infty}$) with the acceptance set $\mathcal{A}_{\rho_{0^+}}
=\{\xi,\>\> \forall \lambda \geq 0, \>\> \lambda \xi \in \mathcal{A}_{\rho}-\rho(0)\}$. \\
On the other hand,
 by monotonicity, the minimal penalty function $\alpha_{0^+}\geq \gamma \alpha\geq 0$ ; so, $\alpha_{0^+}(\mathbf{Q})=0$
on ${\rm {\bf Dom}}(\alpha)$, and
$\alpha_{0^+}(\mathbf{Q})=+\infty$ if not. In other words,
$\alpha_{0^+}$ is the convex indicator of ${\rm {\bf Dom}}(\alpha)$.\\
If $\rho$ is continuous from above on
$\mathbb{L}^{\infty}(\mathbb{P})$, then the same type of dual
characterization holds for $\rho _{0^+}$ but in terms of $\mathcal{M%
}_{1,ac}$. So, $\alpha_{0^+}(\mathbb{Q})=0$ on ${\rm {\it
Dom}}(\alpha)$, and
$\alpha_{0^+}(\mathbb{Q})=+\infty$ if not.\vspace{1mm}\\
We could have proved directly the continuity from above of $\rho
_{0^+}$, since $\rho_{0^{+}}$ is the non-decreasing limit of
continuous from above risk measures $(\rho_{\gamma}-\gamma
\rho(0))$. $\square$
\begin{remark} A nice illustration of this result
can be obtained when considering the entropic risk measure
$e_{\gamma}$. In this case, it comes immediately that
$e_{0^+}(\Psi)=\sup_{\mathbb{Q}}\big\{\mathbb{E}_{\mathbb{Q}}[-\Psi]\big|\>h(\mathbb{Q}\,|\,\mathbb{P})<+\infty\big\}
=\mathbb{P}-\text{ess}\sup(-\Psi)=\rho_{\rm max}(\Psi)$ where
$\rho_{\rm max}$ is here the
$\mathbb{L}_{\infty}(\mathbb{P})$-worst case measure. This also
corresponds to the weak
super-replication price as defined by Biagini and Frittelli in \cite{BiaginiFrittelli}.\\
Note that this conservative risk measure $e_{0^{+}}(\Psi)$ cannot
be realized as $\mathbb{E}_{\mathbb{Q}_0}[-\Psi]$ for some
$\mathbb{Q}_0 \in \mathcal{M}_{1,ac}$. It is a typical example
where the continuity from below fails.
\end{remark}
\section{Inf-Convolution \label{Subsection static inf-convolution}}
A useful tool in convex analysis is the inf-convolution operation.
While the classical convolution acts on the Fourier transforms by
addition, the inf-convolution acts on Fenchel transforms by
addition as we would see later.
\subsection{Definition and Main Properties \label{subsection properties inf-convolution}}
The inf-convolution of two convex functionals $\phi_A$ and
$\phi_B$ may be viewed
 as the functional value of the minimization program
\begin{equation}
\phi_{A,B}(X)=\inf_{H\in \mathcal{X}}\big\{ \phi _{A}(X-H) +\phi
_{B}( H) \big\}, \label{general inf-convolution pg}
\end{equation}
This program is the functional extension of the classical
inf-convolution operator acting on real convex functions $f\square
g(x)=\inf_y \{f(x-y)+g(y) \}$.
\paragraph{Illustrative example:}
Let us assume that the risk measure $\rho_A$ is the linear one
$q_A(X)=\mathbb{E}_{\mathbb{Q}_A}[-X]$, whose the penalty function
is the functional $\alpha_A(\mathbb{Q})=\>0\quad \text{if}\quad
\mathbb{Q}=\mathbb{Q}_A,\quad =+\infty \quad \text{if not}$. Given
a convex risk measure, $\rho_B$, with penalty functional
$\alpha_B$, we deduce from the definition of the inf-convolution
that
$$q_A\square\rho_B(X)=q_A(-X)-\alpha_B(\mathbb{Q}_A)$$
$\diamond$ Then, $q_A\square\rho_B$ is identically $-\infty$ if
$\alpha_B(\mathbb{Q}_A)=+\infty$.\\
 $\diamond$ If it is not the case, the minimal penalty
function $\alpha_{A,B}$ associated with this measure is:
$$\alpha_{A,B}(\mathbb{Q})=\alpha_B(\mathbb{Q}_A)+\alpha_A(\mathbb{Q})=\alpha_B(\mathbb{Q})+\alpha_A(\mathbb{Q})$$
$\diamond$ Moreover, the infimum is attained in the
inf-convolution program by any $H^*$ such that
$$\alpha_B(\mathbb{Q}_A)=\mathbb{E}_{\mathbb{Q}_A}[-H^*]-\rho_B(H^*)$$
that is $H^*$ is optimal for the maximization program defining the
$\alpha_B$. \\
In terms of subdifferential, we have the first order condition:
$\quad \mathbb{Q}_A\in
\partial\rho_B(H^*)$.
\subsubsection{Inf-Convolution and Duality}
In our setting, convex functionals are generally convex risk
measures, but we have also been concerned by the convex indicator
of convex subset, taking infinite values. In that follows, we
already assume that convex functionals $\phi$ we consider are
proper (i.e. not identically $+ \infty$) and in general closed or
lower semicontinuous (in the sense that the level sets
$\{X|\>\phi_B(X)\leq c\},\>\>c\in \mathbb{R}$ are weak*-closed). To be
consistent with the risk measure notations we define their Fenchel
transforms on $\mathcal{X}'$ as
$$\beta(q)=\sup_{X\in \mathcal{X}}\{q(-X)-\phi(X)\}.$$
When the linear form $q$ is related to an additive finite measure
$\mathbf{Q}\in \mathbf{M}_{1,f}$, we use the notation
$q_{\mathbf{Q}}(X)=\mathbb{E}_{\mathbf{Q}}[X]$. For a general
treatment of inf-convolution of convex functionals, the interested
reader may refer to the highlighting paper of
\index[aut]{Borwein}Borwein and \index[aut]{Zhu}Zhu
\cite{Borwein-Zhu}. The following theorem
extends these results to the inf-convolution of convex functionals
whose one of them at least is a convex risk measure:
\begin{theorem}
\label{Theorem modification de rho1 par A_rho2} Let $\rho _{A}$ be
a convex risk measure with penalty function $\alpha _{A}$ and
$\phi_B$ be a proper closed convex functional with Fenchel
transform $\beta$. Let $\rho _{A}\square \phi _{B}$ be the
inf-convolution of $\rho _{A}$ and $\phi _{B}$ defined as
\begin{equation}
\label{eq: def inf-convolution measures}
 X \rightarrow \rho _{A}\square \phi _{B}(X ) =\inf_{H\in
\mathcal{X}}\big\{ \rho _{A}(X -H) +\phi_{B}( H) \big\}
\end{equation}
and assume that $\rho _{A}\square \phi _{B}( 0) >-\infty $.
Then,\\
$\bullet$ $\rho _{A}\square \phi _{B}$ is a convex risk measure
which is finite for all $X\in \mathcal{X}$.\\ $\bullet$ The
associated penalty function $\alpha _{A,B}$ takes the value
$+\infty$ for any $q$ outside of $\mathbf{M}_{1,f}$, and
$$\begin{array}{ll}
& \forall \mathbf{Q} \in \mathbf{M}_{1,f} \quad \alpha _{A,B} \big( %
\mathbf{Q}\big) =\alpha _{A}\big( \mathbf{Q}\big) +\beta_{B}\big( %
q_{\mathbf{Q}}\big),\\
\mbox{\rm and} \quad & \exists \mathbf{Q} \in \mathbf{M}_{1,f}
\quad
\mbox{\rm s.t.}\quad \alpha _{A}\big( \mathbf{Q}\big) +\beta_{B}\big(q_{ %
\mathbf{Q}}\big)<\infty.
\end{array}
$$
$\bullet$ Moreover, if the risk measure $\rho _{A}$ is continuous
from below, then $\rho _{A}\square \phi _{B}$ is also continuous
from below.
\end{theorem}
\textbf{Proof:} We give here the main steps of the proof of this
theorem.\\
$\diamond$ The monotonicity and translation invariance properties
of $\rho _{A}\square \phi _{B}$ are immediate from the definition,
since at least one
of the both functionals have these properties.\\
$\diamond$ The convexity property simply comes from the fact that,
for any $X_{A}$, $X_{B}$, $H_{A}$ and $H_{B}$ in $\mathcal{X}$ and
any $\lambda
\in [ 0,1] $, the following inequalities hold as $\rho _{A}$ and $%
\rho_{B}$ are convex functionals,
$$
\begin{array}{ccc}
\rho _{A}\big((\lambda X_{A}+( 1-\lambda )X_{B})-( \lambda H_{A}+(
1-\lambda ) H_{B}) \big) &\leq& \lambda \rho _{A}\big(X
_{A}-H_{A}\big) +( 1-\lambda) \rho
_{A}\big(X_{B}-H_{B}\big) \\
\phi_{B}\big(\lambda H_{A}+( 1-\lambda) H_{B}\big) &\leq &\lambda
\phi_{B}( H_{A}) +( 1-\lambda) \phi_{B}\big(H_{B}\big).
\end{array}
$$
By adding both inequalities and taking the infimum in $H_{A}$ and
$H_{B}$ on the left-hand side and separately in $H_{A}$ and in
$H_{B}$ on the right-hand side, we obtain:
\begin{equation*}
\rho _{A}\square \phi_{B}\big( \lambda X_{A}+( 1-\lambda)X
_{B}\big) \leq \lambda\, \rho _{A}\square \phi_{B}(X_{A}) +(
1-\lambda) \rho _{A}\square \rho _{B}(X_{B}).
\end{equation*}
$\diamond$ Using Equation \eqref{penalty 2}, the
associated penalty function is given, for any $\mathbf{Q}\mathbb{\in }%
\mathbf{M}_{1,f}$, by
$$
\begin{array}{rlll}
\alpha _{A,B}\big( \mathbf{Q}\big)  &=\sup_{X\in
\mathcal{X}}\big\{ \mathbb{E}_{\mathbf{Q}}[-X] -\rho _{A,B}(
X) \big\} \\
&=\sup_{\Psi \in \mathcal{X}}\big\{ \mathbb{E}_{\mathbf{Q}}[-X]
-\inf_{H\in \mathcal{X}}\big\{ \rho _{A}(X -H) +\phi_{B}( H)
\big\} \big\}\\ \vspace{2mm}
& =\sup_{X\in \mathcal{X}}\sup_{H\in \mathcal{X}}\big\{ \mathbb{E}_{%
\mathbf{Q}}\big[-(X -H) \big] +\mathbb{E}_{\mathbf{Q}%
}[-H] -\rho _{A}(X -H) -\phi_{B}( H)
\big\}\\
{\rm by \>letting} \quad \widetilde{X}\triangleq X -H\in \mathcal{X}\\
&=\sup_{\widetilde{X}\in \mathcal{ X}}\sup_{H\in \mathcal{X}}\big(
\mathbb{E}_{\mathbf{Q}}[-\widetilde{X}] -\rho _{A}( \widetilde{X})
+\mathbb{E}_{\mathbf{ Q}}[-H] -\phi_{B}( H) \big)\vspace{1mm}
 =\alpha _{A}\big( \mathbf{Q}\big)
+\beta_{B}\big( q_{\mathbf{Q}}\big).
\end{array}$$
When $q \not\in \mathbf{M}_{1,f}$, the same equalities hold true.
Since $\rho_A$ is a convex risk measure, $\alpha _{A}(q)=+\infty$,
and since $\beta$ is a proper functional, $\beta(q)$ is dominated
from below; so, $\alpha _{A,B}(q)=+\infty$. This equality $\alpha
_{A,B}=\alpha_A+\beta_B$ holds even when $\alpha _{A}$ and
$\beta_{B}$ they take infinite values.\\
 $\diamond$ The continuity
from below is directly obtained upon considering an increasing
sequence of $\left(X_{n}\right) \in \mathcal{X}$ converging to
$X$. Using the monotonicity property, we have
\begin{eqnarray*}
\inf_{n}\rho _{A}\square \phi_{B}\left(X_{n}\right)
&=&\inf_{n}\inf_{H}\left\{ \rho _{A}\left(X_{n}-H\right) +\phi
_{B}\left( H\right) \right\} \\
&=&\inf_{H}\inf_{n}\left\{ \rho _{A}\left(X_{n}-H\right) +\phi
_{B}\left( H\right) \right\} =\inf_{H}\left\{ \rho_{A}\left(X
-H\right)
+\phi_{B}\left( H\right) \right\} \\
&=&\rho _{A}\square \phi_{B}\left(X\right).\hspace{5mm}
\square\vspace{-3mm}
\end{eqnarray*}
We can now give an inf-convolution interpretation of the convex
risk measure $\nu^{\mathcal{H}}$ generated by a convex set
$\mathcal{H}$ as in Corollary \ref{re:worst risk measure and
hedging} as the inf-convolution of the convex indicator function
of $\mathcal{H}$, and the worst-case risk measure. This
regularization may be applied at any proper convex functional.
\begin{proposition}\label{proposition regularisation par worst}
\rm{[Regularization by inf-convolution with
$\rho_{\text{worst}}$]} Let
$\rho_{\text{worst}}(X)=\sup_{\omega}(-X(\omega))$ be the worst
case risk measure.\\
$i)$ $\rho_{\text{worst}}$ is a neutral element for the infimal
convolution of convex risk measures.\\
$ii)$ Let $\mathcal{H}$ be
a convex set such that $\inf\{m\>|\,\exists\>\xi \in
\mathcal{H}\}> -\infty.$ The convex risk measure generated by
$\mathcal{H}$, $\nu^{\mathcal{H}}$ is the inf-convolution of the
convex indicator functional of $\mathcal{H}$ with the worst case
risk measure,
$$\nu^{\mathcal{H}}=\rho_{\text{worst}}\square l^{\mathcal{H}}$$
$iii)$ More generally, let $\phi$ be a proper convex functional,
such that for any $H$, $\phi(H)\geq -\sup_\omega H(\omega)-c.$\\
The infimal convolution of $\rho_{\text{worst}}$ and $\phi$,
$\rho_{\phi}=\rho_{\text{worst}}\square \phi$ is the largest
convex risk measure dominated by $\phi$.\\
$iv)$ Let $\beta$ the penalty functional associated with $\phi$.
Then, the penalty functional associated with $\rho_{\phi}$ is the
functional $\alpha_{\phi}$, restriction of $\beta$ at the set
$\mathcal{M}_{1,f}$,
$$
\begin{array}{llll}
\alpha_{\phi}(q)&=\beta(q)+l^{\mathcal{M}_{1,f}}(q)&\\
&=\beta(q_{\mathbf{Q}})\quad \text{if }\quad q_{\mathbf{Q}}\in
\mathcal{M}_{1,f},&\quad +\infty \quad \text{if not}.
\end{array}
$$
$v)$ Given a general risk measure $\rho_A$ such that
$\rho_A\square \phi(0)>-\infty$, then
$$\rho_A\square \phi=\rho_A\square\rho_{\text{worst}}\square
\phi=\rho_A\square\rho_{\phi}.$$
\end{proposition}
{\bf Proof:} We start by proving that
$\rho\square\rho_{\text{worst}}=\rho$. By definition,
$$
\begin{array}{llll}
\rho\square\rho_{\text{worst}}(X)&=\inf_Y\{\sup_\omega(-Y(\omega))+\rho(X-Y)\}
=\inf_Y\big\{\rho\big(X-(Y-\sup_\omega(-Y(\omega)))\big)\big\}\\
&=\inf_{Y\geq 0}\{\rho(X-Y)\}=\rho(X)
\end{array}
$$
To conclude, we have used the cash invariance of $\rho$ and the
fact that $\rho(X-Y)\geq
\rho(X)$ whenever $Y\geq 0$.\\
ii) has been proved in Corollary \ref{re:worst
risk measure and hedging}.\\
iii) By Theorem \ref{Theorem modification de rho1 par A_rho2},
$\rho_{\phi}=\rho_{\text{worst}}\square \phi$ is a convex risk
measure. Since $\rho_{\text{worst}}$ is a neutral element for the
inf-convolution of risk measure, any risk measure $\rho$ dominated
by $\phi$ is also dominated by $\rho_{\text{worst}}\square \phi$
since $\rho=\rho_{\text{worst}}\square \rho \leq
\rho_{\text{worst}}\square \phi=\rho_\phi$. Hence the result.\hspace{5mm}$\square$\vspace{2mm}\\
Therefore, in the following, we only consider the infimal
convolution of convex risk measures. The following result makes
more precise Theorem \ref{Theorem modification de rho1 par A_rho2}
and plays a key role in our analysis.
\begin{theorem}\label{Sandwich theorem}\rm{[Sandwich Theorem]} Let $\rho_A$ and $\rho_B$ be two convex risk measures.\\
Under the assumptions of Theorem \ref{Theorem modification de rho1
par A_rho2} (i.e. $\rho_{A,B}(0)=\rho_A \square \rho_B(0)>-\infty$),\\
$i)$ There exists $\mathbf{Q}\in \partial \rho_{A,B}(0)$ such that,
for any $X$ and any $Y$,
$$ \rho_A \square\rho_B(0)\leq
\big(\rho_A(X)-\mathbb{E}_{\mathbf{Q}}[-X]\big)+\big(\rho_B(Y)-\mathbb{E}_{\mathbf{Q}}[-Y]\big).$$
$ii)$ Assume $\rho_{A,B}(0)\geq c$. There is an affine function,
$a_{\mathbf{Q}}(X)=-\mathbb{E}_{\mathbf{Q}}[-X]+r$, with
$\mathbf{Q}\in
\partial \rho_{A,B}(0)$, satisfying
\begin{equation}
\label{eq:sandwich} \rho_A(.)\geq a_{\mathbf{Q}}\geq
-\rho_B(-.)+c.
\end{equation}
Moreover, for any $\overline{X}$ such that $\rho_A(\overline{X})
+\rho_B(-\overline{X})=\rho_{A,B}(0)$, $\>\mathbf{Q}\in \partial
\rho_B(-\overline{X})\cap \partial \rho_A(\overline{X})$. The
inf-convolution is said to be exact at
$\overline{X}$.\\
$iii)$ {\rm Interpretation of the Condition $\rho_A \square
\rho_B(0)>-\infty$.}\\
The following properties are equivalent:\\
$\diamond$ $\rho_A \square \rho_B(0)>-\infty.$\\
$\diamond$ The sandwich
property \eqref{eq:sandwich} holds for some affine function
$a_{\mathbf{Q}}(X)=-\mathbb{E}_{\mathbf{Q}}[-X]+r$.\\
$\diamond$ There
exists  $\mathbf{Q}\in
\mathbf{Dom}(\alpha_A)\cap\mathbf{Dom}(\alpha_B).$
$\diamond$Let
$\rho^A_{0^+}$ (resp. $\rho^B_{0^+}$) be the conservative risk
measure associated with $\rho^A$ (resp. $\rho^B$). Then
$$\rho^A_{0^+}(X)+\rho^B_{0^+}(-X) \geq 0.$$
\end{theorem}
Before proving this Theorem, let us make the following comment:
the inf-convolution risk measure $\rho _{A,B}$, given in Equation
(\ref{eq: def inf-convolution measures}) may also be defined, for
instance, as the value functional of the program
\begin{equation*}
\rho _{A,B}\left( \Psi \right)=\rho_A \square \rho_B(\Psi)=\rho_A \square \nu^{\mathcal{A}_{\rho_B}}(\Psi)
=\inf \left\{ \rho _{A}\left( \Psi
-H\right) ,H\in \mathcal{A}_{\rho _{B}}\right\},
\end{equation*}
where $\nu^{\mathcal{A}_{\rho_B}}$ is the risk measure with acceptance set $\mathcal{A}_{\rho_B}$. This emphasizes again
the key role played the risk measures generated by a convex set, if needed.\vspace{2mm}\\
{\bf Proof:} $i)$ By Theorem \ref{Theorem modification de rho1 par
A_rho2}, the convex risk measure $\rho_{A,B}$ is finite; so its
subdifferential $\partial \rho_{A,B}(0)$ is non empty. More
precisely, there exists $\mathbf{Q}_0\in
\partial \rho_{A,B}(0)$ such that $\rho_{A,B}(X)\geq
\rho_{A,B}(0)+\mathbb{E}_{\mathbf{Q}_0}(-X).$ In other words,
$$\rho_{A,B}(0)\leq \rho_{A,B}(X)+\mathbb{E}_{\mathbf{Q}_0}(X)\leq
\rho_{A}(X-Y)-\mathbb{E}_{\mathbf{Q}_0}[-(X-Y)]+\rho_B(Y)-\mathbb{E}_{\mathbf{Q}_0}[-Y].$$
$ii-a)$ Assume  that $ \rho_A \square\rho_B(0)\geq c$. Applying
the previous inequality at $Y=-Z$, and $X=U+Y=U-Z$, we have
$$\rho_{A}(U)-\mathbb{E}_{\mathbf{Q}_0}[-U]\geq -\rho_B(-Z)-\mathbb{E}_{\mathbf{Q}_0}[Z]+\rho_{A,B}(0).$$
Then,
$$-\alpha_A(\mathbf{Q}_0):=\inf_U\{\rho_{A}(U)-\mathbb{E}_{\mathbf{Q}_0}[-U]\}
\geq
\alpha_B(\mathbf{Q}_0)+\rho_{A,B}(0):=\sup_Z\{-\rho_B(-Z)-\mathbb{E}_{\mathbf{Q}_0}[Z]+
\rho_{A,B}(0)\}.$$ By Theorem \ref{Theorem modification de rho1
par A_rho2} this inequality is in fact an equality. Picking
$r=\alpha_A(\mathbf{Q}_0)$, and defining
$a_{\mathbf{Q}_0}(X)=\mathbb{E}_{\mathbf{Q}_0}[-X]+r$ yield to an
affine function that separates $\rho_A$ and $-\rho_B(-.)+c$.\\
$ii-b)$ Finally, when
$\rho_A(\overline{X})+\rho_B(-\overline{X})=\rho_{A,B}(0)$, by the
above inequalities, we obtain
$-\rho_B(-\overline{X})-\mathbb{E}_{\mathbf{Q}_0}[-\overline{X}]\geq
-\rho_B(-Z)-\mathbb{E}_{\mathbf{Q}_0}[Z]$. In other words,
$\mathbf{Q}_0$ belongs to $\partial \rho_B(-\overline{X})$. By
symmetry, $\mathbf{Q}_0$ also belongs to $\partial
\rho_A(\overline{X})$.\vspace{1mm}\\
$iii)$
$\diamond$ The implication $(1)\> \Rightarrow \>(2) \>\Rightarrow \>(3)$ is clear, using the results i) and ii)
of this Theorem.\\
$\diamond$ Very naturally, one obtains $(3)\> \Rightarrow \>(2)$
and $(3)\> \Rightarrow \>(1)$ as the existence of $\mathbf{Q}_0\in
\mathbf{Dom}(\alpha_A)\cap\mathbf{Dom}(\alpha_B)$ implies that for
any $X$, $\rho_A(X) \geq
\mathbb{E}_{\mathbf{Q}_0}[-X]-\alpha_A(\mathbf{Q}) \qquad {\rm
and}\qquad \rho_B(-X) \geq
\mathbb{E}_{\mathbf{Q}_0}[X]-\alpha_B(\mathbf{Q})$. Considering
$r=\sup\{\alpha_A(\mathbf{Q}));\alpha_B(\mathbf{Q})\} $, one obtains
$(2)$. Moreover, $\rho_A(X)+\rho_B(-X) \geq
-(\alpha_A(\mathbf{Q})+\alpha_B(\mathbf{Q}))$ and taking the
infimum with respect to $X$,
$\rho_A \square \rho_B (0) >- \infty$, i.e. the property $(1)$.\\
$\diamond$ Let us now look at the following implication $(2)
\>\Rightarrow \>(4)$. We first observe that (2), i.e.,
$\rho_A(X)\geq -\mathbf{E}_{\mathbf{Q}_0}[-X]+r$ implies
$\rho^A_{0^+}(X)\geq -\mathbf{E}_{\mathbf{Q}_0}[-X]$, and
$\rho_B(-X)\geq \mathbf{E}_{\mathbf{Q}_0}[-X]-r$ implies
$\rho^B_{0^+}(-X)\geq \mathbf{E}_{\mathbf{Q}_0}[-X]$. Therefore,
we obtain (4) as
$\rho^A_{0^+}(X)+\rho^B_{0^+}(-X) \geq 0.$\\
$\diamond$ The converse implication $(4) \>\Rightarrow \>(2)$ is
obtained by applying the sandwich property \eqref{eq:sandwich} to
$\rho^A_{0^+}$ and $\rho^B_{0^+}$.\hspace{5mm}$\square$
\begin{remark}[On risk measures on
$\mathbb{L}_{\infty}(\mathbb{P})$] Let us consider the
inf-convolution between two risk measures $\rho_A$ and $\rho_B$,
where one of them, for instance $\rho_A$, is continuous from below
(and consequently from above) and therefore is defined on
$\mathbb{L}_{\infty}(\mathbb{P})$. In this case, as the
inf-convolution maintains the property of continuity from below
(see Theorem \ref{Theorem modification de rho1 par A_rho2}), the
risk measure $\rho_A \square \rho_B$ is also continuous from below
and therefore is a risk measure on
$\mathbb{L}_{\infty}(\mathbb{P})$, having a dual representation on
$\mathcal{M}_{1,ac}(\mathbb{P})$.
\end{remark}
\subsubsection{$\gamma$-Tolerant Risk Measures and Inf-Convolution \label{section static dilated
inf-conv}}
In this subsection, we come back to the particular class of
$\gamma$-tolerant convex risk measures $\rho_{\gamma}$ to give an
explicit solution to the exact inf-convolution. Recall that this
family of risk measures is generated from a root risk measure
$\rho$ by the following transformation $\rho _{\gamma
}(\xi_T)=\gamma\rho_\gamma\big(\frac{\xi_T}{\gamma}\big)$ where
$\gamma$ is the risk tolerance coefficient with respect to the
size of the exposure. These risk measures satisfy the following
semi-group property for the inf-convolution:
\begin{proposition}
\label{Theorem dilated risk measure} Let $\big( \rho _{\gamma
},\gamma
>0\big) $ be the family of $\gamma$-tolerant risk measures issued of $\rho$. Then, the
following properties hold:\\
$i)$ For any $\gamma_A,\gamma_B>0$,\hspace{3mm} $\rho
_{\gamma_A}\square \rho _{\gamma_B}=\rho _{\gamma_A +\gamma_B}$.\\
$ii)$ Moreover, $F^{*}=\frac{\gamma_B}{\gamma_A+\gamma_B}X$ is an optimal structure for the minimization program:%
\begin{equation*}
\rho _{\gamma_A+\gamma_B}( X) =\rho_{\gamma_A}\square
\rho_{\gamma_B}( X) =\inf_{F}\big\{ \rho _{\gamma_A}( X-F) +\rho
_{\gamma_B}( F) \big\} =\rho _{\gamma_A}\big( X-F^{*}\big) +\rho
_{\gamma_B}\big( F^{*}\big).
\end{equation*}
The inf-convolution is said to be exact at $F^*$.\\
$iii)$ Let $\rho $ and $\rho ^{\prime }$ be two convex risk
measures. Then, for any $\gamma >0$, $\rho _{\gamma }\square \rho
_{\gamma }^{\prime }=( \rho \square \rho ^{\prime }) _{\gamma
}$.\\
$iv)$ Assume $\rho(0)=0$ and $\rho'(0)=0$. When $\gamma=+\infty$,
this relationship still holds: $\rho _{\infty}\square \rho
_{\infty}^{\prime }=( \rho \square \rho ^{\prime }) _{\infty}$.\\
$v)$ If $\rho _{0^+}\square \rho _{0^+}^{\prime }(0)>-\infty$, we
also have $\rho _{0^+}\square \rho _{0^+}^{\prime }= ( \rho
\square \rho ^{\prime }) _{0^+}.$
\end{proposition}
\textbf{Proof:} Both $i)$ and $iii)$ are immediate consequences of
the
definition of infimal convolution.\\
$ii)$ We first study the stability property of the functional
$\rho_{\gamma}$ by studying the optimization program $
\rho_{\gamma_A}(X-F)+\rho_{\gamma_B}(F) \rightarrow \min_F$
restricted to the family $\{\alpha X, \alpha \in \mathbb{R}\}$.
Then, given the expression of the functional $\rho_{\gamma}$, a
natural candidate becomes
$F^{*}=\frac{\gamma_B}{\gamma_A+\gamma_B}X$, since
$$\rho_{\gamma
_{A}}\big(X-F^{*}\big) +\rho_{\gamma _{B}}\big( F^{*}\big)=(\gamma
_{A}+\gamma _{B})\rho\big(\frac{1}{\gamma _{A}+\gamma _{B}}
X\big)= \rho_{\gamma _{C}}(X).
$$
$iv)$ The asymptotic properties are based on the non increase of
the map $\gamma\rightarrow\rho_\gamma$. Then, when $\gamma$ goes
to infinity, pass to the limit is equivalent to take the infimum
w.r. of $\gamma$ and change the order of minimization,
in such way that pass to the limit is justified.\\
$v)$ When $\gamma$ goes to $0$, the problem becomes a minimax
problem, and we only obtain the inequality.\\
When the finite assumption holds, by Theorem \ref{Theorem
modification de rho1 par A_rho2}, the minimal penalty function of
$\rho _{0^+}\square \rho _{0^+}^{\prime }$ is $\alpha
_{0^+}+\alpha _{0^+}^{\prime }$. By the properties of conservative
risk measures, $\alpha _{0^+}$ is the convex indicator of
$\mathbf{Dom}(\alpha)$. So, $\alpha _{0^+}+\alpha _{0^+}^{\prime
}=l^{\mathbf{Dom}(\alpha)\cap \mathbf{Dom}(\alpha)'}$. On the
other hand, the minimal penalty function of $( \rho \square \rho
^{\prime }) _{0^+}$ is the indicator of
$\mathbf{Dom}(\alpha+\alpha')$. Since, $\alpha$ is dominated by
the same minimal bound $-\rho(0)$,
$\mathbf{Dom}(\alpha+\alpha')=\mathbf{Dom}(\alpha)\cap\mathbf{Dom}(\alpha')$
Both risk measures have same minimal penalty functions. This
completes the proof. \hspace{5mm}$\square$
\subsubsection{An Example of Inf-Convolution: the Market Modified Risk Measure \label{Subsection
Market modified r.m.}}
We now consider a particular inf-convolution which is closely
related to Subsection \ref{section risk measures and hedging} as
it also deals with the question of optimal hedging.\\
More precisely, the following minimization problem
$$\inf_{H \in \mathcal{V}_T} \rho(X-H)$$ can be seen as an hedging
problem, where $\mathcal{V}_T$ corresponds to the set of hedging
instruments. It somehow consists of restricting the risk measure
$\rho$ to a particular set of admissible variables and is in fact the inf-convolution $\rho
\square \nu ^{ \mathcal{V}_{T}}$. Using Proposition \ref{proposition regularisation par worst}, it can also be seen as the
inf-convolution $\rho \square l^{ \mathcal{V}_{T}} \square
\rho_{\rm worst}$. The main role of $\rho_{\rm worst}$ is to transform the convex indicator
$l^{\mathcal{V}_T}$, which is not a convex risk measure (in
particular, it is not translation invariant), into the
convex risk measure $\nu^{\mathcal{V}_T}$.\vspace{2mm}\\
The following corollary is an immediate extension of Theorem \ref{Theorem modification de rho1 par A_rho2}
as it establishes that the value functional of
the problem, denoted by $\rho^m$, is a convex risk measure,
called {\it market modified risk measure}.
\begin{corollary}
\label{corollary modification de rho par H} Let $\mathcal{V}_{T}$
be a convex subset of $\mathbb{L}_{\infty}(\mathbb{P})$ and $\rho$ be a convex risk measure with
penalty function $\alpha $ such
that $ \inf \big\{ \rho( -\xi_T) ,\xi_T\in \mathcal{V}_{T}\big\}
>-\infty$. The inf-convolution of $\rho $ and $\nu ^{\mathcal{V_T}}$,
\index[not]{$\rho^{m}$} $\rho ^{m}\equiv \rho \square \nu ^{\mathcal{V}_{T}}$, also defined as
\begin{equation}\label{inf-convolution deux mesures de risque}
\rho ^{m}(\Psi) \equiv \inf \big\{ \rho (\Psi-\xi_T)\,\big|\,\xi_T
\in \mathcal{V}_T \big\} = \rho \square l^{\mathcal{V}_T} (\Psi)
\end{equation}
is a convex risk measure, called \index[sub]{market modified risk
measure}market modified risk measure, with minimal penalty
function defined on $\mathbf{M}_{1,ac}(\mathbb{P})$,
$\alpha^m(\mathbf{Q})=\alpha(\mathbf{Q})+\alpha^{{\mathcal{V_T}}}(\mathbf{Q})$.\\
Moreover, if $\rho$ is continuous from below, $\rho^m$ is also continuous from below.
\end{corollary}
This corollary makes precise the direct impact on the risk measure
of the agent of the opportunity to invest optimally in a financial
market.
\begin{remark}
Note that the set $\mathcal{V}_T$ is rather general. In most
cases, additional assumptions will be added and the framework will
be similar to those described in Subsection \ref{subsection
hedging description}.
\end{remark}
{\bf Acceptability and market modified risk measure:} The market
modified risk measure has to be related to the notion of
acceptability introduced by \index[aut]{Carr}Carr,
\index[aut]{Geman}Geman and \index[aut]{Madan}Madan in
\cite{CarrGemanMadan}. In this paper, they relax the strict notion
of hedging in the following way: instead of imposing that the
final outcome of an acceptable position, suitably hedged, should
always be non-negative, they simply require that it remains
greater than an acceptable position. More precisely, using the
same notations as in Subsection \ref{subsection hedging
description} and denoting by $\mathcal{A}$ a given acceptance set
and by $\rho_{\mathcal{A}}$ its related risk measure, we can
define the convex risk measure:
$$\bar{\nu}^{\mathcal{H}}(X)=\inf \big\{m \in \mathbb{R},\>\> \exists \theta \in \mathcal{K}\>\> \exists A \in \mathcal{A}\>:\>
m+X+G(\theta)\geq A \>\> \mathbb{P}\>a.s. \big\}$$
To have a clearer picture of what this risk measure really is, let us first fix $G(\theta)$.
In this case, we simply look at $\rho_{\mathcal{A}}(X+G(\theta))$. Then, the risk measure $\bar{\nu}^{\mathcal{H}}$ is
defined by taking the infimum of $\rho_{\mathcal{A}}(X+G(\theta))$ with respect to $\theta$,
$$\bar{\nu}^{\mathcal{H}}(X)= \inf_{\theta \in \mathcal{K}}\rho_{\mathcal{A}}(X+G(\theta))
=\inf_{H \in \mathcal{H}}\rho_{\mathcal{A}}(X-H)$$
Therefore, the risk
measure $\bar{\nu}^{\mathcal{H}}$ is in fact the particular market modified risk measure
$\rho^m= \nu^{\mathcal{H}} \square \rho_{\mathcal{A}}$. We obtain
directly the following result of
\index[aut]{F\"{o}llmer}F\"{o}llmer and \index[aut]{Schied}Schied
\cite{FoellmerSchied02b} (Proposition 4.98): the minimal penalty
function of this convex risk measure $\bar{\nu}^{\mathcal{H}}$ is
given by
$$\bar{\alpha}^{\mathcal{H}}(\mathbf{Q})=\alpha^{\mathcal{H}}(\mathbf{Q})+\alpha(\mathbf{Q})$$
where $\alpha^{\mathcal{H}}$ is the minimal penalty function of
$\nu^{\mathcal{H}}$ and $\alpha$ is the minimal penalty function
of the convex risk measure with acceptance set $\mathcal{A}$.
\section{Optimal Derivative Design\label{Section
utility exp}}
In this section, we now present our main problem, that of
derivative optimal design (and pricing)\index[sub]{optimal
design}. The framework we generally consider involves two economic
agents, at least one of them being exposed to a non-tradable risk.
The risk transfer \index[sub]{risk transfer} between both agents
takes place through a structured contract denoted by $F$ for an
initial price $\pi$. The problem is therefore to design the
transaction, in other words, to find the structure $F$ and its
price $\pi$. This transaction may occur only if both agents find
some interest in doing this transaction. They express their
satisfaction or interest in terms of the expected utility of their
terminal wealth after the transaction, or more generally in terms
of risk measures.
\subsection{General Modelling}
\subsubsection{Framework \label{framework simple}}
Two economic agents, respectively denoted by $A$ and $B$, are
evolving in an uncertain universe modelled by a standard
measurable space $( \Omega ,\Im)$ or, if a reference probability
measure is given, by a probability space $( \Omega ,\Im,
\mathbb{P})$. In the following, for the sake of simplicity in our
argumentation, we will make no distinction between both
situations. More precisely, in the second case, all properties
should hold
$\mathbb{P}-a.s.$.\\
Both agents are taking part in trade talks to improve the
distribution and management of their own risk. The nature of both
agents can be quite freely chosen. It is possible to look at them
in terms of a classical insured-insurer relationship, but from a
more financial point of view, we may think of agent A as a market
maker or a trader managing a particular book and of agent B as a
traditional investor
or as another trader.\\
More precisely, we assume that at a future time horizon $T$, the
value of agent A's terminal wealth, denoted by $X_T^{A}$, is
sensitive to a non-tradable risk. Agent B may also have her own
exposure $X_T^{B}$ at time $T$. Note that by "terminal wealth", we
mean the terminal value at the time horizon $T$ of all capitalized
cash flows paid or received between the initial time and $T$; no
particular sign constraint is imposed. Agent A wants to issue a
structured contract (financial derivative, insurance contract...)
$F$ with maturity $T$ and forward price $\pi$ to reduce her
exposure $X^A_T$. Therefore, she calls on agent B. Hence, when a
transaction occurs, the terminal wealth of the agent A and B are
$$
W^A_T=X_T^{A}-F+\pi,\qquad W^B_T=X_T^{B}+F-\pi.
$$
As before, we assume that all the quantities we consider belong to
the Banach space $\mathcal{X}$, or, if a reference probability
measure is given, to $\mathbb{L}_{\infty}(\mathbb{P})$.\vspace{1mm}\\
The problem is therefore to find the optimal structure of the risk
transfer $(F,\pi)$ according to a given choice criterion, which is
in our study a convex risk measure. More precisely, assuming that
agent A (resp. agent B) assesses her risk exposure using a convex
risk measure $\rho_A$ (resp. $\rho_B$), agent $A$'s objective is
to choose the optimal structure $\big( F,\pi \big) $ in order to
minimize the risk measure of her final wealth
$$\rho_A(X_T^{A}-F+\pi)\rightarrow \inf_{F\in \mathcal{X},\pi }.$$ Her constraint
is then to find a counterpart. Hence, agent $B$ should have an
interest in doing this transaction. At least, the $F$-structure
should not worsen her risk measure. Consequently, agent $B$ simply
compares the risk measures of two terminal wealth, the first one
corresponds to the case of her initial exposure $X_T^{B}$ and the
second one to her new wealth if she enters the $F$-transaction,
$$
\rho_B(X_T^{B}+F-\pi)\leq \rho_B(X_T^{B}).
$$
\subsubsection{Transaction Feasibility and Optimization Program}
\index[sub]{transaction feasibility} The optimization program as
described above as
\begin{equation}
\inf_{F\in \mathcal{X},\pi}\rho_A(X_T^{A}-F+\pi) \qquad
\mbox{subject to} \quad \rho_B(X_T^{B}+F-\pi)\leq \rho_B(X_T^{B})
\label{optimization pg with constraint}
\end{equation}
can be simplified using the cash translation invariance property.
More precisely, binding the constraint imposed by agent B at the
optimum and using the translation invariance property of $\rho_B$,
we find directly the {\it optimal pricing rule} for a structure
$F$:
\begin{equation}
\pi_B( F) =\rho_B\big( X_T^{B}\big)-\rho_B\big(X_T^{B}+F\big).
\label{optimal pricing rule}
\end{equation}
This pricing rule is an indifference pricing rule for agent B. It
gives for any structure $F$ the maximum amount agent B is ready to
pay in order to enter the transaction.\vspace{2mm}\\
Note also that this optimal pricing rule together with the cash
translation invariance property of the functional $\rho_A$ enable
us to rewrite the optimization program (\ref{optimization pg with
constraint}) as follows, without any need for a Lagrangian
multiplier:
\begin{equation*}
\inf_{F\in \mathcal{X}}\big\{\rho_A \big(X_T^{A}-F\big)+
\rho_B\big(X_T^{B}+F\big)-\rho_B\big( X_T^{B}\big)\big\}
\end{equation*}
or to within the constant $\rho_B( X_T^{B})$ as:
\begin{equation}
R_{AB}(X_T^A,X_T^B)=\inf_{F\in
\mathcal{X}}\{\rho_A\big(X_T^{A}-F\big)+
\rho_B\big(X_T^{B}+F\big)\}. \label{program reduit}
\end{equation}
\paragraph{Interpretation in Terms of Indifference Prices}
This optimization program (Program (\ref{program reduit})) can be
reinterpreted in terms of the indifference prices, using the
notations introduced in the exponential utility framework in
Subsection \ref{subsection transaction prices}. To show this, we
introduce the constants $\rho_A\big(X_T^{A}\big)$ and
$\rho_B\big(X_T^{B}\big)$ in such a way that Program (\ref{program
reduit}) is equivalent to:
$$\inf_{F\in \mathcal{X}}\big\{\rho_A\big(X_T^{A}-F\big)-\rho_A\big(X_T^{A}\big)+\rho_B\big(X_T^{B}+F\big)
-\rho_B\big(X_T^{B}\big)\}.$$ Then, using the previous comments,
it is possible to interpret
$\rho_A\big(X_T^{A}-F\big)-\rho_A\big(X_T^{A}\big)$ as
$\pi_A^{s}\big(F|X_T^{A}\big)$, i.e. the seller's indifference
pricing rule for $F$ given agent A's initial exposure $X_T^{A}$,
while $\rho_B\big(X_T^{B}+F\big)-\rho_B\big(X_T^{B}\big)$ is
simply the opposite of $\pi_B^{b}(F|X_T^{B})$, the buyer's
indifference pricing rule for $F$ given agent B's initial exposure
$X_T^{B}$. For agent A, everything consists then of choosing the
structure as to minimize the difference between her (seller's)
indifference price (given $X^A_T$) and the (buyer's) indifference
price imposed by agent B:
\begin{equation}
\inf_{F\in \mathcal{X}}\big\{\pi_A^{s}\big(F|X_T^{A}\big) -
\pi_B^{b}\big(F|X_T^{B}\big)\big\}\leq 0.
\end{equation}
Note that for $F \equiv 0$, the spread between both transaction
indifference prices is equal to $0$. Hence, the infimum is always
non-positive. This is completely coherent with the idea that the
optimal transaction obviously reduces the risk of agent $A$. The
transaction may occur since the minimal
seller price is less than the maximal buyer price.\\
For agent $A$, everything can also be expressed as the following
maximization program
\begin{equation}
\sup_{F\in
\mathcal{X}}\big\{\pi_B^{b}\big(F|X_T^{B}\big)-\pi_A^{s}\big(F|X_T^{A}\big)\big\}.
\end{equation}
The interpretation becomes then more obvious since the issuer has
to optimally choose the structure in order to maximize the
"ask-bid" spread associated with transaction.
\paragraph{Relationships with the Insurance Literature and the Principal-Agent Problem}
\index[sub]{principal-agent problem} The relationship between both
agents is very similar to a Principal-Agent framework. Agent A
plays an active role in the transaction. She chooses the "payment
structure" and then is the "Principal" in our framework. Agent B,
on the other hand, is the "Agent" as she simply imposes a price
constraint to the Principal and in
this sense is rather passive.\vspace{1mm}\\
Such a modelling framework is also very similar to an insurance
problem: Agent A is looking for an optimal "insurance" policy to
cover her risk (extending here the simple notion of loss as
previously mentioned). In this sense, she can be seen as the
"insured". On the other hand, Agent B accepts to bear some risk.
She plays the same role as an "insurer" for Agent A. In fact, this
optimal risk transfer problem is closely related to the standard
issue of optimal policy design in insurance, which has been widely
studied in the literature (see for instance \index[aut]{Borch}
Borch \cite{Borch}, \index[aut]{B\"{u}hlman} B\"{u}hlman
\cite{Buhlmann70}, \cite{Buhlmann80} and \cite{Buhlmann84},
\index[aut]{B\"{u}hlman} B\"{u}hlman and \index[aut]{Jewell}
Jewell \cite{BuhlmannJewell}, \index[aut]{Gerber} Gerber
\cite{Gerber}, \index[aut]{Raviv} Raviv \cite{Raviv}). One of the
fundamental characteristics of an insurance policy design problem
is the sign constraint imposed on the risk, that should represent
a loss. Other specifications can be mentioned as moral hazard or
adverse selection problems that have to be taken into account when
designing a policy (for more details, among a wide literature,
refer for instance to the two papers on the relation
Principal-Agent by \index[aut]{Rees} Rees \cite{Rees1} and
\cite{Rees2}). These are related to the potential influence of the
insured on the considered risk.\\
Transferring risk in finance is somehow different. Risk is then
taken in a wider sense as it represents the uncertain outcome. The
sign of the realization does not a priori matter in the design of
the transfer. The derivative market is a good illustration of this
aspect: forwards, options, swaps have particular payoffs which are
not directly related to any particular loss of the contract's
seller.
\subsection{Optimal Transaction}
This subsection aims at solving explicitly the optimization
Program (\ref{program reduit}):
\begin{equation*}
R_{AB}(X_T^A,X_T^B)=\inf_{F\in
\mathcal{X}}\{\rho_A\big(X_T^{A}-F\big)+
\rho_B\big(X_T^{B}+F\big)\}.
\end{equation*}
The value functional $R_{AB}(X_T^{A},X_T^{B})$ can be seen as the
{\it residual risk measure} after the $F$-transaction, or
equivalently as a measure of the risk remaining after the
transaction. It obviously depends on both initial exposures
$X_T^{A}$ and $X_T^{B}$ since the transaction consists of an
optimal redistribution of the respective risk of both agents.\\
Let us denote by $\widetilde{F} \equiv X_T^{B}+F \in \mathcal{X}$.
The program to be solved becomes
\begin{equation*}
R_{AB}\big(X_T^{A},X_T^{B}\big)= \inf_{\widetilde{F} \in
\mathcal{X}}\{\rho_{A}\big(X_T^{A}+X_T^{B}-\widetilde{F}\big)
+\rho_{B}\big( \widetilde{F}\big)\},
\end{equation*}
or equivalently, using Section \ref{Subsection static
inf-convolution}, it can be written as the following
inf-convolution problem
\begin{equation}
R_{AB}\big(X_T^{A},X_T^{B}\big)=\rho_A \square \rho_B
(X_T^{A}+X_T^{B}). \label{inf-convolution pg}
\end{equation}
As previously mentioned in Theorem \ref{Theorem modification de
rho1 par A_rho2}, the condition $\rho_A \square \rho_B(0)
>-\infty$ is required when considering this inf-convolution
problem. This condition is equivalent to $\forall \xi \in
\mathcal{X},\>\rho^A_{0^+}(\xi)+\rho^B_{0^+}(-\xi) \geq 0$ (Theorem \ref{Sandwich theorem} iii)). \\
This property has a nice economic interpretation, since it says
that the inf-convolution program makes sense if and only if {\em
for any derivative $\xi$, the conservative seller price of the
agent $A$, $-\rho^A_{0^+}(\xi)$, is less than the conservative
buyer price of the agent B, $\rho^B_{0^+}(-\xi)$.} \vspace{2mm}\\
In the following, we assume such a condition to be satisfied.
The problem is not to study the residual risk measure as
previously but to characterize the optimal structure $\tilde{F}^*$
or $F^*$ such that the inf-convolution is exact at this point.\\
To do so, we first consider a particular framework where the
optimal transaction can be explicitly identified. This corresponds
to a well-studied situation in economics where both agents belong
to the same family.
\subsubsection{Optimal Transaction between Agents with Risk Measures in the Same
Family}
More precisely, we now assume that both agents have $\gamma$-tolerant risk measures
$\rho_{\gamma_A}$ and $\rho_{\gamma_B}$ from the same root risk measure $\rho$ with risk
tolerance coefficients $\gamma_A$ and $\gamma_B$, as introduced in Subsection
\ref{subsection dilated static risk measure}. In this framework, the optimization program
(\ref{inf-convolution pg}) is written as follows:
\begin{equation*}
R_{AB}\big(X_T^{A},X_T^{B}\big)=\rho_{\gamma_A} \square
\rho_{\gamma_B} (X_T^{A}+X_T^{B})
\end{equation*}
In this framework, the optimal risk transfer is consistent with
the so-called Borch's theorem. In this sense, the following result
can be seen as an extension of this theorem since the framework we
consider here is different from that of utility functions. In his
paper \cite{Borch}, \index[aut]{Borch} Borch obtained indeed, in a
utility framework, optimal exchange of risk, leading in many cases
to familiar linear quota-sharing of total pooled losses.
\begin{theorem}[Borch \cite{Borch}]
\label{Cas exp: proposition Borch} The residual risk measure
\index[sub]{residual risk measure} after the transaction is given
by:
$$R_{AB}\big(X_T^{A},X_T^{B}\big)=\inf_{F\in \mathcal{X}}\big\{\rho_{\gamma _{A}}\big(X_T^{A}-F\big)+\rho_{\gamma
_{B}}\big(X_T^{B}+F\big)\big\} = \rho_{\gamma
_{C}}(X_T^{A}+X_T^{B}) \quad \mbox{with} \quad
\gamma_C=\gamma_A+\gamma_B.$$ The optimal structure is given as a
proportion of the initial exposures $X_T^{A}$ and $X_T^{B}$,
depending only on the risk tolerance coefficients of both agents:
\begin{equation}
\label{cas exp: optimal structure cas simple}
 F^{*}=\frac{\gamma
_{B}}{\gamma _{A}+\gamma _{B}}X_T^{A}-\frac{\gamma _{A}}{\gamma
_{A}+\gamma _{B}}X_T^{B}\qquad \text{(to within a constant).}
\end{equation}
\end{theorem}
The equality in the equation \eqref{cas exp: optimal structure cas
simple} has to be understood $\mathbb{P}\> a.s.$ if the space of
structured products is
$\mathbb{L}_{\infty}(\mathbb{P})$.\vspace{2mm}\\
\textbf{Proof:}
The optimization program (\ref{program reduit}) to be solved (with
$\tilde{F} \equiv X_T^{B}+F \in \mathcal{X}$) is
\begin{equation*}
R_{AB}\big(X_T^{A},X_T^{B}\big)=\inf_{\tilde{F} \in
\mathcal{X}}\big(\rho_{\gamma
_{A}}\big(X_T^{A}+X_T^{B}-\widetilde{F}\big) +\rho_{\gamma
_{B}}\big( \widetilde{F}\big) \big).
\end{equation*}
Using Proposition \ref{Theorem dilated risk measure}, the optimal
structure $\tilde{F}^{*}$ is
$\widetilde{F}^{*}=\frac{\gamma_B}{\gamma_A+\gamma_B}(X_T^{A}+X_T^{B})$.
The
result is then obtained by replacing $\widetilde{F}^{*}$ by $F^{*}-X_T^{B}$. $\square$\vspace{2mm}\\
{\bf Comments and properties:} $i)$ Both agents are transferring a
part of their initial risk according to their relative tolerance.
The optimal risk transfer underlines the symmetry of the framework
for both agents. Moreover, even if the issuer, agent A has no
exposure, a transaction will occur between both agents. The
structure $F$ enables them to exchange a part of their respective
risk. Note that if none of the agents is initially exposed, no
transaction will occur. In this sense, the transaction has a
non-speculative
underlying logic.\\
$ii)$ Note also that the composite parameter $\gamma _{C}$ is
simply equal to the sum of both risk tolerance coefficients
$\gamma _{A}$ and $\gamma _{B}$. This may justify the use of risk
tolerance instead of risk aversion where harmonic mean has to be
used.
\subsubsection{Individual Hedging as a Risk Transfer \label{static individual hedging}}
In this subsection, we now focus on the individual hedging
problem of agent A and see how this problem can be
interpreted as a particular risk transfer problem. The question of
optimal hedging has been widely studied in the literature under
the name of ''hedging in incomplete markets and pricing via
utility maximization'' in some particular framework. Most of the
studies have considered exponential utility functions. Among the
numerous papers, we may quote the papers by \index[aut]{Frittelli}
Frittelli \cite{Frittelli00a}, \index[aut]{El Karoui}El Karoui and
\index[aut]{Rouge}Rouge \cite{ElKarouiRouge},
\index[aut]{Delbaen}\index[aut]{Grandits}\index[aut]{Rheinlander}
\index[aut]{Samperi}\index[aut]{Schweizer}\index[aut]{Stricker}Delbaen
et al. \cite{Delbaenetal}, \index[aut]{Kabanov}Kabanov and
\index[aut]{Stricker}Stricker \cite{KabanovStricker}
or \index[aut]{Becherer}Becherer \cite{Becherer}.\vspace{2mm}\\
We assume that agent A assesses her risk using a
($\mathbb{L}_{\infty}(\mathbb{P})$) risk measure $\rho_A$. She can
(partially) hedge her initial exposure $X$ using instruments from
a convex subset $\mathcal{V}^{A} _{T}$ (of
$\mathbb{L}_{\infty}(\mathbb{P})$). Her objective is to minimize
the risk measure of her terminal wealth.
\begin{equation}
\inf_{\xi \in \mathcal{V}^{A}_{T}}\rho_A\big(X_T^{A}-\xi \big).
\label{hedging pb agent A}
\end{equation}
The $\mathbb{L}_{\infty}(\mathbb{P})$ framework has been carefully
described in Subsection \ref{subsection hedging description}. In
particular, to have coherent transaction prices, we assume
in the following that the market is arbitrage-free.\\
As already mentioned in Subsection \ref{Subsection Market modified
r.m.}, the opportunity to invest optimally in a financial market
has a direct impact on the risk measure of the agent and
transforms her initial risk measure $\rho_A$ into the market
modified risk measure $\rho_A^m=\rho_A \square \nu ^{A}$.\\
This inf-convolution problem makes sense if the condition
$\rho_A^m(0)>-\infty$ is satisfied. The hedging problem of agent
$A$ is identical to the previous risk transfer problem
(\ref{program reduit}), agent $B$ being now the financial market
with the associated risk measure $\nu^{A}$.\vspace{-2mm}
\paragraph{Existence of an Optimal Hedge}
The question of the existence of an optimal hedge can be answered using different approaches.
One of them is based on analysis techniques and we present it in this subsection. In the following, however,
when introducing dynamic risk measures, we will consider other methods leading to a more constructive answer.\vspace{1mm}\\
In this subsection, we are interested in studying the existence of
a solution for the hedging problem of agent A (Program
(\ref{hedging pb agent A})) or equivalently for the
inf-convolution problem in $\mathbb{L}_{\infty}(\mathbb{P})$. The
following of existence can be obtained:
\begin{theorem}
\label{hedging pb: existence theorem} Let $\mathcal{V}_{T}$ be a
convex subset of $\mathbb{L}_{\infty}(\mathbb{P})$ and $\rho$ be a
convex risk measure on $\mathbb{L}_{\infty}(\mathbb{P})$
continuous from below, such that $\inf_{\xi \in \mathcal{V}
_{T}}\rho(-\xi) >-\infty $.\\
Assume the convex set $\mathcal{V}_{T}$ bounded in $\mathbb{L}^{\infty}( \mathbb{P%
}) $. The infimum of the hedging program
$$\rho^{m}\big(X\big) \triangleq \inf_{\xi \in
\mathcal{V}_{T}}\rho\big(X-\xi \big)$$
is ''attained'' for a random variable $\xi ^{*}_T$ in $\mathbb{L}^{\infty}( \mathbb{P}%
) $, belonging to the closure of $\mathcal{V}_{T}$ with respect to
the a.s. convergence.
\end{theorem}
\textbf{Proof:} First note that the proof of this theorem
relies on arguments similar to those used by \index[aut]{Kabanov}
Kabanov and \index[aut]{Stricker} Stricker \cite{KabanovStricker}.
In particular, a key argument is the Komlos Theorem
(\index[aut]{Komlos}Komlos \cite{Komlos}):
\begin{lemma}[Komlos]
\label{Komlos theorem}Let $( \phi _{n}) $ be a sequence in $%
\mathbb{L}^{1}( \mathbb{P}) $ such that
$\sup_{n}\mathbb{E}_{\mathbb{P}}( | \phi _{n}| ) <+\infty$.
Then there exists a subsequence $( \phi _{n^{\prime }}) $ of $%
( \phi _{n}) $ and a function $\phi ^{*}\in L^{1}( \mathbb{P}%
) $ such that for every further subsequence $\big( \phi
_{n^{\prime \prime }}\big) $ of $( \phi _{n}) $,
the Cesaro-means of these subsequences converge to $\varphi^*$, that is
$$\lim_{N\rightarrow \infty }\frac{1}{N}\sum\limits_{n^{\prime
\prime }=1}^{N}\phi _{n^{\prime \prime }}( \omega) =\phi ^{*}(
\omega ) \text{ \qquad for almost every }\omega \in \Omega. $$
\end{lemma}
We first show that the set $S_r=\{\xi \in \mathbb{L}_{\infty}(
\mathbb{P})\big|\> \rho(X-\xi)\leq r\}$ is closed for the
weak*-topology. To do that, by the Krein-Smulian theorem
(\cite{FoellmerSchied02b} Theorem A.63), it is sufficient to show
that $S_r\cap \{\xi; \|\xi\|_{\infty}\leq C\}$ is closed in
$\mathbb{L}^{\infty}( \mathbb{P})$.\\
Let $( \xi _{n}\in \mathcal{V}_{T}) $ be a sequence bounded by $C$,
converging in $\mathbb{L}^{\infty}$ to $\xi^*$. A subsequence
still denoted by $\xi_n$ converges a.s. to $\xi^*$. Since
$\rho_{A} $ is continuous from below, $\rho$ is continuous w.r. to
pointwise convergence of bounded sequences and then $\xi^*$
belongs to $S_r$. $S_r$ is weak*-closed.\\
Given the assumption that $( \xi _{n}) $ is $\mathbb{L}^{\infty}$%
-bounded, we can apply Komlos lemma:
 therefore, there
exists a subsequence $( \xi _{j_{k}}\in
\mathcal{V}^{A}_{T}) $ such that the Cesaro-means, $\widetilde{\xi }%
_{n}\triangleq \frac{1}{n}\sum\limits_{k=1}^{n}\xi _{j_{k}}$
converges almost surely to $\xi ^{*}\in \mathbb{L}^{\infty}(
\mathbb{P}) $. Note that $\widetilde{\xi }_{n}$ belongs to
$\mathcal{V}^{A}_{T}$ as a convex combination of elements of
$\mathcal{V}^{A}_{T}$. So $\xi ^{*}$ belongs to the $a.s.$ closure
of $\mathcal{V}^{A}_{T}$. Since $\rho_{A} $ is continuous from
below, $\rho$ is continuous w.r. to pointwise convergence of
bounded sequences.
\begin{equation*}
\lim_{n}\sup \rho_{A} \big( X-\widetilde{\xi }_{n}\big) \leq
\rho_{A} \big( X-\xi ^{*}\big) =\rho_{A} \big( \lim_{n}\big(
X-\widetilde{\xi }_{n}\big) \big) \leq \lim_{n}\inf \rho_{A} \big(
X-\widetilde{\xi }_{n}\big).
\end{equation*}
Then, $\rho_{A} ^{m}( X) \leq \rho_{A}( X-\xi
^{*}) \leq \lim_{n}\inf \rho_{A} \big(
\frac{1}{n}\sum\limits_{k=1}^{n}( X-\xi
_{j_{k}}) \big) \leq \lim_{n}\inf \frac{1}{n}\sum\limits_{k=1}^{n}%
\rho_{A}( X-\xi _{j_{k}})$ by Jensen inequality.
Finally, given the convergence of $\rho_{A} \big( X-\xi
_{j_{k}}\big) $ to $\rho_{A} ^{m}( X) $, the Cesaro-means also
converge and $\rho_{A}( X-\xi ^{*}) =\inf_{\xi \in \mathcal{V}^{A}_{T}}\rho_{A}
( X-\xi)$. $\qquad \square$
\subsubsection{$\gamma$-Tolerant Risk Measures: Derivatives Design with Hedging Opportunities\label{section Borch
dilate}}
We now consider the situation where both agents $A$ and $B$ have a
$\gamma$-dilated risk measure, defined on
$\mathbb{L}_{\infty}(\mathbb{P})$ and continuous from above.
Moreover, they may reduce their risk by transferring it between
themselves but also by investing in the
financial market, choosing optimally their financial investments.\\
The investment opportunities of both agents are described by two
convex subsets \index[not]{$\mathcal{V}^{A}_{T}$}
$\mathcal{V}^{A}_{T}$ and \index[not]{$\mathcal{V}_{T}^{B}$}
$\mathcal{V}^{B}_{T}$ of $\mathbb{L}_{\infty}(\mathbb{P})$. In
order to have coherent transaction prices, we assume that the
market is arbitrage-free. In our framework, this can be expressed
as the existence of a probability measure which is {\bf
equivalent} to $\mathbb{P}$ in both sets of probability measures
$\mathcal{M}_{\mathcal{V}^{i}_T}=\big\{\mathbb{Q} \in
\mathcal{M}_{1,e}(\mathbb{P}); \forall \xi \in \mathcal{V}^{i}_T,
\mathbb{E}_{\mathbb{Q}}[-\xi] \leq 0\big\}$ for $i=A,B$.
Equivalently,
$$\exists \mathbb{Q} \sim \mathbb{P} \quad \mbox{s.t.}\quad \mathbb{Q} \in
\mathcal{M}_{\mathcal{V}^{A}_T} \cap
\mathcal{M}_{\mathcal{V}^{B}_T}.$$ This opportunity to invest
optimally in a financial market reduces the risk of both agents.
To assess their respective risk exposure, they now refer to market
modified risk measures $\rho _{\gamma_A}^{m}$ and $\rho
_{\gamma_B}^{m}$ defined if $J=A,B$ as $$ \rho _{\gamma_A}^{m}(
\Psi) =\rho _{\gamma_A}\square \nu ^{A}( \Psi ) \qquad
\text{and\qquad }\rho _{\gamma_B}^{m}(\Psi) =\rho
_{\gamma_B}\square \nu ^{B}( \Psi ).
$$
Let us consider directly the optimal risk transfer problem with these market modified risk measures, i.e.
\begin{equation}\label{pg general avec Fbar et acceptance set}
R_{AB}^{m}\big(X_T^A+X_T^B\big) = \inf_{F\in
\mathcal{X}}\big\{\rho _{\gamma_A}^{m}\big( X_T^A-F\big) +\rho
_{\gamma_B}^{m}\big(X_T^B+F\big) \big\}
\end{equation}
The details of this computation will be given in the next
subsection, when considering the general framework. The residual
risk measure $R_{AB}^{m}\big(X_T^A+X_T^B\big) $ defined in
equation (\ref{pg general avec Fbar et acceptance set}) may be
simplified using the commutativity property of the inf-convolution
and the semi-group property of $\gamma$-tolerant risk measures:
\begin{eqnarray*}
R_{AB}^{m}\big(X_T^A+X_T^B\big) &=& \rho_{\gamma_A}^m \square
\rho_{\gamma_B}^m
\big(X_T^A+X_T^B\big)\\
&=& \rho _{\gamma_A}\square \nu ^{A} \square \rho
_{\gamma_B}\square \nu ^{B}\big(X_T^A+X_T^B\big)\\
&=&\rho _{\gamma_A}\square \rho
_{\gamma_B}\square \nu ^{A} \square \nu ^{B}\big(X_T^A+X_T^B\big)\\
&=& \rho _{\gamma_C}\square \nu ^{A}\square
\nu^{B}\big(X_T^A+X_T^B\big).
\end{eqnarray*}
where $\rho _{\gamma_C}$ is the $\gamma$-tolerant risk measure
associated with the risk tolerance coefficient $\gamma _{C}=\gamma
_{A}+\gamma
_{B}$.\vspace{2mm}\\
This inf-convolution program makes sense under the initial
condition $\rho _{\gamma_A}^{m}\square \rho_{\gamma_B}^{m}(0) >
-\infty$. Such an assumption is made. The following theorem gives
the optimal risk transfer in different situations depending on the
access both agents have to the financial markets.
\begin{theorem}
\label{theorem Borch dilate}
Let both agents have $\gamma$-tolerant risk measures with respective
risk tolerance coefficients $\gamma_A$ and $\gamma_B$.\\
(a) If both agents have the same access to the financial market
from a {\bf cone}, $\mathcal{V}_T$, then an optimal structure,
solution of the minimization Program $\eqref{pg general avec Fbar
et acceptance set}$ is given by:
\begin{equation*}
F^{*}=\frac{\gamma _{B}}{\gamma _{A}+\gamma _{B}}X_T^A -
\frac{\gamma_A}{\gamma_A+\gamma_B}X_T^B.
\end{equation*}
(b) Assume that both agents have different access to the financial
market via {\bf two convex sets} $\mathcal{V}^{A}_T$ and
$\mathcal{V}^{B}_T$. Suppose $\xi ^{*}=\eta _{A}^{*}+\eta
_{B}^{*}$ is an optimal solution of the Program $ \inf_{\xi \in
{\mathcal{V}_{T}^{(A+B)} }}\rho _{\gamma_C}\big(X_T^A+X_T^B-\xi
\big) $ with $\eta _{A}^{*}\in \mathcal{V}_{T}^{( A) }$, $\eta
_{B}^{*}\in \mathcal{V}_{T}^{( B) }$ and $\mathcal{V}_{T}^{(A+B)}
= \big \{\xi^A_T+\xi^B_T\>| \>\xi^A_T \in \mathcal{V}_{T}^{(A)},\>
\xi^B_T \in \mathcal{V}_{T}^{(B)}\big\}$. Then
\begin{equation*}
F^{*}=\frac{\gamma _{B}}{\gamma _{A}+\gamma _{B}}X_T^A - \frac{\gamma_A}{\gamma_A+\gamma_B}X_T^B
-\frac{\gamma _{B}}{\gamma
_{A}+\gamma _{B}}\eta _{A}^{*}+\frac{\gamma _{A}}{\gamma _{A}+\gamma _{B}}%
\eta _{B}^{*}
\end{equation*}
is an optimal structure.
Moreover,\\
$i)$ \noindent $\eta _{B}^{*}$ is an optimal hedging portfolio of
$\big(
X_T^B+F^{*}\big) $ for Agent $%
B$%
\begin{equation*}
\frac{1}{\gamma _{B}}\rho _{\gamma_B}\big(X_T^B+F^{*}- \eta _{B}^{*}\big) =%
\frac{1}{\gamma _{B}}\inf_{\xi _{B}\in \mathcal{V}_{T}^{( B)
}}\rho _{\gamma_B}\big(X_T^B+F^{*}-\xi _{B}\big) =\frac{1}{\gamma
_{C}}\rho _{\gamma_C}\big(X_T^A+X_T^B-\xi ^{*}\big).
\end{equation*}
$ii)$ \noindent $\eta _{A}^{*}$ is an optimal hedging portfolio of $\big(
X_T^A-F^{*}\big) $ for Agent $A$%
\begin{equation*}
\frac{1}{\gamma _{A}}\rho _{\gamma_A}\big( X_T^A-\big( F^{*}+ \eta
_{A}^{*}\big) \big) =\frac{1}{\gamma _{A}}\inf_{\xi _{A}\in \mathcal{V}%
_{T}^{( A) }}\rho _{\gamma_A}\big( X_T^A-\big( F^{*}+\xi _{A}\big)
\big) =\frac{1}{\gamma _{C}}\rho _{\gamma_C}\big( X_T^A+X_T^B-\xi
^{*}\big).
\end{equation*}
\end{theorem}
{\bf Proof:} To prove this theorem, we proceed in several steps:\\
\underline{Step $1$:}\\
Let us first observe that $$ R_{AB}^{m}\left(X_T^A+X_T^B\right)
=\rho_{\gamma _{C}}\left(X_T^A+X_T^B-\xi ^{*}\right)
=\inf_{\widetilde{F}\in \mathcal{X}}\left(\rho_{\gamma _{A}}\left(X_T^A+X_T^B-\widetilde{F}%
-\xi ^{*}\right) +\rho_{\gamma _{B}}\left( \widetilde{F}\right)
\right),$$ where $\widetilde{F}=F+X_T^B-\xi_B$. Given Proposition
\ref{Theorem dilated risk measure}, we
obtain directly an expression for the optimal ''structure'' $\widetilde{F}%
^{*}$ as: $\widetilde{F}^{*}=\frac{\gamma _{B}}{\gamma _{A}+\gamma
_{B}}\left(X_T^A+X_T^B-\xi ^{*}\right) =\frac{\gamma _{B}}{\gamma
_{C}}\left(X_T^A+X_T^B-\xi ^{*}\right)$. Moreover,
$\rho_{\gamma_B}(\widetilde{F})=\frac{\gamma_B}{\gamma_C}(X_T^A+X_T^B-\xi^*)$.
\vspace{1mm}\\
\underline{Step $2$:}\\
Rewriting in the reverse order, we naturally set $F^{*}=
\widetilde{F}^{*}-X_T^B+\eta _{B}^{*}$. We then want to prove that
$\eta
_{B}^{*}$ is an optimal investment for agent $B$.\\
For the sake of simplicity in our notation, we consider
$G^{X}\left( \xi _{A},\xi _{B},F\right) \triangleq \rho_{\gamma
_{A}}\left(X_T^A-F-\xi _{A}\right) +\rho_{\gamma
_{B}}\left(X_T^B+F-\xi
_{B}\right)$.\\
Given the optimality of $\xi ^{*}=\eta _{A}^{*}+\eta _{B}^{*}$ and $%
\widetilde{F}^{*}=F^{*}+X_T^B-\eta _{B}^{*}$, we have
\begin{eqnarray*}
R_{AB}^{m}\left(X_T^A+X_T^B\right) &=&G^{X}\left( \eta
_{A}^{*},\eta
_{B}^{*},F^{*}\right) \\
&=&\inf_{F\in \mathcal{X},\xi _{A}\in \mathcal{V}_{T}^{\left(
A\right) },\xi _{B}\in \mathcal{V}_{T}^{\left( B\right)
}}G^{X}\left( \xi _{A},\xi _{B},F\right) \leq \inf_{\xi _{B}\in
\mathcal{V}_{T}^{\left( B\right) }}G^{X}\left( \eta _{A}^{*},\xi
_{B},F^{*}\right) \leq G^{X}\left( \eta _{A}^{*},\eta
_{B}^{*},F^{*}\right).
\end{eqnarray*}
Then $\eta _{B}^{*}$ is optimal for the problem $\rho_{\gamma
_{B}}\left( F-\xi _{B}\right) \rightarrow \inf_{\xi _{B}\in
\mathcal{V}_{T}^{\left( B\right) }} $. The optimality of $\eta
_{A}^{\ast }$ can be proved using the same arguments. $\square$
\begin{remark}
(a) We first assume that both agents have the same access to the
financial market from a \textit{cone} $\mathcal{H}$. Given the
fact that the risk measure generated by $\mathcal{H}$ is coherent
and thus invariant by dilatation, the market modified risk
measures of both agents are generated from the root risk measure
$\rho \square \nu ^{\mathcal{H}}=\rho^{\mathcal{H}}$ as
$\rho^{\mathcal{H}}_A=\rho _{\gamma_A}\square \nu ^{\mathcal{H}}=
\rho _{\gamma_A}\square \nu _{\gamma _{A}}^{%
\mathcal{H}}=\big( \rho \square \nu ^{\mathcal{H}}\big) _{\gamma _{A}}=
\rho^{\mathcal{H}}_{\gamma_A}%
\text{ \quad and~\quad
}\rho^{\mathcal{H}}_B=\rho^{\mathcal{H}}_{\gamma_B}$.\\
(b) In a more general framework, when both agents have different
access to the financial market, the convex set
$\mathcal{V}_{T}^{(A+B)}$ associated with the risk measure
$\nu^{(A+B)}=\nu^A\square\nu^B$ plays the same role as the set
${\cal H}$ above, since $\rho _{\gamma_C}\square \nu
^{A}\square\nu ^{B}\big(X_T^A+X_T^B\big) =\rho _{\gamma_C}\square
\nu^{(A+B)}\big( X_T^A+X_T^B\big)$.
\end{remark}
{\bf Comments}: Note that when both agents have the same access to the financial market, it
is optimal to transfer the same proportion of the initial
risk as in the problem without market. This result is very strong
as it does not require any specific assumption either for the
non-tradable risk or the financial market. Moreover, the optimal
structure $F^*$ does not depend on the financial market. The
impact of the financial market is simply visible through the
pricing rule, which depends on the market modified risk measure of
agent B.\\
Standard diversification will also occur in exchange economies as
soon as agents have proportional penalty functions. The regulator
has to impose very different rules on agents as to generate risk
measures with non-proportional penalty functions if she wants to
increase the diversification in the market. In other words,
diversification occurs when agents are very different one from the
other. This result supports for instance the intervention of
reinsurance companies on financial markets in order to increase
the diversification on the reinsurance market.
\subsubsection{Optimal Transaction in the General Framework\label{Section optimal design general}}
We now come back to our initial problem of optimal risk transfer between agent A and agent B, when now they both
have access to the financial market to hedge and diversify their respective portfolio.
\paragraph{General framework}
As in the dilated framework, we assume that both their
risk measures $\rho_A$ and $\rho_B$ are defined on
$\mathbb{L}_{\infty}(\mathbb{P})$ and are continuous from above.
The investment opportunities of both agents are described by two
convex subsets $\mathcal{V}^{A}_{T}$ and $\mathcal{V}^{B}_{T}$ of
$\mathbb{L}_{\infty}(\mathbb{P})$ and the financial market is
assumed to be arbitrage-free.
\begin{enumerate}
\item This opportunity to invest optimally in a financial market
reduces the risk of both agents. To assess their respective risk
exposure, they now refer to market modified risk measures $\rho
_{A}^{m}$ and $\rho _{B}^{m}$ defined if $J=A,B$ as $\rho
_{J}^{m}( \Psi)\triangleq \inf_{\xi_{J}\in \mathcal{V}_{T}^{( J)
}}\rho _{J}( \Psi -\xi _{J}) $. As usual, we assume that $\rho
_{J}^{m}( 0) >-\infty $ for the individual hedging programs to
make sense. Thanks to Corollary \ref{corollary modification de rho
par H},
\begin{equation*}
\rho _{A}^{m}( \Psi) =\rho _{A}\square \nu ^{A}( \Psi ) \qquad
\text{and\qquad }\rho _{B}^{m}(\Psi) =\rho _{B}\square \nu ^{B}(
\Psi ).
\end{equation*}
\item Consequently, the optimization program related to the $F$%
-transaction is simply
\begin{equation*}
\inf_{F,\pi }\rho _{A}^{m}\big( X^A_T-F+\pi \big) \qquad \text{%
subject to \qquad }\rho _{B}^{m}\big(X_T^B + F-\pi \big) \leq \rho
_{B}^{m}\big(X_T^B\big).
\end{equation*}
As previously, using the cash translation invariance property and binding
the constraint at the optimum, the pricing rule of the $F$-structure is
fully determined by the buyer as
\begin{equation}
\pi ^{*}( F) =\rho _{B}^{m}\big(X_T^B\big) -\rho _{B}^{m}\big(
X_T^B+F\big).  \label{Equation prix avec marche}
\end{equation}
It corresponds to an ''indifference'' pricing rule from the agent $B$'s market modified risk measure.
\item Using again the cash translation invariance property, the optimization
program simply becomes
\begin{equation*}
\inf_{F}\big\{ \rho _{A}^{m}\big( X_T^A-F\big) +\rho
_{B}^{m}\big(X_T^B+F\big)\big\} -\rho _{B}^{m}\big(X_T^B\big)
\triangleq R_{AB}^{m}\big(X_T^A+X_T^B\big)-\rho
_{B}^{m}\big(X_T^B\big).
\end{equation*}
With the functional $R_{AB}^{m}$, we are in the framework of
Theorem
 \ref{Theorem modification de rho1
par A_rho2}.
\begin{eqnarray}
R_{AB}^{m}\big(X_T^A+X_T^B\big) &=&\inf_{F}\big\{\rho
_{A}^{m}\big( X_T^A-F\big) +\rho _{B}^{m}\big(X_T^B+F\big) \big\}
\label{pg general avec Fbar} \\
&=& \inf_{\widetilde{F}}\big\{\rho _{A}^{m}\big(
X_T^A+X_T^B-\widetilde{F}\big) +\rho _{B}^{m}\big(\widetilde{F}\big) \big\} = \rho_{A}^{m}\square \rho _{B}^{m} \big(X_T^A+X_T^B \big) \notag\\
&=& \rho_{A}\square \nu ^{A}\square \rho _{B}\square \nu
^{B}\big(X_T^A+X_T^B \big). \label{equation RABm general}
\end{eqnarray}
The value functional $R_{AB}^{m}$ of this program, resulting from
the inf-convolution of four different risk measures, may be
interpreted as the \emph{residual risk measure} after all
transactions. This inf-convolution problem makes sense if the
initial condition $\rho_A^m \square \rho^m_B(0)>-\infty$ is
satisfied. \item Using the previous Theorem \ref{Theorem
modification de rho1 par A_rho2} on the stability of convex risk
measure, provided the initial condition is satisfied, $R_{AB}^{m}$
is a convex risk measure with the penalty function
$\alpha_{AB}^{m}=\alpha _{A}^{m}+\alpha _{B}^{m} =\alpha _{A}
 +\alpha _{B} +\alpha^{\mathcal{V}_T^{A}}+\alpha^{\mathcal{V}_T^{B}}$.
\end{enumerate}
{\bf Comments:}  The general risk transfer problem
can be viewed as a game involving four different agents if the access to
the financial market is different for agent $A$ and agent $B$ (or three otherwise).
As a consequence,
we end up with an inf-convolution problem involving four different risk measures, two per agents.
\paragraph{Optimal design problem}
Our problem is to find an optimal structure $F^{\ast }$ realizing
the minimum of the Program (\ref{pg general avec Fbar}):
\begin{equation*}
R_{AB}^{m}\big(X_T^A+X_T^B\big) =\inf_{F}\big\{\rho _{A}^{m}\big(
X_T^A-F\big) +\rho_{B}^{m}\big(X_T^B+F\big) \big\}
\end{equation*}
Let us first consider the following simple inf-convolution problem
between a convex risk measure $\rho_B$ and a linear function $q_A$
as introduced in Subsection \ref{subsection properties
inf-convolution}:
\begin{equation}
\label{eq:linear convolution} q_A \square \rho_B(X)=\inf_F
\{\mathbb{E}_{\mathbb{Q}_A}[-(X-F)]+\rho_B(F)\}.
\end{equation}
\begin{proposition}
The necessary and sufficient condition to have an optimal solution
$F^*$ to the linear inf-convolution problem (\ref{eq:linear
convolution}) is expressed in terms of the subdifferential of
$\rho_B$ as $\mathbb{Q}_A \in
\partial \rho_B(F^*)$.
\end{proposition}
This necessary and sufficient corresponds to the first order
condition of the optimization problem. More generally, the
following result is obtained:
\begin{theorem}[Characterization of the optimal]
\label{Theorem CNS existence} Assume that $\rho
_{A}^{m}\square \rho_{B}^{m}(0) > -\infty$.\\
The inf-convolution program
\begin{equation*}
R^m_{AB}\big(X_T^A+X_T^B\big) =\inf_{F}\big\{ \rho^m_{A}\big(
X_T^A-F\big) +\rho^m_{B}\big(X_T^B+F\big) \big\}
\end{equation*}
is exact at $F^*$ if and only if there exists $\mathbf{Q}_{AB}^{X}
\in \partial R^m_{AB}(X_T^A+X_T^B)$ such that $\mathbf{Q}_{AB}^{X}
\in \partial \rho_A^m(X_T^A-F^*) \cap
\partial \rho_B^m(X_T^B+F^*)$.\\
In other words, the necessary and sufficient condition to have an
optimal solution $F^{*}$ to the inf-convolution program is that
there exists an optimal additive measure $\mathbf{Q}_{AB}^{X}$ for
$\big(X_T^A+X_T^B,R^m_{AB}\big)$ such that $X_T^B+F^{*}$ is
optimal for $\big(\mathbf{Q}_{AB}^{X},\alpha^m_{B}\big) $ and
$X_T^A-F^{*}$ is optimal for $\big( \mathbf{Q}_{AB}^{X},\alpha^m_{A}\big) $.\\
\end{theorem}
Both notions of optimality are rather intuitive as they simply
translate the fact that the dual representations of the risk
measure on the one hand, and of the penalty function on the other
hand, are exact respectively at a given additive measure and at a
given
exposure.\\
A natural interpretation of this theorem is that both
agents agree on the measure $\mathbf{Q}_{AB}^{X}$ in order to
value their respective residual risk. This agreement enables the
transaction.\vspace{1mm}\\
\textbf{Proof:}\\
Let us denote by $\mathbf{Q}_{AB}^{X}$ the optimal additive
measure for $\big(X_T^A+X_T^B,R^m_{AB}\big) $. In this case,
$\mathbf{Q}_{AB}^X \in \partial R^m_{AB}(X_T^A+X_T^B)$. As
mentioned in Subsection \ref{risk measure basic properties}, the
existence of such an additive measure is guaranteed as soon as the
penalty function is defined by Equation $(\ref{penalty 1
(acceptance
set)})$. This justifies the writing of the theorem in terms of additive measures rather
than in terms of probability measures.\\
$i)$ In the proof, we denote by $X \triangleq
X_T^A+X_T^B$ and by $\Psi ^{c}$, the centered random
variable $\Psi $ with respect to the given additive measure $\mathbf{Q}%
_{AB}^{X}$ optimal for $\big( X,R_{AB}\big) $: $\Psi ^{c}=\Psi
-\mathbb{E}_{\mathbf{Q}_{AB}^{X}} [\Psi].$ So, by definition,
\begin{eqnarray*}
-R_{AB}\big( X^{c}\big) &=&\alpha _{A}\big( \mathbf{Q}%
_{AB}^{X}\big) +\alpha _{B}\big( \mathbf{Q}_{AB}^{X}\big) \\
&=&\sup_{F}\big\{ -\rho _{A}\big( X^{c}-F^{c}\big) \big\}
+\sup_{F}\big\{ -\rho _{B}\big( F^{c}\big)
\big\} \\
&\geq &-\inf_{F}\big\{ \rho _{A}\big( X^{c}-F^{c}\big) +\rho
_{B}\big( F^{c}\big) \big\} =-R_{AB}\big( X^{c}\big).
\end{eqnarray*}
In particular, all inequalities are equalities and
\begin{equation*}
\sup_{F}\big\{ -\rho _{A}\big( X^{c}-F^{c}\big) \big\}
+\sup_{F}\big\{ -\rho _{B}\big( F^{c}\big) \big\} =\sup_{F}\big\{
-\rho _{A}\big( X^{c}-F^{c}\big) -\rho _{B}\big( F^{c}\big)
\big\}.
\end{equation*}
Hence, $F^{*}$ is optimal for the inf-convolution problem, or equivalently
for the program on the right-hand side of this equality, if and only if $%
F^{*}$ is optimal for both problems $\sup_{F}\big\{ -\rho
_{B}\big( F^{c}\big) \big\} $ and $\sup_{F}\big\{ -\rho _{A}\big(
X^{c}-F^{c}\big) \big\} $.\\
The second formulation is a straightforward application of Theorem
\ref{Sandwich theorem} ii), considering the problem not at $0$ but
at $X_T^A+X_T^B$. $\square$\vspace{2mm}\\
In order to obtain an explicit representation of an optimal
structure $F^*$, some technical methods involving a localization
of convex risk measures have to be used. This is the aim of the
second part of this chapter, which is based upon some technical
results on BSDEs. Therefore, before localizing convex risk
measures and studying our optimal risk transfer in this new
framework, we present in a separate section some quick recalls on
BSDEs, which is essential for a good understanding of the second
part on dynamic risk measures.
\begin{center}
\section*{Part II: Dynamic Risk Measures
\label{section local risk measures} }
\end{center}
We now consider {\it dynamic convex risk
measures}\index[sub]{dynamic convex risk measure}. Quite recently,
many authors have studied dynamic version of static risk measures,
focusing especially on the question of law invariance of these
dynamic risk measures: among many other references, one may quote
the papers by \index[aut]{Cvitanic}Cvitanic and
\index[aut]{Karatzas}Karatzas \cite{CvitanicKaratzas99},
\index[aut]{Wang}Wang \cite{Wang},\index[aut]{Scandolo}Scandolo
\cite{Scandolo}, \index[aut]{Weber}Weber \cite{Weber},
\index[aut]{Artzner} \index[aut]{Delbaen} \index[aut]{Eber}
\index[aut]{Heath}Artzner et al. \cite{ArtznerDelbaenEberHeath2},
\index[aut]{Cheridito}Cheridito, \index[aut]{Delbaen}Delbaen and
\index[aut]{Kupper}Kupper \cite{CheriditoDelbaenKupper1}
\cite{CheriditoDelbaenKupper2} or \cite{CheriditoDelbaenKupper3},
\index[aut]{Detlefsen}Detlefsen and \index[aut]{Scandolo}Scandolo
\cite{DetlefsenScandolo}, \index[aut]{Frittelli}Frittelli and
\index[aut]{Gianin}Gianin \cite{FrittelliGianin04},
\index[aut]{Frittelli}Frittelli and \index[aut]{Scandolo}Scandolo
\cite{FrittelliScandolo}, \index[aut]{Gianin}Gianin \cite{Gianin},
\index[aut]{Riedel}Riedel \cite{Riedel},
\index[aut]{Roorda}Roorda, \index[aut]{Schumacher}Schumacher and
\index[aut]{Engwerda}Engwerda \cite{RoordaSchumacherEngwerda} or
the lecture notes of \index[aut]{peng}Peng \cite{Peng03}. Very
recently, extending the work of \index[aut]{El Karoui}El Karoui
and \index[aut]{Quenez}Quenez \cite{ElKarouiQuenez95},
\index[aut]{Kl\"{o}ppel}Kl\"{o}ppel and
\index[aut]{Schweizer}Schweizer have related dynamic indifference
pricing and BSDEs in
\cite{KloppelSchweizer}.\vspace{2mm}\\
In this second part, we extend the axiomatic approach adopted in
the static framework and introduce some additional axioms for the
risk measures to be time-consistent. We then relate the dynamic
version of convex risk measures to BSDEs. The associated dynamic
risk measure is called $g$-conditional risk measure, where $g$ is
the BSDE coefficient. We will see how the properties of both the
risk measure and the coefficient $g$ are intimately connected. In
particular, one of the key axioms in the characterization of the
dynamic convex risk measure will be the translation invariance, as
we will see in Section \ref{section axiomatic approach dynamic},
and this will impose the
$g$-coefficient of the related BSDE to depend only on $z$.\\
In the last two sections, we come back to the essential point of
this chapter, the optimal risk transfer problem. We first derive
some results on the inf-convolution of dynamic convex risk
measures and obtain the optimal structure as a solution to the
inf-convolution problem.\vspace{1mm}\\
The idea behind our approach is to find a trade-off between static
and very abstract risk measures as to obtain tractable risk
measures. Therefore, we are more interested in tractability issues
and interpretations of the dynamic risk measures we obtain rather
than the ultimate general results in BSDEs.
\section{Some recalls on Backward Stochastic Differential Equations}
In the rest of the chapter, we take into account more information
on the risk structure. In particular, we assume the $\sigma$-field
$\cal F$ generated by a d-dimensional Brownian motion between
$[0,T]$. Since any bounded ${\cal F}_T$-measurable variable is an
stochastic integral w.r. to the Brownian motion, the risk measures
of interest have to be robust with respect of this localization
principle. To do that, we consider a family of risk measures
described by backward
stochastic differential equations (BSDE).\\
In this section, we introduce general BSDEs, defining them,
recalling some key results on existence and uniqueness of a
solution and presenting the comparison theorem. Complete proofs
and additional useful results are given in the Chapter dedicated
to BSDEs.
\subsection{General Framework and Definition}
Let $\big(\Omega ,{\cal F},\mathbb{P}\big)$
\index[not]{$\mathbb{P}$} be a probability space on which is
defined a $d$-dimensional Brownian motion \index[not]{$W$}
$W:=(W_t;t\leq T_H)$, where \index[not]{$T_H$} $T_H>0$ is the time
horizon of the study. Let us consider the natural Brownian
filtration ${\cal F} _{t}^0=\sigma \big(W_{s};0\leq s\leq t ;t\geq
0\big)$ and \index[not]{$\mathcal{F}_t$} $({\cal F}_t;t\leq T_H)$
its completion with the
$\mathbb{P}$-null sets of ${\cal F}$.\\
Denoting by \index[not]{$\mathbb E$} $\mathbb E$ the expected
value with respect to $\mathbb P$, we introduce the following
spaces which will be important in the formal setting of BSDEs.
Since the time horizon may be sometimes modified, the definitions
are referring
to a generic time $T\leq T_H$.\\
$\bullet$ \index[not]{$L_n^2 \big({\cal F}_t\big)$} $L_n^2
\big({\cal F}_t\big)= \{\eta\, : \, {\cal F}_t-\mbox{measurable}
\, {\mathbb R}^n- \mbox{valued random variable s.t.} \; {\mathbb E}(|\eta|^2)<\infty \}$.\\
$\bullet$ \index[not]{${\cal P}_n(0, T)$} ${\cal P}_n(0,
T)=\{(\phi_t;0 \leq t\leq T)\,  : \, \mbox{progressively
measurable process with values in }\;\mathbb{R}^n \}$\\
$\bullet$ \index[not]{${\cal S}_n^2 (0, T)$} ${\cal S}_n^2 (0, T)
= \{(\phi_t;0 \leq t\leq T)\, : \, \phi\in {\cal P}_n\,
\mbox{s.t.} \;
{\mathbb E}[\sup_{t\leq T}|Y_t|^2]<\infty\}$ .\\
$\bullet$ \index[not]{${\cal H}_n^{2}(0,T)$} ${\cal H}_n^{2}(0,T)
= \{(\phi_t;0 \leq t \leq T) : \, \, \phi\in {\cal P}_n\,
 \,\mbox{s.t.} \;{\mathbb E}[\int_{0}^{T}|Z_s|^2ds]<\infty \}$. \\
$\bullet$ \index[not]{${\cal H}_n^{1}(0,T)$} ${\cal H}_n^{1}(0,T)
= \{(\phi_t;0 \leq t \leq T) : \, \, \phi\in {\cal P}_n\,
 \,\mbox{s.t.} \;{\mathbb E}[(\int_{0}^{T}|Z_s|^2ds)^{1/2}]<\infty \}$.
\vspace{2mm}\\
Let us give the definition of the one-dimensional BSDE; the
multidimensional case is considered in the book's chapter
dedicated to BSDEs.
\begin{definition}
Let $ \xi_T \in L^{2}(\Omega ,{\cal F}_{T},\mathbb{P})$ be a
$\mathbb{R}$-valued terminal condition and $g$ a coefficient
${\cal P}_1\otimes{\cal B}({\mathbb R})\otimes{\cal B}({\mathbb
R}^d)$-measurable. A solution for the BSDE associated with
$(g,\xi_T)$ is a pair of progressively measurable processes
$(Y_t,Z_t)_{t\leq T}$, with values in
$\mathbb{R}\times\mathbb{R}^{1\times d}$ such that:
\begin{equation}
\label{pp} \Big\{
\begin{array}{l}
(Y_t) \in {\cal S}_1^2(0,T), \,\,(Z_t) \in {\cal H}_{1\times d}^{2}(0,T)\\
Y_{t}= \xi_T
+\int_{t}^{T}g(s,Y_{s},Z_{s})ds-\int_{t}^{T}Z_{s}dW_{s},\,\, 0
\leq t\leq T\, .
\end{array}
 \Big.
 \end{equation}
The following differential form is also useful
\begin{equation}
\label{eq:difpp} -dY_{t}= g(t,Y_{t},Z_{t})dt-Z_{t}dW_{t},\quad
Y_T=\xi_T.
\end{equation} \end{definition}
{\bf Conventional notation}: To simply the writing of the BSDE, we
adopt the following notations: the Brownian motion $W$ is
described as a column vector $(d,1)$ and the $Z$ vector is
described as a row vector $(1,d)$ such that the notation $ZdW$ has
to be understood as a matrix product with (1,1)-dimension.
\begin{remark}
If $\xi_T$ and $ g(t,y,z)$ are deterministic, then $ Z_t \,\equiv
\, 0$, and  $(Y_t)$ is the solution of ODE
$$
\frac{dy_t}{dt} = - g (t,y_t,0), \quad \quad y_T = \xi_T\, .
$$
If the final condition $\xi_T$ is random, the previous solution is
${\cal F}_{T}$-measurable, and so non adapted. So we need to
introduce the martingale $\int_{0}^{t}Z_{s}dW_{s}$ as a control
process
 to obtain an adapted solution.
\end{remark}
\subsection{Some Key Results on BSDEs}
Before presenting key results of BSDEs, we first summarize the
results concerning the existence and uniqueness of a solution. The
proofs are given in the Chapter dedicated to BSDEs with some
complementary results.
\subsubsection{Existence and Uniqueness Results}
In the following, we always assume  the necessary condition on the
terminal condition $\xi_T \in L^2(\Omega ,{\cal
F}_{T},\mathbb{P})$.
\begin{enumerate}
\item {\sc \index[not]{$(H1)$} (H1): The standard case (uniformly Lipschitz)}:
$(g(t,0,0);0 \leq t\leq T)$ belongs to ${\cal H}^{2}(0,T)$ and $g$
uniformly Lipschitz continuous with respect to $(y,z)$, $i.e.$
there exists a constant $C\geq 0$ such that
$$ d{\mathbb P}\times dt-a.s.\quad \forall (y,y^{\prime },z, z^{\prime })\quad|g( \omega ,t,y,z)
-g(\omega ,t,y',z')| \leq C( | y-y'| +|z-z'|).$$ Under these
assumptions, \index[aut]{Pardoux} Pardoux and \index[aut]{Peng}
Peng \cite{PardouxPeng90} proved in 1990 the existence and
uniqueness of a solution.
\item {\sc \index[not]{$(H2)$} (H2) The continuous case with linear growth}: there exists
a constant $C\geq 0$ such that
$$d{\mathbb P}\times dt-a.s.\quad \forall (y,z)\quad
|g( \omega ,t,y,z) | \leq k(1+ | y| + |z|).$$ Moreover we assume
that $d{\mathbb P}\times dt \>a.s.$, $g(\omega,t,.,.)$ is
continuous in $ (y,z)$. Then, there exist a maximal and a minimal
solutions (for a precise definition, please refer to the Chapter
dedicated to BSDEs), as proved by \index[aut]{Lepeltier} Lepeltier
and \index[aut]{San Martin} San Martin in 1998
\cite{LepeltierSanMartin98}.
\item {\sc \index[not]{$(H3)$} (H3) The continuous case with quadratic growth in $z$}:
In this case, the assumption of square integrability on the
solution is too strong. So we only consider bounded solution and
obviously terminal condition $\xi_T\in L_{\infty}$. We also
suppose that there exists a constant $k\geq 0$ such that
$$d{\mathbb P}\times dt-a.s.\quad \forall (y,z)\quad|g( \omega ,t,y,z) |
\leq k( 1+ | y| + |z|^2).$$ Moreover we assume that $d{\mathbb
P}\times dt -a.s.$, $g(\omega,t,.,.)$ is continuous in
$ (y,z)$.\\
Then there exist a maximal and a minimal bounded solutions as
first proved  by \index[aut]{Kobylansky} Kobylansky
\cite{Kobylanski00} in 2000 and extended by \index[aut]{Lepeltier}
Lepeltier and \index[aut]{San Martin} San Martin
\cite{LepeltierSanMartin98} in 1998. The uniqueness of the
solution was proved by \index[aut]{Kobylansky} Kobylansky
\cite{Kobylanski00} under the additional conditions that the
coefficient $g$ is differentiable in $(y,z)$ on a compact interval
$[-K, \, K] \times \mathbb{R}^d$ and that there exists $ c_1
> 0 \;{\rm and} \; c_2
> 0 $ such that:
\begin{equation}
\label{uniqueness-quadratic} \frac{\partial  g}{\partial z} \leq
c_1 (1 + |z|), \quad \quad \frac{\partial  g}{\partial y} \leq c_2
(1 + |z|^2)
\end{equation}
\end{enumerate}
\subsubsection{Comparison Theorem}
We first present an important tool in the study of one-dimensional
BSDEs: the so-called \emph{comparison theorem}. It is the
equivalent of the maximum principle when working with PDEs.
\begin{theorem}[Comparison Theorem]\label{comparison Theorem}
Let $(\xi^{1}_T,g^{1})$ and $(\xi^{2}_T, g^{2})$ be two pairs
(terminal condition, coefficient) satisfying one of the above
conditions (H1,H2,H3) (but the same for both pairs). Let
$(Y^1,Z^{1})$ and $(Y^2,Z^{2})$ be the maximal associated
solutions.\\
$(i)$ We assume that $\quad \xi^{1}_T\leq \xi^{2}_T$, $\mathbb{P}
-a.s.$ and that $\quad d{\mathbb P}\times dt-a.s.\quad
\forall(y,z)\quad  g^{1}(\omega, t,y,z) \leq g^{2}(\omega, t,y, z)
$. Then we have
$$
Y^1_t \leq  Y^2_t \; \; \mbox{a.s.} \; \forall \, t\in[0,T]
$$
$(ii)$ {\bf Strict inequality} Moreover, under (H1), if in
addition $Y_t^1 =Y_t^2$ on $B\in {\cal F }_t$, then a strict
version of this result holds as
$$a.s. \; {\rm on} \; B \quad \xi^1_T = \xi^2_T \;,\; \forall s\geq t, \; Y_s^1 = Y_s^2
\quad {\rm and} \quad g^1(s, Y_s^1,Z_s^1) = g^2(s,Y_s^2, Z_s^2)
\quad d{\mathbb P}\times ds-a.s. \; {\rm on}\; B\times [t,T]$$
\end{theorem}
\section{Axiomatic Approach and $g$-Conditional Risk Measures
\label{section axiomatic approach dynamic}}
In this section, we give a general axiomatic approach for dynamic
convex risk measures and see how they are connected to the
existing notions of consistent convex price systems and non-linear
expectations, respectively introduced by \index[aut]{El Karoui} El
Karoui and \index[aut]{Quenez} Quenez \cite{ElKarouiQuenez96} and
\index[aut]{Peng} Peng \cite{Peng97}. Then, we relate the dynamic
risk measures with BSDEs and focus on the properties of the
solution of some particular BSDEs associated with a convex
coefficient $g$, called g-conditional risk measures.
\subsection{Axiomatic Approach}
Following the study of static risk measures by
\index[aut]{F\"{o}llmer} F\"{o}llmer and \index[aut]{Schied}
Schied \cite{FoellmerSchied02a} and \cite{FoellmerSchied02b}, we
now propose a common axiomatic approach to dynamic convex risk
measures, non-linear expectations and convex price systems and
non-linear.
\begin{definition} \label{def:axiomatic} Let
$\big(\Omega,\mathcal{F},\mathbb{P},({\cal F}_{t};t \geq 0)\big)$
be a filtered probability space. A dynamic $L^2$-operator
($L^{\infty}$-operator)\index[sub]{dynamic operator}
\index[not]{$\mathcal{Y}$} ${\cal Y}$ with respect to $\big({\cal
F}_{t};t \geq 0\big)$ is a family of continuous semi-martingales
which maps, for any bounded stopping time $T$, a $L^2({\cal F}_T)$
(resp. $L^{\infty}({\cal F}_T))$ -variable $\xi_T$ onto a process
$\big({\cal Y}_t(\xi_T);t\in[0,T]\big)$. Such an operator is said
to be
\begin{enumerate}
\item {\sc \index[not]{$(P1)$} (P1) Convex}: For any stopping times $S\leq T$,  for any
$(\xi^1_T,\xi^2_T)$, for any $0\leq \lambda\leq 1$,
$${\cal Y}_S(\lambda\xi^1_T+(1-\lambda)\xi^2_T)\leq \lambda{\cal Y}_S(\xi^1_T)+
(1-\lambda){\cal Y}_S(\xi^2_T)\quad {\mathbb P}-a.s.$$ \item {\sc
\index[not]{$(P2)$} (P2) Monotonic}: For any stopping times $S\leq
T$, for any $(\xi^1_T,\xi^2_T)$ such that $\xi^1_T\geq \xi^2_T\>
a.s.$,
\begin{description}
\item {\sc \index[not]{$(P2+)$} (P2+)}: \hspace{3mm}  the operator is increasing if
${\cal Y}_S(\xi^1_T) \geq {\cal Y}_S(\xi^2_T)\quad a.s.$ \item
{\sc \index[not]{$(P2-)$} (P2-)}: \hspace{3mm} the operator is
decreasing   ${\cal Y}_S(\xi^1_T) \leq {\cal Y}_S(\xi^2_T)\quad
a.s.$
\end{description}
\item {\sc \index[not]{$(P3)$} (P3) Translation invariant}: For any stopping times
$S\leq T$ and any $\eta_S\in {\cal F}_S$, for any
$\xi_T$,\vspace{2mm}\\
{\sc \index[not]{$(P3+)$} (P3+)} \hspace{3mm}  ${\cal
Y}_S(\xi_T+\eta_S)={\cal Y}_S(\xi_T)-\eta_S\quad
a.s.$,\hspace{6mm} {\sc \index[not]{$(P3-)$} (P3-)}\hspace{3mm}
${\cal Y}_S(\xi_T+\eta_S)={\cal Y}_S(\xi_T)-\eta_S\quad a.s.$
\item {\sc \index[not]{$(P4)$} (P4) Time-consistent}: For $S\leq T\leq U$ three
bounded stopping times, for any $\xi_U$ \vspace{2mm}\\
{\sc \index[not]{$(P4+)$} (P4+)} \hspace{3mm} ${\cal Y}_S\big(
\xi_U\big)={\cal Y}_S\big({\cal Y}_{T}\big( \xi_U\big)\big)\quad
a.s.$,\hspace{6mm} {\sc \index[not]{$(P4-)$} (P4-)} \hspace{3mm}
${\cal Y}_S\big( \xi_U\big)= {\cal Y}_S\big(-{\cal
Y}_{T}\big(\xi_U\big)\big)\quad a.s.$ \item {\sc
\index[not]{$(P5)$} (P5) Arbitrage-free}: For any stopping times
$S\leq T$, and for any $(\xi^1_T,\xi^2_T)$ such that
$\xi_T^1\geq\xi_T^2$,
$${\cal Y}_S(\xi_T^1)={\cal Y}_S(\xi_T^2) \quad {\rm on}\quad A_S=\{S<T\}
\Longrightarrow\> \xi_T^1=\xi_T^2\quad a.s.\>{\rm on}\>A_S.$$
\item {\sc \index[not]{$(P6)$} (P6) Conditionally invariant}: For any stopping times
$S\leq T$ and any $B\in {\cal F}_S$, for any $\xi_T$,
$${\cal Y}_S({\mathbf 1}_B\,\xi_T)={\mathbf 1}_B{\cal Y}_S(\xi_T) \quad a.s.$$
\item{\sc \index[not]{$(P7)$} (P7) Positive homogeneous}: For any stopping times
$S\leq T$, for any $\lambda_S \geq 0$ ($\lambda_S \in {\cal F}_S$)
and for any $\xi_T$,
$${\cal Y}_S(\lambda_S \xi_T)=\lambda_S {\cal Y}_S(\xi_T) \quad a.s.$$
\end{enumerate}
\end{definition}
First, note that the property (P5) of no-arbitrage implies that
the monotonicity property (P2)
is strict.\vspace{1mm}\\
Most axioms have two different versions, depending on the sign
involved. Making such a distinction is completely coherent with
the previous observations in the static part of this chapter about
the relationship between price and risk measure: since the
opposite of a risk measure is a price, the axioms with a "$+$"
sign are related to the characterization of a price system, while
the axioms with a "$-$" sign are
related to that of a dynamic risk measure.\\
In \cite{ElKarouiQuenez96}, when studying pricing problems under
constraints, \index[aut]{El Karoui} El Karoui and
\index[aut]{Quenez} Quenez defined
 a {\it consistent convex (forward) price system} \index[sub]{consistent convex price system} as a convex (P1), increasing
(P2+), time-consistency (P4+) dynamic operator ${\cal P}_t$,
without arbitrage (P5). Time-consistency (P4+) may be view as a
dynamic
programming principle.\\
At the same period, Peng introduced the notion of {\it non-linear
expectation} \index[sub]{non-linear expectation} as a {\bf
translation invariance} (P3+) convex price system, satisfying the
conditional invariance property (P6) which is very intuitive in
this framework, (see for instance, \index[aut]{Peng} Peng
\cite{Peng97}). Note that (P6) of conditional invariance implies
some additional assumptions on the operator ${\cal Y}$: in
particular for any $t$, ${\cal Y}_t(0)=0$. In the
following, we denote the non-linear expectation by ${\cal E}$.\\
Now on, we focus on dynamic convex risk measures, where now only
the properties with the "$-$" sign hold.
\begin{definition}
A dynamic operator satisfying the axioms of convexity (P1),
decreasing monotonicity (P2-), translation invariance (P3-),
time-consistency (P4-) and arbitrage-free (P5) is said to be a
{\it dynamic convex risk measure}. It will be denoted by
\index[not]{$\mathcal{R}$} ${\cal R}$ in
the following.\\
If ${\cal R}$ also satisfies the positive homogeneity property
(P7), then it is called a {\it dynamic coherent risk measure}.
\end{definition}
Note that a non-linear expectation defines a dynamic risk measure
conditionally invariant and centered.
\begin{remark}
It is not obvious to find a negligible set ${\cal N}$ such that
for any bounded stopping time $S$ and any bounded $\xi_T$,
$\forall \> \omega \notin {\cal N}$, $\xi_T \rightarrow {\cal
R}^g_S(\omega,\xi_T)$ is a static convex risk measure. The
negligible sets may depend on the variable $\xi_T$ itself.
\end{remark}
\paragraph{Dynamic Entropic Risk Measure}
A typical example is the {\it dynamic entropic risk measure}
\index[sub]{dynamic entropic risk measure}, obtained by
conditioning the static entropic risk measure. For any $\xi_T$
bounded: \index[not]{$e_{\gamma,t}$}
$$
e_{\gamma}(\xi_T)= \gamma \ln \mathbb{E}\big[\exp(-\frac{1}{\gamma
}\xi_T)\big] \qquad \Rightarrow \qquad e_{\gamma ,t}(\xi_T)
=\gamma \ln \mathbb{E}\big[ \exp (-\frac{1}{\gamma }\xi_T) |{\cal
F} _{t}\big].
$$
Since, $\xi_T$ is bounded, $e_{\gamma ,t}(\xi_T)$ is bounded for
any $t$. Therefore, this dynamic operator defined on $L^\infty$
satisfies the properties of dynamic convex risk measures.
Convexity, decreasing monotonicity, translation invariance,
no-arbitrage are obvious; the time-consistency property $(P4-)$
results from the transitivity of conditional expectation:
$$\forall t \geq 0 \>,\> \forall h>0,\quad e_{\gamma ,t}\big(
\xi_T\big)= e_{\gamma,t}\big(-e_{\gamma ,t+h}\big( \xi_T\big)\big)
\quad a.s.$$ We give the easy proof of this identity to help the
reader to understand the $(-)$ sign in
the formula.\\
{\bf Proof:}
$$\begin{array}{llll}
Y_{t+h} &\equiv& e_{\gamma ,t+h}\big( \xi_T\big)=\gamma \ln
\mathbb{E}\big[\exp (-\frac{1}{\gamma }\xi_T) |{\cal
F} _{t+h}\big]\\
e_{\gamma,t}\big(-e_{\gamma ,t+h}\big( \xi_T\big)\big)&=&\gamma
\ln \mathbb{E}\big[\exp (-\frac{1}{\gamma }(-Y_{t+h}))
|{\cal F} _{t}\big]\\
 &=& \gamma \ln
\mathbb{E}\Big[\exp (\frac{1}{\gamma }\gamma \ln
\mathbb{E}\big[\exp (-\frac{1}{\gamma }\xi_T) |{\cal F}
_{t+h}\big])
|{\cal F} _{t}\Big]\\
 &=&\gamma \ln
\mathbb{E}\Big[\mathbb{E}\big[\exp (-\frac{1}{\gamma }\xi_T)
|{\cal F} _{t}\big] \Big].\hspace{36mm}\Box
\end{array}
$$
Moreover, it is possible to relate the dynamic entropic risk
measure $e_{\gamma,t}$ with the solution of a BSDE, as follows:
\begin{proposition}
The dynamic entropic measure
$\big(e_{\gamma,t}(\xi_T);t\in[0,T]\big)$ is solution of the
following BSDE with the quadratic coefficient $g\big(
t,z\big)=\frac{1}{2\gamma }\big\| z\big\| ^{2}$ and terminal
bounded condition $\xi_T$.
 \begin{equation}
\label{quadratic} -de_{\gamma ,t}\big(\xi_T\big) =\frac{1}{2\gamma
}\big\| Z_{t}\big\| ^{2}dt- Z_{t}dW_{t}\qquad e_{\gamma
,T}\big(\xi_T\big) =-\xi_T.
\end{equation}
\end{proposition}
{\bf Proof:} Let us denote by
$M_{t}( \xi_T) =\mathbb{E}\Big[\exp \big( -\frac{1}{%
\gamma } \xi_T\big)|{\cal F} _{t}\Big]$. As $M$ is a positive and
bounded continuous martingale, one can use the multiplicative
decomposition to get $dM_t = \frac{1}{\gamma
}M_t\big(Z_t\,dW_{t}\big)$ where $\big( Z_t;t\geq 0\big) $ is  a
$1\times d $ dimensional square-integrable process.  By It\^{o}'s
formula applied to the function $\gamma\ln(x)$, we obtain the
Equation (\ref{quadratic}).\\
Note that the conditional expectation of the quadratic variation
$\mathbb{E}\big[\int_t^T|Z_s|^2ds|{\cal
F}_t]=\mathbb{E}\big[e_{\gamma,t}(\xi_T)-\xi_T|{\cal F}_t\big]$ is
bounded and conversely if the Equation (\ref{quadratic}) has the
solution (Y,Z) such that $Y_T$ and
$\mathbb{E}\big[\int_t^T|Z_s|^2ds|{\cal F}_t]$ are bounded, then
$Y$ is bounded. This point will be detailed in Theorem
\ref{theorem dual representation}. $\Box$
\vspace{2mm}\\
This relationship between the dynamic entropic risk measure and
BSDE can be extended to general dynamic convex risk measures as we
will see in the rest of this section.
\subsection{Dynamic Convex Risk Measures and BSDEs}
This section is about the relationship between dynamic convex risk
measures and BSDEs. More precisely, we are interested in the
correspondence between the properties of the "BSDE" operator and
that of the coefficient.\\
We consider the dynamic operator generated by the maximal solution
of a BSDE:
\begin{definition}
Let $g$ be a standard coefficient. The \index[sub]{$g$-dynamic
operator} $g$-dynamic operator, denoted by
\index[not]{$\mathcal{Y}^{g}$} $\mathcal{Y}^g$, is such that
$\mathcal{Y}^g_t(\xi_T)$ is the maximal solution of the
BSDE$(g,\xi_T)$.
\end{definition}
As a consequence, the adopted point of view is different from that
of the section dedicated to recalls on BSDEs where the terminal
condition
of the BSDE was fixed.\\
It is easy to deduce properties of the $g$-dynamic operator from
those of the coefficient $g$. The converse is more complex and
this study has been initiated by Peng when considering
$g$-expectations (\cite{Peng97}). Our characterization is based
upon the following lemma:
\begin{lemma}[Coefficient Uniqueness]\label{coefuniq}
Let $g^1$ and $g^2$ be two regular coefficients, such that
uniqueness of solution for the BSDE$(g^1)$ holds. Let $\mathcal{Y}^{g^i}$ be $g^i$-dynamic
operator($i=1,2$). Assume that
$$\forall (T,\xi_T),\quad \>d\mathbb{P}\times dt - a.s.\qquad
\mathcal{Y}^{g^1}_t(\xi_T)=\mathcal{Y}^{g^2}_t(\xi_T).$$ a) If the
coefficients  $g^1$ and $g^2$ simply depend on $t$ and $z$, then
$\>d\mathbb{P}\times dt \>- a.s.\quad\forall
z\quad g^1(t,z)=g^2(t,z)$.\\
b) In the general case, the same identity holds provided the
coefficients are continuous w.r. to $t$.
$$\>d\mathbb{P}\times dt \>- a.s.\quad\forall (y,z)\quad
g^1(t,y,z)=g^2(t,y,z).$$
\end{lemma}
{\bf Proof:} Suppose that both coefficients $g^1$ and $g^2$
generate the same solution $Y$ (but a priori different processes
$Z^1$ and $Z^2$) for the BSDEs $(g^1,\xi_T)$ and $(g^2,\xi_T)$,
for {\bf any} $\xi_T$ in the appropriate space ($L^2$ or
$L^{\infty}$). Given the uniqueness of the decomposition of the
semimartingale $Y$, the martingale parts and the finite variation
processes of the both decompositions of $Y$ are indistinguishable.
In particular, $\int_0^t Z^1_sdW_s=\int_0^t Z^2_sdW_s=\int_0^t
Z_sdW_s,\>a.s.$ and $\int_0^t g^1(s,Y_s,Z^1_s)ds=\int_0^t
g^2(s,Y_s,Z^2_s)ds$. Therefore, $\int_0^t
g^1(s,Y_s,Z_s)ds=\int_0^t g^2(s,Y_s,Z_s)ds$. \\
A priori, these equalities only hold for processes $(Y,Z)$
obtained through BSDEs.\\
{\bf a)} Assume that $g^1$ and $g^2$ do not depend on $y$. Given a
bounded adapted process $Z$, we consider the following locally
bounded semimartingale $U$ as $dU_t =g^1(t,Z_t)dt
-Z_tdW_t;\>U_0=u_0$. $(U,Z)$ is the solution of the
BSDE$(g^1,U_{T\wedge\tau})$ where $\tau$ is a stopping time s.t.
$U_{T\wedge\tau}$ is bounded. By uniqueness, $\int_0^{t\wedge\tau}
g^1(s,Z_s)ds=\int_0^{t\wedge\tau} g^2(s,Z_s)ds$. As shown in b)
below, this equality implies that $g^1(s,Z_s)=g^2(s,Z_s),\> a.s.\>
ds\times d\mathbb{P}$. Thanks to the continuity of $g^1$ and $g^2$
w.r. to $z$, we can only consider denumerable rational $z$ to show
that, with $dt\times d{\mathbb P}$ probability one, for any $z$,
$\>g^1(s,z)=g^2(s,z)$. \\
{\bf b1)} In the general case, given a bounded process $Z$, we
consider {\bf a} solution of the following forward stochastic
differential equation $Y$ as $dY_t =g^1(t,Y_t,Z_t)dt
-Z_tdW_t;\>Y_0=y_0$ and the stopping time $\tau_N$ defined as the
first time, when $|Y|$ crosses the level $N$.\\
The pair of processes $(Y_{\tau_N\wedge t},Z_t{\mathbf
1}_{]0,\tau_N]}(t))$ is solution of the BSDE with bounded terminal
condition $\xi_T=Y_{\tau_N\wedge T}$. Thanks to the previous
observation, the pair of processes $(Y_{\tau_N\wedge
t},Z_t{\mathbf 1}_{]0,\tau_N]}(t))$ is also solution of the
BSDE($g^2,Y_{\tau_N\wedge T}$). Hence, both processes $\int_0^t
g^1(s,Y_s,Z_s)ds=\int_0^t g^2(s,Y_s,Z_s)ds$ are indistinguishable
on $]0,\tau_N\wedge T]$. Since $\tau_N$ goes to infinity with $N$,
the equality holds at any time,
for any bounded process $Z$.\\
{\bf b2)} Assume $g^1(s,y,z)$ and $g^2(s,y,z)$ continuous w.r. to
$s$. Let $z$ be a given vector. Let $Y^z_{t+h}$ be a forward
perturbation of a general solution $Y$, at the level $z$ between
$t$ and $t+h$,
$$Y^{z}_{u}=Y_t+\int_t^{u}g^1(s,Y^{z}_s,z)ds-\int_t^{u}zdW_s \qquad \forall u \in [t,t+h].$$
By assumption, $(Y^{z}_{u},z)$ is also solution of the BSDE($g^2,
Y^{z}_{t+h})$, for $u\in[t,t+h]$, and
$\int_t^{u}g^1(s,Y^{z}_s,z)ds=\int_t^{t}g^2(s,Y^{z}_s,z)ds$.
Hence, by continuity,
$\frac{1}{h}\mathbb{E}\Big[Y^{z,i}_{t+h}-Y_t|{\cal F}_t\Big]$ goes
in $L^1$ to $g^i(t,Y_t,z)$ ($i=1,2$) with $h \rightarrow 0$. Then
$g^1(t,Y_t,z)=g^2(t,Y_t,z)$ for any solution $Y_t$ of the BSDE,
i.e. for any v.a ${\cal F}_t$-measurable. $\square$
\paragraph{Comments} \index[aut]{Peng}Peng (\cite{Peng97} and \cite{Peng03}) and
\index[aut]{Briand}\index[aut]{Coquet}\index[aut]{M\'{e}min}\index[aut]{Peng}Briand et al.
\cite{Briand-Coquet-Memin-Peng} have been among the first to look at the dynamic
operators to deduce local properties through the coefficient $g$ of the associated BSDE,
when considering non-linear expectations. More recently, \index[aut]{Jiang}Jiang has considered the
applications of $g$-expectations in finance in his PhD thesis \cite{Jiang}.\\
ii) In \cite{Briand-Coquet-Memin-Peng}, Briand et al. proved a
more accurate result for $g$ Lipschitz. More precisely, let $g$ be
a standard coefficient such that $\mathbb{P}$-a.s., $t \longmapsto
g(t,y,z)$ is continuous and $g(t,0, 0)\in {\cal S}^2$. Let us fix
$(t,y,z) \in [0,T]\times \mathbb{R} \times \mathbb{R}^d$ and
consider for each $n \in \mathbb{N}^*$, $\{(Y_s^{n}, Z_s^n); s \in
[ t, t_n=t + \frac{1}{n}]\}$ solution of the BSDE $(g,X_n)$ where
the terminal condition $X_{t_n}$ is given by $ X_{t_n} = y + z
\big(W_{t_n} - W_t \big)$. Then for each $(t,y,z)\in [0,T]\times
\mathbb{R}\times \mathbb{R}^d$, we have
$$
L^2 -\lim_{n \to \infty } n \big( Y_t^n - y \big) = g (t,y,z).
$$
Some properties automatically hold for the dynamic operator
$\mathcal{Y}^g$ simply because it is the maximal solution of a
BSDE. Some others can be obtained by imposing conditions on the
coefficient $g$:
\begin{theorem}\label{Theorem verification BSDE-risk
measure} Let $\mathcal{Y}^g$ be the $g$-dynamic operator.\\
a) Then, $\mathcal{Y}^g$ is increasing monotonic (P2+),
time-consistent (P4+) and arbitrage-free (P5).\\
b) Moreover, under the assumptions of Lemma \ref{coefuniq},
\begin{enumerate} \item $\mathcal{Y}^g$ is
conditionally invariant (P6) if and only if for any $t \in [0,T],
z\in {\mathbb R}^n$, $g(t,0,0)=0$. \vspace{-1mm} \item
$\mathcal{Y}^g$ is translation invariant (P3+) if and only if $g$
does not depend on $y$. \vspace{-1mm}\item $\mathcal{Y}^g$ is
homogeneous if and only if $g$ is homogeneous;
\end{enumerate}
c) For properties related to the order, the following implications
simply hold: \vspace{-1mm}
\begin{enumerate}
\item If $g$ is convex, then $\mathcal{Y}^g$ is convex
(P1).\vspace{-1mm} \item If $g^1\leq g^2$ , then
$\mathcal{Y}^{g^1}\leq \mathcal{Y}^{g^2}$.
\end{enumerate}
Therefore, if $g$ is a convex coefficient depending only on $z$,
\index[not]{$\mathcal{R}^{g}$} $\mathcal{R}^g(\xi_T) \equiv
\mathcal{Y}^g(-\xi_T)$ is a dynamic convex risk measure, called
{\emph g-conditional risk measure}. \index[sub]{g-conditional risk
measure}
\end{theorem}
Note that $\mathcal{Y}^g$ is a consistent convex price system and
moreover, if for any $t \in [0,T]$, $g(t,0)=0$, then
$\mathcal{Y}^g$ is a non-linear expectation, called
g-expectation.\vspace{2mm}\\
\textbf{Proof}: a) $\bullet$ The strict version of the comparison
Theorem \ref{comparison Theorem}
leads immediately to both properties (P2+) and (P5).\\
$\bullet$ Up to now, we have defined and considered BSDEs with a
terminal condition at a fixed given time $T$. It is always
possible to consider it as a BSDE with a time horizon $T_H\geq T$,
even if $T_H$ is a bounded stopping time. Obviously, the
coefficient $g$ has to be extended as $g \mathbf{1}_{[0,T]}$
 and the terminal condition
$\xi_{T_H}=\xi_T$. Therefore the solution $Y_t$ is
constant on $[T,T_H]$.\\
To obtain the {\it time-consistency} property (P4), also called
the flow property, we consider three bounded stopping times $S
\leq T \leq U$ and write the solution of the BSDEs as function of
the terminal date. With obvious notations, we want to prove that
$Y_S(T,Y_T(U,\xi_U))=Y_S(U,\xi_U)\quad a.s.$.\\
By simply noticing that:
$$
\begin{array}{llll}
&Y_S(T,Y_T(U,\xi_U))=Y_T(U,\xi_U)+&\int_S^T g(t,Z_t) dt-
\int_S^T Z_tdW_t\vspace{2mm}\\
&= \xi_U+\int_T^U g(t,Z_t) dt- \int_T^U Z_tdW_t &+\int_S^T
g(t,Z_t) dt- \int_S^T Z_tdW_t,
\end{array} $$
the process which is defined as $Y_t(T,Y_T(U,\xi_U))$ on $[0,T]$
and by $Y_t(U,\xi_U)$ on $]T,U]$ is the maximal solution of the
BSDE $(g,\xi_U,U)$. Uniqueness of the maximal
solution implies (P4). \vspace{2mm}\\
b) The three properties b1), b2) and b3) involve the same type of
arguments
to be proved, so we simply present the proof for b2).\\
Let $g_m(t,y,z)=g(t,y+m,z)$. We simply note that
$Y^m_.=Y_.(\xi_T+m)-m$ is the maximal solution of the BSDE
$(g_m,\xi_T)$. The translation invariance property is equivalent
to the indistinguishability of both processes $Y$ and $Y^m$; by
the uniqueness Lemma \ref{coefuniq}, this property is equivalent
to the identity
$$g(t,y,z)=g_m(t,y,z)=g(t,y+m,z) \qquad a.s.$$
This implies that $g$ does not depend on $y$.\vspace{2mm}\\
c) $\bullet$ c1) For the {\it convexity} property, we consider
different BSDEs: $(Y_t^1,Z_t^1)$ is the (maximal) solution of
$(g,\xi_T^1)$ and $(Y_t^2,Z_t^2)$ is the (maximal) solution of
$(g,\xi_T^2)$. Then, we look at $\widetilde{Y}_t=\lambda
Y_t^1+(1-\lambda)Y_t^2$, with $\lambda \in [0,1]$. We have:
$$-d\widetilde{Y}_t=(\lambda g(t,Y_t^1,Z_t^1)+(1-\lambda)g(t,Y_t^2,Z_t^2))dt-
(\lambda Z_t^1+(1-\lambda)Z_t^2)dW_t \quad ; \quad
\widetilde{Y}_T=\lambda \xi_T^1+(1-\lambda)\xi_T^2.$$ Since $g$ is
convex, we can rewrite this BSDE as:
$$-d\widetilde{Y}_t=(g(t,\widetilde{Y}_t,\widetilde{Z}_t)+
\alpha(t,Y_t^1,Y_t^2,Z_t^1,Z_t^2,\lambda))dt-\widetilde{Z}_tdW_t$$
where $\alpha$ is a a.s. non-negative process. Hence, using the
comparison theorem, the solution $\widetilde{Y}_t$ of this BSDE is
for any $t \in [0,T]$ a.s. greater than the solution $Y_t$ of the
BSDE $(g,\lambda \xi_T^1+(1-\lambda)\xi_T^2)$. It is a
super-solution in the sense of Definition 2.1 of \index[aut]{El
Karoui} El Karoui, \index[aut]{Peng} Peng and \index[aut]{Quenez}
Quenez
\cite{EPQ}.\\
$\bullet$ c2) is a direct consequence of the comparison Theorem
\ref{comparison Theorem}. $\qquad \square$ \vspace{-2mm}
\paragraph{Some additional comments on the relationship between BSDE and
dynamic operators} Since 1995, Peng has focused on finding
conditions on dynamic operators so that they are linear growth
g-expectations. This difficult problem is solved in particular for
dynamic operators satisfying a domination assumption introduced by
\index[aut]{Peng} Peng \cite{Peng97} in 1997 where
$b_{k}(z)=k|z|$. For more details, please refer to his lecture
notes on BSDEs and dynamic operators \cite{Peng03}.
\begin{theorem}
Let $({\cal E}_t;0 \leq t \leq T)$ be a non-linear expectation such that:\\
There exists $|\lambda| \in {\cal H}^2$ and a sufficiently large
real number
$k >0$ such that for any $t \in [0,T]$ and any $\xi_T \in L^2({\cal F}_T)$:\\
$${\cal E}_t^{-b_{k}+|\lambda|}(\xi_T) \leq {\cal E}_t(\xi_T) \leq {\cal E}_t^{b_{k}+|\lambda|}(\xi_T)\quad a.s. $$
and for any $(\xi^1_T,\xi^2_T) \in L^2({\cal F}_T)$: ${\cal
E}_t(\xi^1_T) -{\cal E}_t(\xi^2_T) \leq {\cal
E}_t^{b_{k}}(\xi^1_T-\xi^2_T)$ Then, there exists a function
$g(t,y,z)$ satisfying assumption (H1) such that for any $t \in
[0,T]$,
$$\forall \; \xi_T \in L^2({\cal F}_T),\quad a.s.,\quad \forall t,\quad
{\cal E}_t(\xi_T)={\cal E}_t^g(\xi_T)$$
\end{theorem}
For a proof of this theorem, please refer to \index[aut]{Peng}
Peng \cite{Peng99}.
\paragraph{Infinitesimal Risk Management}
The coefficient of any $g$-conditional risk measure
$\mathcal{R}^g$ can be naturally interpreted as the {\it
infinitesimal risk measure} \index[sub]{infinitesimal risk
measure} over a time interval $[t,t+dt]$ as:
$$\mathbb{E}_{\mathbb{P}}[d\mathcal{R}^g_t | \mathcal{F}_t]=-g(t,Z_t)dt,$$
where $Z_t$ is the local volatility of the $g$conditional risk measure.\\
Therefore, choosing carefully the coefficient $g$ enables to
generate $g$-conditional risk measures that are locally compatible
with the views and practice of the different agents in the market.
In other words, knowing the infinitesimal measure of risk used by
the agents is enough to generate a dynamic risk measure, locally
compatible. In this sense, the $g$-conditional risk measure
may appear more tractable than static risk measures.\\
The following example gives a good intuition of this idea: the
$g$-conditional risk measure corresponding to the mean-variance
paradigm has a $g$-coefficient of the type $g(t,z)=-\lambda_t z +
\frac{1}{2} z^2$. The process $\lambda_t$ can be interpreted as
the correlation with the market
num\'{e}raire.\vspace{2mm}\\
Therefore, $g$-Conditional risk measures are a way to construct a
wide family of convex risk measures on a probability space with
Brownian filtration, taking into account the ability to decompose
the risk through inter-temporal local risk measures $g(t,Z_t)$. \vspace{2mm}\\
In the following, to study $g$-conditional risk measures, we adopt
the same methodology as in the static framework. In particular, we
start by developing a dual representation for these dynamic risk
measures, in terms of the "dual function" of their coefficient.
This study requires some general properties of convex functions on
$\mathbb{R}^n$.
\section{Dual Representation of $g$-Conditional Risk Measures}
\label{subsection dual representation risk measures}
Following the approach adopted in the first part of this chapter
when studying static risk measures, we now focus on a dual
representation for $g$-conditional risk measures. The main tool is
the {\it Legendre-Fenchel transform}  \index[not]{$G$} $G$ of the
coefficient $g$, defined by:
\begin{equation}\label{polar}
G(t,\mu)=\sup_{z\in \mathbf{Q_{\text{rational}}^n}}\left\{\langle \mu,-z \rangle -g(t,z)\right\}.
\end{equation}
The convex function $G$ is also called the {\it polar function} or
the {\it conjugate} of $g$. Provided that $g$ is continuous,
\begin{equation}\label{bipolar}
g(t,z)=\sup_{\mu \mathbf{Q_{\text{rational}}^n}}\big(\langle
\mu,-z \rangle -G(t,\mu)\big).
\end{equation}
More precisely,
\begin{definition} A g-conditional risk ${\cal R}^g$ measure is said
to have a {\it dual representation} if there exists a set
$\mathcal{A}$ of admissible controls such that for any bounded
stopping time $S \leq T$ and any $\xi_T$ in the appropriate space
\begin{equation}
\label{eq:dual} {\cal R}^g_S(\xi_T)={{\rm ess}\sup}_{\mu \in {\cal
A}} {\mathbb E}_{{\mathbb Q}^{\mu}}\Big[-\xi_T-\int_S^T
G(t,\mu_t)dt \big|{\cal F}_S\Big]
\end{equation}
where ${\mathbb Q}^{\mu}$ is a probability measure absolutely continuous with
respect to $\mathbb P$.\\
The dual representation is said to be exact at $\bar{\mu}$ if the ${{\rm ess}\sup}$ is reached for $\bar{\mu}$.
\end{definition}
In order to obtain this representation, several intermediate steps
are needed:
\begin{enumerate}
\item Refine results on Girsanov theorem and the integrability
properties of martingales with respect to change of probability
measures.
\item Refine results from convex analysis on
the Legendre-Fenchel transform and the existence of an optimal
control in both Formulae \eqref{polar} and \eqref{bipolar},
including measurability properties,
\end{enumerate}
The next paragraph gives a summary of the main results that are
needed.
\subsection{Girsanov Theorem and BMO-Martingales}
 Our main reference on Girsanov theorem and BMO-martingales is the book by
\index[aut]{Kazamaki} Kazamaki \cite{Kazamaki}. The exponential
martingale associated with the $d$-dimensional Brownian motion $W$, ${\cal E}(\int_0^t\mu_s
dW_s)=\Gamma_t^\mu=\exp\Big(\int_0^t\mu_s
dW_s-\frac{1}{2}|\mu_s|^2ds\Big)$, solution of the forward
stochastic equation
\begin{equation}
\label{eq:density} d\Gamma^\mu_t=\Gamma^\mu_t\mu_t^*dW_t\>\>,\quad
\Gamma^\mu_0=1
\end{equation}
is a positive local martingale, if $\mu$ is an adapted process
such that
$\int_0^T|\mu_s|^2ds<\infty$.\\
When $\Gamma^\mu$ is a uniformly integrable \index[not]{${\rm
u.i.}$} (u.i.) martingale, $\Gamma^\mu_T$ is the density (w.r. to
$\mathbb P$) of a new probability measure denoted by ${\mathbb
Q}^\mu$. Moreover, if $W$ is a ${\mathbb P}$-Brownian motion, then
$W^\mu_t=W_t-\int_0^t\mu_sds$
is a ${\mathbb Q}^\mu$-Brownian motion.\vspace{1mm}\\
Questions around Girsanov theorem are of two main types. They
mainly consist of:
\begin{itemize}
\item first, finding conditions on $\mu$ so that $\Gamma^{\mu}$ is
a u.i. martingale. \item second, giving so;e precision on the
integrability properties that are preserved under the new
probability measure.
\end{itemize}
The bounded case, that is recalled below, is well-known. The BMO
case is less standard, so we give more details.
\subsubsection{Change of Probability Measures with Bounded Coefficient}
When $\mu$ is bounded, it is well-known that
the exponential martingale belongs to all ${\cal H}^p$-spaces.\\
Moreover, if a process is in $\mathcal{H}^2(\mathbb{P})$, it is in
$\mathcal{H}^{1+\epsilon}(\mathbb{Q}^{\mu})$. In particular, if
$M^Z_t=\int_0^t Z_sdW_s$ is a ${\cal H}^2(\mathbb P)$-martingale,
then $\widehat M^Z_t=\int_0^tZ_sdW^\mu_s$ is a u.i. martingale
under $\mathbb{Q}^{\mu}$, with null
$\mathbb{Q}^{\mu}$-expectation.
\subsubsection{Change of Probability Measures with BMO-Martingale}
The right extension of the space of bounded processes is the space
of BMO processes defined as: \index[not]{${\rm BMO}(\mathbb P)$}
$${\rm BMO}(\mathbb P)=\{\varphi\in{\cal
H}^2\quad s.t \quad \exists{ C}\quad \forall t\quad {\mathbb
E}\Big[\int_t^T|\varphi_s|^2ds|{\cal F}_t\Big]\leq C \>\>a.s.\}$$ The
smallest constant $C$ such that the previous inequality holds is
denoted by $C^*=||\varphi||^2_{\rm BMO}$.\\
In terms of martingale, the stochastic integral
$\int_0^t\varphi_sdW_s$ is said to be a {\rm BMO($\mathbb
P$)}-martingale if and only if the process $\varphi$ belongs to
{\rm BMO($\mathbb P$)}. The following deep result is proved in
\index[aut]{Kamazaki} Kazamaki \cite{Kazamaki} (Section 3.3).
\begin{theorem} \label{th:kazamaki}
Let the adapted process $\mu$ be in {\rm BMO($\mathbb P$)}. Then
\begin{enumerate}
\item The exponential martingale $\Gamma^\mu$ is a u.i. martingale
and defines a new equivalent probability measure ${\mathbb
Q}^\mu$. Moreover, $W^\mu_t=W_t-\int_0^t\mu_sds$ is a ${\mathbb
Q}^\mu$-Brownian motion. \item $M^\mu_t=\int_0^t\mu_s^*dW_s$, and
more generally any {\rm BMO($\mathbb P$)}-martingale
$M^Z_t=\int_0^tZ_sdW_s$, are transformed into continuous processes
$\widehat M^\mu_t=\int_0^t\mu_s^*dW^\mu_s$ and $\widehat
M^Z_t=\int_0^tZ_sdW^\mu_s$ that are {\rm BMO(${\mathbb
Q}^\mu$)}-martingales. \item The {\rm BMO}-norms with respect to
${\mathbb P}$ and ${\mathbb Q}^\mu$ are equivalent:
$$k||Z||_{{\rm BMO({\mathbb Q}^\mu)}} \leq ||Z||_{{\rm BMO({\mathbb P})}}\leq K
||Z||_{{\rm BMO({\mathbb Q}^\mu)}}.$$ The constants $k$ and $K$
only depend on the {\rm BMO}-norm of $\mu$.
\end{enumerate}
\end{theorem}
\index[aut]{Hu} Hu, \index[aut]{Imkeller} Imkeller and
\index[aut]{M\"{u}ller} M\"{u}ller \cite{HuImkellerMuller} were
amongst the first to use the property that the martingale
$dM^Z_t=Z_tdW_t$ which naturally appears in BSDEs associated with
exponential hedging problems, is BMO. Since then, such a property
has been used in different papers, mostly dealing with the
question of dynamic hedging in an exponential utility framework
(see for instance the recent papers by \index[aut]{Mania} Mania,
\index[aut]{Santacroce} Santacroce and \index[aut]{Tevzadze}
Tevzadze \cite{ManiaSantacroceTevradze} and \index[aut]{Mania}
Mania and
\index[aut]{Schweizer} Schweizer \cite{ManiaSchweizer}).\\
In the proposition below, we extend their results to general
quadratic BSDEs.
\begin{proposition}
\label{propo:MZMartBMO} Let (Y,Z) be the maximal solution of the
quadratic (H3) BSDE with coefficient $g$, and
$M^Z=\int_0^.Z_sdW_s$ the stochastic integral $Z.W$
$$dY_t=g(t,Z_t)dt -dM^Z_t,\quad Y_T=\xi_T.$$
Given that by assumption $Y$ is bounded, and $|g(t,0)|^{1/2}\in
{\rm BMO({\mathbb P})}$, $M^Z$ is a ${\rm BMO({\mathbb
P})}$-martingale
\end{proposition}
{\bf Proof:} Let $k$ be the constant such that $|g(t,z)|\leq|g(t,0)|+k|z|^2$.\\
Thanks to Itô's formula applied to the solution $(Y,Z)$ and to the
exponential function:
\begin{eqnarray*}
\exp(\beta\,Y_t)&=&\exp(\beta\,Y_T)+\beta\int_t^T\,
\exp(\beta\,Y_s)g(s,Z_s)ds-\frac{\beta^2}{2}\int_t^T\,
\exp(\beta\,Y_s)|Z_s|^2ds\\
&-&\beta\int_t^T\exp(\beta\,Y_s)Z_sdW_s\\
&=&\exp(\beta\,Y_T)+\beta\int_t^T\,
\exp(\beta\,Y_s)\Big(g(s,Z_s)-\frac{\beta}{2}|Z_s|^2\Big)ds
-\beta\int_t^T\exp(\beta\,Y_s)Z_sdW_s.
\end{eqnarray*}
Given that $\quad\frac{\beta}{2}|Z_s|^2-g(s,Z_s)\geq
(\frac{\beta}{2}-k)|Z_s|^2-|g(s,0)|\geq \varepsilon
|Z_s|^2-|g(s,0)|$ for $\beta\geq (k+\varepsilon)\>$ and taking the
conditional expected value, we obtain:
$$\beta \>\varepsilon\>
{\mathbb E}\Big[\int_t^T\exp(\beta\,Y_s)|Z_s|^2ds|{\cal
F}_t\Big]\leq C+\beta {\mathbb
E}\Big[\int_t^T\exp(\beta\,Y_s)|g(s,0)|ds|{\cal F}_t\Big]\leq C $$
where $C$ is a universal constant that may change from place to
place. Since $\exp(\beta\,Y_s)$ is bounded both from below and
from above, the property holds.$\qquad \square$
\subsection{Some Results in Convex Analysis \label{Subsection results convex analysis}}
Some key results in convex analysis are needed to obtain the dual
representation of $g$-conditional risk measures. They are
presented in the Appendix \ref{Appendix convex analysis} to
preserve the continuity of the arguments in this part. More
details or proofs may be found in \index[aut]{Aubin} Aubin
\cite{Aubin}, \index[aut]{Hiriart-Urruty} Hiriart-Urruty and
\index[aut]{Lemar\'{e}chal} Lemar\'{e}chal
\cite{Hiriart-Urruty-Lemarchal} or \index[aut]{Rockafellar}
Rockafellar \cite{Rockafellar}.
\subsection{Dual Representation of Risk Measures}
We now study the dual representation of $g$-conditional risk
measures. The space of admissible controls depends on the
assumption imposed on the coefficient $g$. We consider
successively both situations (H1) and (H3). There is no need to
look separately at (H2), as, under our assumptions, the condition
(H2) implies the condition (H1) (for more details, please refer to
the Appendix \ref{par:recession function}). The (H1) case has been
solved in \cite{EPQ} but the (H3) case is new.
\begin{theorem} \label{theorem dual representation}
Let $g$ be a convex coefficient satisfying (H1) or (H3) and $G$ be
the associated polar process, $ G(t,\mu)=\sup_{z\in
\mathbf{Q_{\text{rational}}^n}}\left\{\langle \mu,-z \rangle
-g(t,z)\right\}$.\\
$i)$ For almost all $(\omega,t)$, the program $g(\omega,t, z )
=\sup_{\mu\in \mathbf{Q_{\text{rational}}^n}}[\langle
\mu,-z\rangle - G(\omega, t,\mu )\,]$ has an optimal progressively
measurable solution ${\bar \mu}(\omega, t)$ in the subdifferential of $g$ at $z$, $\partial g(\omega,t,z)$.\\
$ii)$  Then ${\cal R}^g$ has the following dual representation,
exact at ${\bar \mu}$,
$$\mathcal{R}^g _t(\xi_T)={\rm esssup}_{\mu \in {\cal A}}{
\mathbb E}_{{\mathbb Q}^{\mu}}\Big[-\xi_T- \int_t^T\>G(s,
\mu_s)ds\big|{\cal F}_t\Big]={ \mathbb E}_{{\mathbb Q}^{\bar
\mu}}\Big[-\xi_T- \int_t^T\>G(s,{\bar \mu}_s)ds\big|{\cal
F}_t\Big]$$ where:
\begin{enumerate}
 \item Under (H1) ($|g(t,z)|\leq |g(t,0)||+k|z|$), ${\cal A}$
is the space of adapted processes $\mu$ bounded by $k$, and
${\mathbb Q}^\mu$ is the associated equivalent probability measure
with density $\Gamma^\mu_T$ where $\Gamma^\mu$ is the exponential
martingale defined in \eqref{eq:density}. \item Under (H3),
($|g(t,z)|\leq |g(t,0)|+k|z|^2$), ${\cal A}$ is the space of
BMO($\mathbb{P}$)-processes $\mu$ and  ${\mathbb Q}^\mu$ is
defined as above.
\end{enumerate}
$iii)$ Let $g(t,.)$ be a strongly convex function (i.e.
$g(t,z)-\frac{1}{2}C|z|^2$ is a convex function). Then the
Fenchel-Legendre transform $G(t,\mu)$ has a quadratic growth in
$\mu$ and the following dual representation holds true:
$${\mathbb E}_{{\mathbb Q}^{\mu}}\Big[\int_t^T
G(s,\mu_s)ds \Big| {\cal F}_t \Big]= {\rm esssup}_{\xi_T}{\mathbb
E}_{{\mathbb Q}^{\mu}}\big[\xi_T \big| {\cal F}_t \big]-{\cal
R}_t^g(\xi_T)={\mathbb E}_{{\mathbb Q}^{\mu}}\big[\bar{\xi}_T
\big| {\cal F}_t \big]-{\cal R}_t^g(\bar{\xi}_T)$$
\end{theorem}
{\bf Proof:} $i)$ Since $g$ is a proper function, the dual
representation of $g$ with its polar function $G$ is exact at
$\bar{\mu} \in \partial g(z)$: $$g(t, z ) =\sup_{\mu\in
\mathbf{Q_{\text{rational}}^n}}[\langle \mu,-z\rangle - G(t,\mu
)\,]=\langle \bar{\mu},-z \rangle - G(t,\bar{\mu}),$$ using
classical results of convex analysis, recalled in the Appendix
\ref{Appendix convex analysis}.\\
The measurability of $\bar{\mu}$ is separately studied in Lemma
\ref{Lemma measurability of mu-bar}
just after this proof.\vspace{2mm}\\
$ii)$ a) Let us first consider a coefficient $g$ with linear
growth (H1); so, $g(t,0)$ is in $\mathcal{H}^2$. By definition,
$-G(t,\mu_t)$ is dominated from above by the square integrable
process $g(t,0)$. Then, let $\mathcal{R}_t^g (\xi_T):=Y_t$ be the
solution of the BSDE $(g,-\xi_T)$,
\begin{equation}
\label{eq:mubsde}
-dY_t=g(t,Z_t)dt-Z_tdW_t=(g(t,Z_t)-\langle \mu_t,-Z_t \rangle)dt-Z_tdW^{\mu}_t,\quad
Y_T=-\xi_T.
\end{equation}
By Girsanov Theorem (Theorem \ref{th:kazamaki}), for $\mu \in
{\cal A}$ the exponential martingale $\Gamma^\mu$ is u.i. and
defines a probability measure  ${\mathbb Q}^\mu$ on ${\cal F}_T$
such that the process $W^{\mu}=W-\int_0^.\mu_sds$ is a ${\mathbb
Q}^\mu$-Brownian motion. Moreover, since $M^Z=\int_0^.Z_sdW_s$ is
in ${\cal H}^2({\mathbb P})$, $\widehat M^Z=\int_0^.Z_sdW^\mu_s$
is a u.i. ${\mathbb Q}^{\mu}$-martingale. Moreover, since $\mu$ is
bounded and $g$ uniformly Lipschitz, the process
$(g(t,Z_t)-Z_t\mu_t)$ belongs to ${\cal H}^2({\mathbb P})$ but
also to ${\cal H}^{1+\epsilon}({\mathbb Q}^\mu)$. So we can use an
integral representation of the BSDE \eqref{eq:mubsde} in terms of
\begin{equation}\label{eq:ineqduality}
Y_t={ \mathbb E}_{{\mathbb Q}^{\mu}}\Big[-\xi_T+
\int_t^T(g(s,Z_s)-\langle \mu_s,-Z_s \rangle)ds\big|{\cal F}_t\Big] \geq { \mathbb
E}_{{\mathbb Q}^{\mu}}\Big[-\xi_T-
\int_t^T\>G(s,\mu_s)ds\big|{\cal F}_t\Big].
\end{equation}
We do not need to prove that the last term is finite. It is enough
to recall that $(-G(s,\mu_s))^+$ is dominated from above by the
$d{\mathbb Q}\times ds$ integrable process $(g(s,0))^+$.\vspace{1mm}\\
b) Let ${\bar \mu}$ be an optimal control, bounded by $k$, such
that $g(t,Z_t)=\langle {\bar \mu}_t,-Z_t \rangle -G(t,{\bar
\mu}_t)$ (see Lemma \ref{Lemma measurability of mu-bar} for
measurability results). Then the process $-G(t,{\bar \mu}_t)$
belongs to ${\cal H}^2({\mathbb P})$ and so to ${\cal
H}^{1+\epsilon}({\mathbb Q}^{\bar \mu})$. By the previous result,
$Y_t={ \mathbb E}_{{\mathbb Q}^{\bar \mu}}\Big[-\xi_T-
\int_t^T\>G(s,{\bar \mu}_s)ds\big|{\cal F}_t\Big].$ So the process
$Y$ is the value function of the maximization dual problem
$Y_t={\rm esssup}_{\mu \in {\cal A}}{ \mathbb E}_{{\mathbb
Q}^{\mu}}\Big[-\xi_T- \int_t^T\>G(s,
\mu_s)ds\big|{\cal F}_t\Big]$.\vspace{2mm}\\
c) We now consider a coefficient $g$ with quadratic growth (H3)
and bounded solution $Y_t$. Using the same notation, we know by
Girsanov Theorem \ref{th:kazamaki} that if $\mu \in {\rm
BMO({\mathbb P})}$, $\Gamma^\mu$ is a u.i. martingale and the
probability measure ${\mathbb Q}^\mu$ is well-defined. The proof
of the dual representation is very similar to that of the previous
case, after solving some integrability questions. It is enough to
notice that
\begin{itemize}
\item by assumption, $|g(.,0)|^{\frac{1}{2}}$ is ${\rm BMO({\mathbb P})}$,
\item by Proposition \ref{propo:MZMartBMO}, $Z$ is ${\rm BMO({\mathbb P})}$,
\item by Girsanov Theorem \ref{th:kazamaki}, for any $\mu\in {\rm BMO({\mathbb P})}$, the processes $\mu$, $Z$
and $|g(.,0)|^{\frac{1}{2}}$ are in ${\rm BMO({\mathbb
Q}^{\mu})}$.
\end{itemize}
So $|g(t,Z_t)|^{\frac{1}{2}}$ and $|\mu_t Z_t|^{\frac{1}{2}}$ are
in ${\rm BMO({\mathbb Q}^{\mu})}$. Moreover, the process
$(-G(t,\mu))^+$ which is dominated from above by $|g(t,0)|$
is a ${\mathbb Q}^{\mu}\times dt$-integrable process. Then the inequality \eqref{eq:ineqduality} holds.\\
d) Let ${\bar \mu}$ be an optimal control, such that
$g(t,Z_t)=\langle {\bar \mu}_t,-Z_t \rangle -G(t,{\bar \mu}_t)$.
Given that $g(t,.)$ has quadratic growth, the polar function
$G(t,.)$ satisfies the following inequality, $G(t,{\bar
\mu}_t)\geq -|g(t,0)|+\frac{1}{4k}|{\bar \mu}_t|^2$. Then, for
small $\varepsilon<\frac{1}{4k}$,
$(\frac{1}{4k}-\varepsilon)|{\bar \mu}_t|^2 \leq G(t,{\bar
\mu}_t)+|g(t,0)|-\varepsilon|{\bar \mu}_t|^2 \leq
|g(t,0)|-g(t,Z_t)+\langle {\bar \mu}_t,-Z_t \rangle
-\varepsilon|{\bar \mu}_t|^2 \leq
|g(t,0)|-g(t,Z_t)+\frac{1}{4\varepsilon}|Z_t|^2\>$. Since both
processes $|g(t,Z_t)|^{1/2}$ and $Z$ are ${\rm BMO({\mathbb P})}$,
${\bar \mu}$ is also ${\rm BMO({\mathbb P})}$, and the other
processes hold nice integrability properties with respect to both
probability measures $\mathbb P$ and ${\mathbb
Q}^{\mu}$ and the integral representation follows.\vspace{1mm}\\
$iii)$ Let $h(t,z)=g(t,z)-\frac{1}{2}C|z|^2$ be the convex
function associated with $g$. Since $g$ is the sum of two convex
functions $h$ and $\frac{1}{2}C|.|^2$, its Fenchel-Legendre
transform $G$ is the inf-convolution of the Fenchel-Legendre
transforms of both $h$ and $\frac{1}{2}C|.|^2$. But the
Fenchel-Legendre transform of the quadratic function
$\frac{1}{2}C|.|^2$ is still a quadratic function,
$\frac{1}{2C}|\mu|^2$ and $G$ has a quadratic growth (as the
inf-convolution of a convex function $H$ with a quadratic
function). Therefore, for a given $\mu \in \> {\rm
BMO(\mathbb{P})}$, there exists $\bar{Z} \in \>{\rm
BMO(\mathbb{P})}$ such that $G(t,\mu_t)=\langle \mu_t,-\bar{Z}_t
\rangle -g(t,\bar{Z}_t)$ (in other words, $\mu \in
\partial(\bar{Z})$).\\
We now introduce the penalty function $\alpha^{\mu}$ defined by
$\alpha^{\mu}_t={\mathbb E}_{{\mathbb Q}^{\mu}}\Big[\int_t^T
G(s,\mu_s)ds \Big| {\cal F}_t \Big]$. Using the above duality
result, we have:
$$\alpha^{\mu}_t={\mathbb E}_{{\mathbb Q}^{\mu}}[\int_t^T(\langle
\mu_s,-\bar{Z}_s \rangle -g(s,\bar{Z}_s))ds \Big| {\cal F}_t
\Big].$$ Since $\bar{\xi}_T=\int_0^T (\langle \mu_s,-\bar{Z}_s
\rangle -g(s,\bar{Z}_s))ds + \int_0^T \bar{Z}_s dW_s^{\mu}$ and
$\mathcal{R}^g_t(\bar{\xi}_T)=\int_0^t (\langle \mu_s,-\bar{Z}_s
\rangle -g(s,\bar{Z}_s))ds + \int_0^t \bar{Z}_s dW_s^{\mu}$, we
finally deduce that:
$$\alpha^{\mu}_t={\mathbb E}_{{\mathbb Q}^{\mu}}\big[\bar{\xi}_T
\big| {\cal F}_t \big]-{\cal R}_t^g(\bar{\xi}_T).$$ Moreover,
using Equation (\ref{eq:ineqduality}), we have:
$$\alpha^{\mu}_S \geq {\rm ess}\sup_{\xi_T}{\mathbb
E}_{{\mathbb Q}^{\mu}}\big[\xi_T \big| {\cal F}_S \big]-{\cal
R}^g(\xi_T).$$ Hence, the result. $\qquad
\square$ \vspace{1mm}\\
The question of the measurability of the optimal solution(s)
$\bar{\mu}$ is considered in the following lemma.
\begin{lemma}\label{Lemma measurability of mu-bar}
Let $g$ be a convex coefficient satisfying (H1) or (H3) and $G$ be the associated polar
function. There exists an progressively measurable optimal
solution ${\bar \mu}$ such that $g(t,Z_t )=\langle {\bar \mu}_t,-Z_t \rangle  -
G(t,{\bar \mu}_t ) \quad a.s.\> d\mathbb{P}\times dt$.
\end{lemma}
{\bf Proof:} For each $(\omega,t) \in \Omega \times [0,T]$, the
sets given by: $\{\mu \in {\mathbb R}^n\>:\> g(\omega, t,Z_t) =
Z_t\mu- G(\omega, t,\mu)\}$ are nonempty. Hence, by a measurable
selection theorem (see for instance \index[aut]{Dellacherie}
Dellacherie and \index[aut]{Meyer} Meyer \cite{dm2} or
\index[aut]{Benes} Benes \cite{Benes}), there exists a ${\mathbb
R}^n$-valued progressively measurable process ${\bar \mu}$ such
that: $g(\omega, t,Z_t )=\langle {\bar \mu}_t,-Z_t \rangle  -
G(\omega, t,{\bar \mu}_t) \quad d\mathbb{P}\times dt-a.s.$. $
\quad \square$
\subsection{$g$-Conditional $\gamma$-Tolerant Risk Measures
and Asymptotics}
In this subsection, we pursue our presentation and study of
g-conditional risk measures using an approach similar to that we
have adopted in the static framework.
\subsubsection{$g$-Conditional $\gamma$-Tolerant Risk Measures}
As in the static framework, we can define dynamic versions for
both coherent and $\gamma$-tolerant risk measures based on the
properties of their coefficients using the uniqueness Lemma \ref{coefuniq}.\\
More precisely, let $\gamma>0$ be a risk-tolerance coefficient. As
in the static framework, where the $\gamma$-dilated of any static
convex risk measure $\rho$ is defined by $\rho_{\gamma}(\xi_T) =
\gamma\rho\Big(\frac{1}{\gamma}\xi_T\Big)$we can define the
g-conditional risk measure, \index[not]{$\mathcal{R}^g_{\gamma}$}
${\cal R}^g_{\gamma}$, {\it $\gamma$-tolerant of $\mathcal{R}^g$},
as the risk measure associated with the coefficient $g_{\gamma}$,
which is the $\gamma$-dilated of $g$: $g_{\gamma}(t,z)=\gamma
g(\frac{1}{t,\gamma}z)$.\\
Note that if $g$ is Lipschitz continuous (H1), $g_{\gamma}$ also
satisfies (H1), and if $g$ is continuous with quadratic growth
(H3) with parameter $k$, then $g$ also satisfies (H3), but with
parameter $\frac{k}{\gamma}$. Note also that the dual function of
$g_{\gamma}$, $G_{\gamma}$, can be expressed in terms of $G$, the dual function of $g$ as
$G_{\gamma}(\mu)=\gamma G(\mu)$.\\
A standard example of $g$-conditional $\gamma$-tolerant risk
measure is certainly the dynamic entropic risk measure
$e_{\gamma,t}(\xi_T)=\gamma\ln\mathbb{E}\big[\exp(-\frac{1}{\gamma}\xi_T)|
\mathcal{F}_t\big]$, which is the $\gamma$-tolerant of
$e_{1,t}$.\vspace{-2mm}
\paragraph{Asymptotic behaviour of entropic risk measure}
Let us look more closely at the dynamic entropic risk measure.
Letting $\gamma$ go to $+\infty$, the BSDE-coefficient
$q_{\gamma}(z)=\frac{1}{2\gamma}|z|^2$ tends to $0$ and we
directly obtain the natural extension of the static case,
$e_{\infty,t}(\xi_T)=\mathbb{E}_{\mathbb{P}}[-\xi_T |
\mathcal{F}_t]$.\\
Letting $\gamma$ tend to $0$, the BSDE coefficient explodes if
$|z|\neq 0$ and intuitively the martingale of this BSDE has to be
equal to $0$. More precisely, since by definition
$\exp(e_{\gamma,t}(\xi_T))=
\mathbb{E}\big[\exp(-\frac{1}{\gamma}\xi_T)|
\mathcal{F}_t\big]^{\gamma}$, $\lim_{\gamma\rightarrow
0}\exp(e_{\gamma,t}(\xi_T))=||\exp(-\xi_T)||^{\infty}_t=\inf\left\{Y
\in \mathcal{F}_t:Y_t \geq \exp(-\xi_T)\right\}$. So we have
$e_{0^+,t}(\xi_T)=||-\xi_T||^{\infty}_t$. This conditional risk
measure is a $g$-conditional risk measure associated with the
indicator function of $\{0\}$. Let us also observe that
$e_{0^+,t}(\xi_T)$ is an adapted non-increasing process without
martingale part.
\subsubsection{Marginal Risk Measure}
In the general $\gamma$-tolerant case, assuming that the
$g$-conditional risk measures are centered (equivalently
$g(t,0)=0$ equivalently $G(t,.)\geq 0$), the same type of results
can be obtained concerning the asymptotic behavior of the
$\gamma$-dilated coefficient and the duality. Then, the limit of
$g_{\gamma}$ when ${\gamma}\rightarrow +\infty$ is the derivative
of $g$ at the
origin in the direction of $z$.\\
\index[not]{$\mathcal{R}^g_{\infty}$}${\cal R}^g_{\infty}$ is the
non-increasing limit of ${\cal R}^g_{\gamma}$ defined by its dual
representation ${\cal R}^g_{\infty,t}(\xi_T)={\rm ess\sup}_{\mu
\in {\cal A}}\{{\mathbb E}_{{\mathbb
Q}^{\mu}}\big[-\xi_T\big|{\cal F}_t\big]\big|\>G(u,\mu_u)=0,\>
\forall u\geq t,\>-a.s.\}$ ; in some cases (in particular, in the
quadratic case when the polar function $G$ has a unique $0$, i.e.
$G(u,0)=0$ is unique), $-{\cal R}^g_{\infty}$ is a linear pricing
rule and can be seen as an extension of the Davis price (see
\index[aut]{Davis} Davis \cite{Davis97}).
\subsubsection{Conservative Risk Measures and Super Pricing \label{Subsection dyn conservative risk measures}}
We now focus on the properties of the $g$-conditional
$\gamma$-tolerant risk measures when the risk tolerance
coefficient goes to zero. To do so, we need some results in convex
analysis regarding the so-called \index[sub]{recession function}
{\it recession function}, defined for any $z \in {\rm Dom(g)}$ by
\index[not]{$g_{0^+}$} $g_{0^+}(z) :=\lim_{\gamma \downarrow
0}\gamma g\big(\frac{1}{\gamma}z\big)= \lim_{\gamma \downarrow
0}\gamma \big(g(y+\frac{1}{\gamma}z)-g(y)\big)$. The key
properties of this function are recalled in the Appendix
\ref{par:recession function}.
\paragraph{\sc \index[sub]{conservative risk measure}Conservative Risk Measure}
$\bullet$ Under assumption $(H1)$, we may assume that $g(t,0)=0$.
Therefore, the polar function $G$ is non negative. Since $g(t,.)$
has a linear growth with constant $k$, the recession function
$g_{0^+}(t,.)$ is finite everywhere with linear growth, and the
domain of the dual function $G$ is bounded by $k$. The
BSDE($g_{0^+},\xi_T$) has a unique solution $Y^{0}_t(\xi_T)\geq
{\cal R}^{g_{\gamma}}_t(\xi_T)$. Using their dual representation
through their polar functions $l_{{\rm Dom(G)}}$ and $\gamma G$,
$$
\begin{array}{ll}
Y^{0}_t(\xi_T)&= {\rm ess\sup}_{\mu \in {\cal A}_k}{\mathbb
E}_{{\mathbb
Q}^{\mu}}\big[-\xi_T-\int_t^T l_{\rm Dom(G)}(u,\mu_u)du\big|{\cal F}_t\big],\\
{\cal R}^g_{\gamma,t}(\xi_T)&={\rm ess\sup}_{\mu \in {\cal
A}_k}{\mathbb E}_{{\mathbb Q}^{\mu}}\big[-\xi_T-\gamma \int_t^T
G(u,\mu_u)du\big|{\cal F}_t\big].
\end{array}
$$
we can take the non decreasing limit in the second line and show
that \index[not]{$\mathcal{R}^{g}_{0^+}$}
$$
\begin{array}{lll}
{\cal R}^g_{0^+,t}(\xi_T)&=\lim_{\gamma\downarrow 0}{\cal
R}^g_{\gamma,t}(\xi_T)=Y^{0}_t(\xi_T)\\
&={\rm ess\sup}_{\mu \in {\cal A}_k}\big\{{\mathbb E}_{{\mathbb
Q}^{\mu}}\big[-\xi_T\big|{\cal F}_t\big]\big|\>G(u,\mu_u)<\infty
\>\forall u\geq t,\>du \>-a.s..\big\}\\
&={\rm ess\sup}_{\mu \in {\cal A}_k \cap {\rm Dom}(G) }{\mathbb
E}_{{\mathbb Q}^{\mu}}\big[-\xi_T\big|{\cal F}_t\big].
\end{array}
$$
$\bullet$ When the coefficient $g$ has a quadratic growth $(H3)$,
the recession function may be infinite on a set with positive
measure and the BSDE($g_{0^+},\xi_T$) is not well-defined.
However, we can still take the limit in the dual representation of
${\cal R}^g_{\gamma,t}$, obtain the same characterization of
${\cal R}^g_{0^+,t}$, and consider ${\cal R}^g_{0^+}$ as a
generalized solution of BSDE whose the coefficient $g_{0^+}$ may
be take infinite values. In particular if, as in the entropic
case, $g_{0^+}=l_{\{0\}}$, $G$ is finite everywhere and any
equivalent probability measure associated with BMO coefficient,
said to be in ${\cal Q}(\rm BMO)$, is admissible. Then,
$${\cal R}^{l_{\{0\}}}_{0^+,t}(\xi_T)={\rm
ess\sup}_{{\mathbb Q}\in{\cal Q}(\rm BMO)}{\mathbb E}_{{\mathbb
Q}}\big[-\xi_T\big|{\cal
F}_t\big]=||-\xi_T||^{\infty}_t=e_{0^+,t}(\xi_T).$$ \vspace{-5mm}
\paragraph{\sc \index[sub]{recession super price system}\index[sub]{super price system}Super Price System }
Note that the conservative risk-measure ${\cal
R}^g_{0^+,t}(\xi_T)= {\rm ess\sup}_{\mu \in {\cal A}\cap{\rm
Dom}(G)}{\mathbb E}_{{\mathbb Q}^{\mu}}\big[-\xi_T\big|{\cal
F}_t\big]$ is the equivalent of the super-pricing rule of $-\xi_T$
(this notion was first introduced by \index[aut]{El Karoui} El
Karoui and \index[aut]{Quenez} Quenez \cite{ElKarouiQuenez95}
under the name "upper hedging price"). When the
$\lambda_t$-translated of ${\rm Dom}(G)_t$ is a vector space, the
recession function $g_{0^+}(t,z)$ is the indicator function of the
orthogonal vector space ${\rm Dom}(G)^\top_t$ plus a linear
function $\langle z,-\lambda_t\rangle$. Then, ${\cal
R}^g_{0^+,t}(-\xi_T)$
is exactly the upper-hedging price associated with hedging portfolios constrained to live in ${\rm Dom}(G)^\top$.\\
The conservative measure is the smallest of coherent risk measure
such that ${\cal R}^g_t(-\xi_T)-{\cal R}^g_t(-\eta_T)\leq {\cal
R}^{\rm coh}_t(-\xi_T+\eta_T)$ for any $(\xi_T, \eta_T)$ in the
appropriate space.
\paragraph{\sc \index[sub]{volume perspective risk measure}Volume Perspective Risk Measure}
It is also possible to associate a coherent risk measure
$\mathcal{R}^{\tilde g}$ with any convex risk measure
$\mathcal{R}^{g}$, using the \index[sub]{perspective
function}perspective function $\tilde g$ of the coefficient $g$,
which is assumed to be normalized for the sake of simplicity
($g(t,0)=0$). The perspective function $\tilde{g}$ is defined as:
$$\tilde{g}(t,\gamma,z)=
\left\{\begin{array}{ccc} &\gamma g(t,\frac{z}{\gamma}) \quad &{\rm if}\> \gamma>0\\
&\lim_{\gamma \rightarrow 0} \gamma
g(t,\frac{z}{\gamma})=g_{0^+}(t,z) \quad &{\rm if}\> \gamma =0
\end{array}
\right.
$$
More details about $\tilde{g}$ can be found in the Appendix
\ref{Perspective function}. As a direct consequence, the $\tilde
g$-conditional risk measure $\mathcal{R}^{\tilde g}$ is a coherent
risk measure.
\section{Inf-Convolution of $g$-Conditional Risk Measures \label{subsection Inf-convolution of g-conditional risk measures}}
In this section, we come back to inf-convolution of risk measures,
when they are g-conditional risk measures. This study is based
upon the
inf-convolution of their respective coefficients.\\
More precisely, we will study for any $t$ the inf-convolution of
the g-conditional risk measures ${\cal R} _{t}^{A}$ and $%
{\cal R} _{t}^{B}$ defined as
\begin{equation}
\big( {\cal R} ^{A}\square {\cal R} ^{B}\big) _{t}\big( \xi_T\big)
={{\rm ess}\inf}_{F_T}\big\{ {\cal R} _{t}^{A}\big(\xi_T-F_T\big)
+{\cal R} _{t}^{B}\big( F_T\big) \big\}   \label{BSDE: pb
inf-convolution}
\end{equation}
where both $\xi_T$ and $F_T$ are taken in the appropriate space
and show that this new dynamic risk measure is under mild
assumptions the (maximal solution) ${\cal R}^{A,B}$ of the BSDE
$(g^{A}\square g^{B},-\xi_T)$ where $(g^{A}\square g^{B})(.,t,z)=
{{\rm ess}\inf}_{z}(g(.,t,x-z)+g(.,t,z)).$ Then, the next step is
to characterize the optimal transfer of risk between both agents A
and B, agent A being exposed to $\xi_T$ at time $T$. Some key
results on the inf-convolution of convex functions are recalled in
the Appendix \ref{Subsection inf-convolution convex functions},
the main argument being summarized in the proposition below:
\begin{proposition}\label{theorem inf-convolution atteint}
Let $g^A$ and $g^B$ be two convex functions of $z$. Under the
following condition
\begin{equation*}
g^A_{0^+}(t,z)+g^B_{0^+}(t,-z)>0,\quad \forall z\neq 0
\end{equation*}
then $g^A \square g^B$ is exact for any z as the infimum is
attained by some $x^*$: $$ g^A\square
g^B(z)=\inf_x\{g^A(z-x)+g^B(x)\}=g^A(z-x^*)+g^B(x^*).$$
\end{proposition}
\subsection{Inf-convolution and Optima}
We now focus on our main problem of inf-convolution of
$g$-conditional risk measures as expressed in Equation (\ref{BSDE:
pb inf-convolution}). The following theorem gives us an explicit
characterization of an optimum for the inf-convolution problem
provided such an optimum exists:
\begin{theorem}
\label{theorem BSDE F* explicite}Let $g^A$ and $g^B$ be two convex
coefficients depending only on $z$ and satisfying the condition of
Proposition \ref{theorem inf-convolution atteint}. For a given
$\xi_T$ in the appropriate space (either $\mathbb{L}^2$ or
$\mathbb{L}_{\infty}$), let $({\cal R}_{t}^{A,B}(\xi_T),Z_{t})$ be
the maximal solution of the BSDE $(g^{A}\square g^{B},-\xi_T)$ and
$\widehat{Z}^{B}_{t}$ be a measurable process such that $\quad
\widehat{Z}^{B}_{t}=\arg \min_{x} \Big\{ g^{A}\big( t,Z_{t}-x\big)
+g^{B}\big( t,x\big) \Big\} \quad
dt \times d\mathbb{P}\>-a.s.$.\\
Then, the following results
hold:\vspace{1mm}\\
$(1)$ For any $t \in [0,T]$ and for any $F_T$ such that both
${\cal R}_t^{A}(\xi_T-F_T)$ and ${\cal R}_t^{B}(F_T)$ are well
defined:
$${\cal R}_{t}^{A,B}\big(\xi_T\big)
\leq {\cal R}_t^{A}(\xi_T-F_T)+ {\cal R}_t^{B}(F_T) \quad
\mathbb{P}\ -a.s.$$ $(2)$ If the process $\widehat{Z}^{B}$ is
admissible, then for any $t \in [0,T]$
$${\cal R}
_{t}^{A,B}\big(\xi_T\big) = ({\cal R} ^{A}\square {\cal R}
^{B}\big) _{t}\big( \xi_T\big) \quad \mathbb{P}-a.s.$$ and the
structure $F^{*}_T$ defined by the forward equation
$$F^{*}_T=\int\limits_{0}^{T}g^{B}\big( t,\widehat{Z}_{t}^{B}\big)
dt-\int\limits_{0}^{T} \widehat{Z}_{t}^{B}dW_{t}$$ is an optimal
solution for the inf-convolution problem:
$$({\cal R} ^{A}\square {\cal R}
^{B}\big) _{t}\big( \xi_T\big)={\cal R}_t^{A}(\xi_T-F^*_T)+ {\cal
R}_t^{B}(F^*_T).$$
\end{theorem}
{\bf Proof}: $(1)$ First, note that the existence of such a
measurable process $\widehat{Z}^{B}_{t}$ is guaranteed by Theorem
\ref{theorem inf-convolution
atteint}.\\
In the following, we consider any $F_T$ such that both ${\cal
R}_t^{A}(\xi_T-F_T)$ and
${\cal R}_t^{B}(F_T)$ are well defined. \\
Let us now focus on ${\cal R}_t^{A}(\xi_T-F_T)+{\cal
R}_t^{B}(F_T)$. It satisfies
\begin{equation*}
\begin{array}{lll}
-d\big( \mathcal{R}_t^A(\xi_T-F_T)+\mathcal{R}_t^B(F_T)\big)&=&
\big(g^{A}(t,Z_{t}^A)+g^{B}(t,Z_t^B)\big)dt-\big(Z_t^A+Z_t^B\big)dW_t\\
 &=&\left(g^{A}(t,Z_{t}-Z_{t}^{B})
+g^{B}(t,Z_{t}^{B})\right)dt-Z_{t}dW_t,
\end{array}
\end{equation*}
and at time $T$, $\mathcal{R}_T^A(\xi_T-F_T)+\mathcal{R}_T^B(F_T)=-\xi_T$.\\
Therefore, $({\cal R}_t^{A}(\xi_T-F_T)+{\cal R}_t^{B}(F_T),Z_t)$
is solution of the BSDE with terminal condition $-\xi_T$, which is
also the terminal condition of the BSDE $(g^{A}\square
g^{B},-\xi_T)$, and a coefficient $g$ written in terms of the
solution $Z_t^B$ of the BSDE $(g^B,F_T)$ as:
$g(t,z)=g^{A}(t,z-Z_t^B)+g^{B}(t,Z_{t}^{B})$. Using the definition
of the inf-convolution, this coefficient is then always greater
than $g^{A}\square g^{B}$. Thus, we can compare ${\cal
R}_t^{A}(\xi_T-F_T)+{\cal R}_t^{B}(F_T)$ with the solution of the
BSDEs $(g^{A}\square g^{B},-\xi_T)$ using the comparison Theorem
(\ref{comparison Theorem})
and obtain the desired inequality.\vspace{1mm}\\
$(2)$ Let now assume that the process $\widehat{Z}^{B}_{t}$ is
admissible, using different notions of admissibility when either
(H1) or (H3) (square
integrability or BMO).\\
Thanks to Theorem \ref{theorem inf-convolution atteint}, we can show that both dynamic risk measures coincide.\\
We now introduce the structure $F^{*}_T$ defined by the forward
equation $ F^{*}_t=\int\limits_{0}^{t}g^{B}\big(
s,\widehat{Z}_{s}^{B}\big) ds-\int\limits_{0}^{t}
\widehat{Z}_{s}^{B}dW_{s}
$.\\
Note first that thanks to the admissibility of the process
$\widehat{Z}^{B}_{t}$, such a structure is well-defined and
belongs to the appropriate space
(either $L^2(\mathcal{F}_T)$ or $L^{\infty}(\mathcal{F}_T)$).\\
Let us also observe that $-F^{*}_t$ is also solution of the BSDE
$(g^B,-F^{*}_T)$ since  $
-F^{*}_t=-F^{*}_T+\int\limits_{t}^{T}g^{B}\big( u,
\widehat{Z}_{u}^{B}\big) dt-\int\limits_{t}^{T} \widehat{Z}
_{u}^{B}dW_{u} $.
 By uniqueness, this process is ${\cal
 R}_{t}^{B}(F^{*}_T)$.\vspace{2mm}\\
Since ${\cal R}_t^{A}(\xi_T-F^{*}_T)+{\cal R}_t^{B}(F^{*}_T)$ is
solution of the BSDE with coefficient written as
$g^{A}(t,Z_t-\widehat{Z}_t^{B}) +g^{B}(t,\widehat{Z}_t^{B})$ and
terminal condition $-\xi_T$ and given that $\big( g^{A}\square
g^{B}\big) \big( t,Z_{t}\big) =g^{A}\big(
t,Z_{t}-\widehat{Z}_{t}^{B}\big) +g^{B}\big( t,\widehat{Z}
_{t}^{B}\big)$, by uniqueness, we also have $\forall t\geq 0$,
${\cal R} _{t}^{A,B}\big( \xi_T\big) =\big( {\cal R} ^{A}\square
{\cal R} ^{B}\big) _{t}\big( \xi_T\big) \quad \mathbb{P}\ a.s
$.\\
The proof also gives the optimality for the Problem (\ref{BSDE: pb
inf-convolution}) of the structure
$F^{*}_T=\int\limits_{0}^{T}g^{B}\big( t,\widehat{Z}_{t}^{B}\big)
dt-\int\limits_{0}^{T} \widehat{Z}_{t}^{B}dW_{t} $. $\qquad
\square$
\begin{remark}[On uniqueness on the optimum]
Note that the optimal structure $F^*_T$ is determined to within a
constant because of the translation invariance property (P3-)
satisfied by both risk measures ${\cal R}_{t}^{A}$ and ${\cal
R}_{t}^{B}$ since:
$$
\begin{array}{ccc}
{{\rm ess}\inf}_{F_T}\big\{ {\cal R}
_{t}^{A}\big(\xi_T-(F_T+m)\big) +{\cal R} _{t}^{B}\big( F_T+m\big)
\big\} &=&{{\rm ess}\inf}_{F_T}\big\{ {\cal R}
_{t}^{A}\big(\xi_T-F_T\big)+m +{\cal R} _{t}^{B}\big( F_T\big)-m
\big\}\\
&=&\big( {\cal R} ^{A}\square {\cal R} ^{B}\big) _{t}\big(
\xi_T\big).
\end{array}
$$
Note also that $F^*_T$ is optimal for all the optimal structure
problems for all stopping times $S$ such that $0 \leq S \leq T
\quad a.s.$.
\end{remark}
The following Theorem gives some sufficient conditions ensuring
the admissibility of the process $\widehat{Z}^B_t$:
\begin{theorem}{\rm [Exact Inf-convolution]}
Let $g^B$ be a strongly convex coefficient. For any convex
function $g^A$, the inf-convolution $g^A \square g^B$ is convex
with quadratic growth (H3), so in particular, if $g^A$ satisfies
(H3).\\
In this case, the process $\widehat{Z}^B_t$, defined in Theorem
\ref{theorem BSDE F* explicite}, is in ${\rm BMO(\mathbb{P})}$.
\end{theorem}
Note that in this case, the optimal structure $F^*_T$, defined in
Theorem \ref{theorem BSDE F* explicite}, is quasi-bounded as it
belongs to the BMO-closure of $\mathbb{L}_{\infty}$ as defined by
Kazamaki \cite{Kazamaki} (chapter 3).\vspace{2mm}\\
{\bf Proof}:
From the duality Theorem \ref{theorem dual representation}, the
optimal control $\mu^*$ of $G^{A,B}$, the polar function of $g^A
\square g^B$, is in ${\rm BMO(\mathbb{P})}$. From the
inf-convolution, we deduce that this is also the optimal control
for $G^A$ and $G^B$ in the following sense:
$$\begin{array}{ll}
g^A(t,Z_t-\widehat{Z}^B_t)=\langle \mu^*_t,-(Z_t-\widehat{Z}^B_t)
\rangle -G^A(t,\mu^*_t),\\
g^B(t,\widehat{Z}^B_t)=\langle \mu^*_t,-\widehat{Z}^B_t \rangle
-G^B(t,\mu^*_t).
\end{array}$$
Therefore, both $Z_t-\widehat{Z}^B_t$ and $\widehat{Z}^B_t$ are in
${\rm BMO(\mathbb{P})}$ (from Proposition \ref{propo:MZMartBMO})
and the process $\widehat{Z}^B_t$ is admissible. $\qquad
\square$\vspace{2mm}\\
{\bf Comments:}\\
$(i)$ Just as in the static framework, we obtain the same result
when considering $g$-conditional $\gamma$-tolerant risk measures.
The Borch theorem is therefore extremely robust since the quota
sharing of the
initial exposure remains an optimal way of transferring the risk between different agents.\\
$(ii)$ Under some particular assumptions, the underlying logic of
the transaction is non-speculative since there is no interest for
the first agent to transfer some risk or equivalently to issue a
structure if she is not initially exposed. This result is
completely consistent with the result we have already obtained in
the static framework.
\subsection{Hedging Problem}
As in subsection \ref{static individual hedging}, we consider the
hedging problem of a single agent. She wants to hedge her terminal
wealth $X_T$ by optimally investing on financial market and
assesses her risk using a general g-conditional risk measure
${\cal R}^g$.
\subsubsection{Framework}
We consider the same framework as that introduced in Subsection
\ref{Dynamic hedging - static} when looking at the question of
dynamic hedging in the static part. More precisely, we assume that
$d$ basic securities are traded on the market. Their forward
(non-negative) vector price process $S$ follows an It\^{o}
semi-martingale with a uniformly bounded drift coefficient and an
invertible and bounded volatility matrix $\sigma_t$. Under
$\mathbb{P}$,
\begin{equation}\label{eq:dynamics}
\frac{dS_t}{S_t}=\sigma_t(dW_t+\lambda_t dt)\quad; \quad S_0
\quad{\rm given}.
\end{equation}
To avoid arbitrage, we assume {\bf (AAO)}: there exists a
probability measure $\mathbb{Q}$, equivalent to $\mathbb{P}$, such
that $S$ is a $\mathbb{Q}$-local martingale. From the completeness
of this basic arbitrage-free market, we deduce the uniqueness of
$\mathbb{Q}$, which is usually called the
risk-neutral probability measure.\\
The agent can invest in dynamic strategies $\theta$, i.e.
$d$-predictable
processes and $(G_t(\theta)=(\theta.S)_t)$ denotes the associated gain process.\\
We assume that not all strategies are admissible and that, for
instance, the agent has some restriction imposed on the
transaction size. These constraints create some market
incompleteness in the framework we consider.
$\Theta_T^S=\{G_T(\theta)\>|\>\theta.S \> \>{\rm is \> bounded
\>by\> below}\>,\>\theta \in \mathcal{K}\}$ is the set of
admissible hedging gain processes. $\mathcal{K}$ is a convex
subset of ${\rm BMO}(\mathbb{P})$ such that any admissible
strategies $\theta$ is in $\mathcal{K}$ (equivalently, $\forall
\>t,\>\> \theta_t \in \mathcal{K}_t$).\\
\subsubsection{Hedging Problem} At time $0$, the hedging problem of the agent can
be expressed as the determination of an optimal admissible
strategy $\theta$ as to minimize the initial g-conditional risk
measure of her terminal wealth
\begin{equation}
\inf_{\theta \in \mathcal{K}}{\cal R}_0\big(X_T-G_T(\theta) \big).
\label{hedging pb agent A - dynamic 0}
\end{equation}
The value functional of this program is the dynamic market
modified risk measure of agent A, denoted by $\mathcal{R}^{m}$.
Using the previous results, we can obtain the following
proposition:
\begin{proposition}
i) Let $l_{\sigma_t^*}(\mathcal{K}_t)=l_{\widehat{\mathcal{K}_t}}$
be the indicator function of the convex set
$\widehat{\mathcal{K}_t}=\sigma_t^*\mathcal{K}_t$. Provided that
the inf-convolution $g \square l_{\sigma_t^*}(\mathcal{K}_t)
(Z_t)$ is well-defined, the residual risk measure
$\mathcal{R}^{m}$ is given as the maximal solution of the
following BSDE:
\begin{equation*}
-d\mathcal{R}_t^{m}(X)=g^{m}(t,Z_t)dt-Z_T dW_t \quad ; \quad
\mathcal{R}^{m}_T(X)=-X_T
\end{equation*}
where $g^{m}$ is the restriction of the coefficient $g$ to the
admissible set: $g^{m}(t,Z_t)=g \square
l_{\sigma_t^*}(\mathcal{K}_t) (Z_t)$.\\
ii) If $g$ is strongly convex, then this hedging problem has a
solution.\\
In particular, in the entropic case,
$g^{m}(t,z)=\frac{1}{2\gamma}d_{\frac{1}{\gamma}}(z,\widehat{K_t})^2$
where $\gamma$ is the risk tolerance coefficient and
$d_{\frac{1}{\gamma}}(z,\mathcal{K})$ is the distance function to
$\mathcal{K}$. The optimal investment strategy $\theta^{\star}$ is
the projection on $\mathcal{K}$
of $Z_t$, solution of the BSDE $(g^{m},-X_T)$.\\
The terminal value $G_T(\theta^{\star})$ of the associated
portfolio is given by:
$$G_T(\theta^{\star})=x + \int_0^T (\theta_t^{\star})^*\sigma_t \lambda_t dt +
\int_0^T (\theta_t^{\star})^*\sigma_t  dW_t.$$
\end{proposition}
\subsubsection{Comments}
{\bf Generalized BSDEs:} In the static framework, we expressed the
hedging problem as an inf-convolution between the seminal risk
measure of the agent and the risk measure $\nu^{\mathcal{H}}$
generated by $\mathcal{H}$, the convex set of constrained terminal
gains, or more generally the inf-convolution between the seminal
risk measure of the agent
and the convex indicator of $\mathcal{H}$ (Proposition \ref{proposition regularisation par worst}).\\
From a dynamic point of view, the set $\mathcal{H}$ can be seen as
the set of all dynamic terminal values of portfolios with some
constraint on the strategies. Everything can be formulated in the
same way. Note that the natural candidate for
$\mathcal{R}^{\mathcal{H}}$ would be the inf-convolution between
the dynamic worst case risk measure and the convex indicator of
$\mathcal{H}$: $l^{\mathcal{H}} \square \lim_{\gamma \rightarrow
\infty}(\frac{1}{2\gamma}|z|^2)$. This infimum is always strictly
positive. Moreover, it is an increasing process at the limit. To
model this "limit BSDE", an increasing process has to be
introduced (for more details, please refer to \index[aut]{El
Karoui}El Karoui and \index[aut]{Quenez}Quenez
\cite{ElKarouiQuenez96} and \index[aut]{Cvitanic}Cvitanic and
\index[aut]{Karatzas}Karatzas \cite{CvitanicKaratzas99}). As a
consequence, the dynamic version of the risk measure generated by
$\mathcal{H}$ cannot be seen exactly as the solution of a standard
BSDE, as previously defined, in the sense that the
coefficient can take infinite values. \\
This is however not such a problem here as we really focus on the
inf-convolution. Therefore, we can simply consider the restriction
of the seminal risk measure to a particular set. The powerful
regularization impact of the inf-convolution is again visible
here.\vspace{2mm}\\
{\bf Hedging problem at any time $t$:} Solving the hedging problem
at time $0$ leads to the characterization of a particular
probability measure, which can be called {\it calibration
probability measure} as the prices of any hedging instruments made
with respect to this measure coincide
with the observed market prices on which all agents agree.\\
Solving the hedging problem at any time $t$ is equivalent to
solving the same problem at time $0$ as soon as the prices of
these hedging instruments at this time $t$ are given as the
expected value of their discounted future cash flows under the
optimal calibration probability measure determined at time $0$.
This optimal probability measure is very robust as it
remains the pricing measure for hedging instruments between $0$ and $T$.\\
Therefore, we can introduce the same problem at any time $t$: $$
{\rm ess}\inf_{\theta \in \mathcal{K}}{\cal
R}_t\big(X_T-G_T(\theta) \big)= \mathcal{R}_t^{m}(X_T).
$$
BSDEs time-consistency and uniqueness are key arguments to show
that if $\theta$ is optimal for the problem at time $0$, then
$\theta$ is optimal for the optimization program at any time
$t$.\vspace{2mm}\\
{\bf Dynamic Entropic Framework} The entropic hedging problem,
lying at the core of this book, has been intensively studied in
the literature. But only a few papers are using a BSDEs framework.
After the seminal paper by \index[aut]{El Karoui}El Karoui and
\index[aut]{Rouge}Rouge \cite{ElKarouiRouge}, different authors
have used BSDEs to solve this problem under various assumptions
(see in particular \index[aut]{Sekine}Sekine \cite{Sekine02},
\index[aut]{Mania}\index[aut]{Santacroce}\index[aut]{Tevzadze}Mania
et al. \cite{ManiaSantacroceTevradze} and more recently
\index[aut]{Hu} Hu, \index[aut]{Imkeller}Imkeller and
\index[aut]{M\"{u}ller}M\"{u}ller
\cite{HuImkellerMuller} or \index[aut]{Mania}Mania and \index[aut]{Schweizer}Schweizer \cite{ManiaSchweizer}).\vspace{2mm}\\
Another approach, different from what we have mentioned above, has
been used to solve the hedging problem involves the dual
representation for the dynamic entropic risk measure as given by
Theorem \ref{theorem dual representation}:
$$e_{\gamma ,t}(\Psi)=\sup_{\mu \in {\cal
A}^q}{ \mathbb E}_{{\mathbb Q}^{\mu}}\Big[-\Psi-\gamma \int_t^T
\frac{|\mu_s|^2}{2}ds | \mathcal{F}_t\Big].
$$
Therefore, the hedging problem at any time $t$ can be rewritten
as:
\begin{equation*}
{\rm ess}\inf_{\theta \in \mathcal{K}}{\rm ess\sup}_{\mu \in {\cal
A}^q} \Big\{{\mathbb E}_{{\mathbb Q}^{\mu}}[-X_T+G_T(\theta) |
\mathcal{F}_t] -\gamma h(\mathbb{Q}^{\mu}|\mathbb{P})\Big\}
\end{equation*}
and it may be solved by using dynamic programming arguments.
\section{Appendix: Some Results in Convex Analysis \label{Appendix convex analysis}}
We now present some key results in convex analysis
that will be useful to obtain the dual representation of
$g$-conditional risk measures. More details or proofs may be found
in \index[aut]{Aubin} Aubin \cite{Aubin},
\index[aut]{Hiriart-Urruty} Hiriart-Urruty and
\index[aut]{Lemar\'{e}chal} Lemar\'{e}chal
\cite{Hiriart-Urruty-Lemarchal} or \index[aut]{Rockafellar} Rockafellar \cite{Rockafellar}.\\
All the notations and definitions we introduce are consistent with the notations of risk measures. They may differ
from the standard framework of convex analysis (especially regaring the sign).\vspace{2mm}\\
Even if the coefficient of the BSDE is finite, we are also
interested in convex functions taking infinite values. The main
motivation is the definition of its convex polar function $G$. In
that follows, as in \cite{Hiriart-Urruty-Lemarchal}, we always
assume that the considered functions are not identically $+\infty$
and are bounded from below by a affine function (note that this
assumption is rather general and does not necessarily require that
the functions are convex). The domain of a function $g$ is defined
as the nonempty set \index[not]{${\rm Dom}(g)$} ${\rm
Dom}(g)=\{z\>:\>g(z)<+\infty\}$. The \index[sub]{epigraph}
epigraph of convex function is the subset of ${\mathbb R}^n\times
{\mathbb R}$ as: ${\rm epi\, g}=\{(x,\lambda)\>|\>g(x)\leq
\lambda\}$. When the convex functions are lower semicontinuous
\index[not]{${\rm lsc}$}(lsc), ${\rm epi g}$ is closed, and they
are said to be closed.
\subsection{Duality}
\subsubsection{Legendre-Fenchel Transformation} Let $g$ be a convex
function. The polar function $G$ is defined on ${\mathbb R}^n$ by
\begin{equation}
\label{eq:polar}
G(\mu)=\sup_z(\langle\mu,z\rangle-g(z))=\sup_{z\in {\rm Dom}(g)
}(\langle\mu,-z\rangle-g(z)).
\end{equation}
The function $G$ is a closed convex function, which can take
infinite values. The conjugacy operation induces a symmetric
one-to-one correspondence in the class of all closed convex
functions on $\mathbb{R}^n$ \vspace{-3mm}:
$$g(z)=\sup_\mu(\langle\mu,-z\rangle-G(\mu)),\quad G(\mu)=\sup_z(\langle\mu,-z\rangle-g(z)).$$
\paragraph{\sc Convex set and duality}
Given a nonempty subset $S\subset {\mathbb R}^n$, the {\it
indicator function} (in the convex analysis terminology) of S,
$l_S\>:\>{\mathbb R}^n\rightarrow \mathbb{R}^+\cup\{+\infty\}$, is
defined by:
$$l_S(z)=0\quad{\rm if}\quad z\in S\qquad {\rm and} \quad +\infty \quad{\rm if\>not}.$$
$l_S$ is convex (closed), iff S is convex
(closed) since ${\rm epi}\>l_S\>=S\times \mathbb{R}^+$.\\
The polar function of $l_S$ is the \index[sub]{support function}
{\it support function} of $-S$:
$$\sigma_S(z):=\sup_{s\in S}\langle s,-z\rangle=\sup_s\{\langle
 s,-z\rangle-l_S(s)\}.$$
The support function is closed, convex, homgeneous function:
$\sigma_S(\lambda z)=\lambda \sigma_S(z)$ for all $\lambda>0$. Its
epigraph and its domain are convex cones.
\subsubsection{Subdifferential and Optimization} The
sub-differential of the convex function $g$ in $z$, whose the
elements are called {\it subgradient} of $g$ at $z$, is the set
\index[not]{$\partial g$}$\partial g(z)$ defined as:
\begin{equation}\label{eq: sous-gradient}
\partial g(z)=\{\> \mu\>|\>g(x)\geq g(z)-\langle \mu,x-z\rangle,\quad \forall x\}=
\{\>\mu\>|g(z)-\langle \mu,-z\rangle \geq G(\mu)\}.
\end{equation}
If $z \notin Dom(g),\>\> \partial g(z)=\emptyset$. But if $z$ is
in the interior of ${\rm Dom}(g)$, the subgradient $\partial g(z)$
is non-empty (see Section E in \cite{Hiriart-Urruty-Lemarchal} or
Chapter 23 in \cite{Rockafellar}); in fact, it is enough that $z$
belongs to the relative interior of ${\rm Dom}(g)$, where ${\rm
ridom}(g)$ is defined in Section A in
\cite{Hiriart-Urruty-Lemarchal} and in Chapter 6 in
\cite{Rockafellar}. In particular, {\em if $g$ is finite, then
$\partial g(z)$ is nonempty for any $z$}. When $\partial g(z)$ is
reduced to a single point, the function is said to be
differentiable in $z$. Note that when the function $g$ is the
indicator function of the convex set $C$, the sub-differential of
$g$ in $z\in C$
is the positive normal cone $N^+_C(z)$ to $C$ at $z$, $N^+_C(z)=\{s\in {\mathbb R}^n\>|\>\forall y\in C \quad -\langle s,y-z\rangle\}\leq 0\}$.\vspace{2mm}\\
Subgradients are solutions of minimization programs as $ \inf_z
\left(g(z)-\langle \mu,-z \rangle \right) \>\> (=-G(\mu)), $ or its
dual program, $\inf_{\mu}\left(G(\mu)-\langle \mu,-z \rangle
\right) \>\> (=-g(z))$. The precise result is the following (see
Section E in \cite{Hiriart-Urruty-Lemarchal}): Let $g$ be a closed
convex function and $G$ its polar function.
\begin{itemize}
\item [$\bullet$] $\widehat{\mu} \in \partial g(\widehat{z})$
$\Longleftrightarrow$ $\widehat{\mu}$ is the optimal for the
following minimization program, that is
$-g(z)=\inf_{\mu}\left(G(\mu)-\langle \mu,-z \rangle \right)=
G(\widehat \mu)-\langle \widehat\mu,-z \rangle $. \item [$\bullet$]
$\widehat{z} \in \partial G(\widehat{\mu})$ $\Longleftrightarrow$
$\widehat{z}$ is optimal for the following minimization program,
that is in  $ -G(\mu)=\inf_z \left(g(z)-\langle \mu,-z \rangle
\right)=g(\widehat z)-\langle \mu,-\widehat{z} \rangle.$
\end{itemize}
In the following, when working with BSDEs, we will denote by
$z\mu$ the scalar product between the line vector $z$ and the
column vector $\mu$.
\subsection{Recession function}
\subsubsection{Recession Function \label{par:recession function}}
The recession function associated with a closed convex function
$g$ is the homogeneous convex function defined for $z\in {\rm
Dom}(g)$ by $g_{0^+}(z) :=\lim_{\gamma \downarrow 0}\gamma
g\big(\frac{1}{\gamma}z\big)= \lim_{\gamma \downarrow 0}\gamma
\big(g(y+\frac{1}{\gamma}z)-g(y)\big)$. This function $g_{0^+}$ is
the smallest homogeneous function $h$ such that for any $z,y\in
{\rm Dom}(g)$, $g(z)-g(y)\leq h(z-y)$. When $g(z)\leq c+k|z|$,
$g_{0^+}(z)\leq k|z|$ is a finite convex function,
and the function $g$ is Lipschitz-continuous function with Lipschitz coefficient $k$ since $g(z)-g(y)\leq k|z-y|$.\\
{\em This property explains why any convex coefficient of BSDE satisfying the assumption $(H2)$ in fact satisfies $(H1)$.}\\
Let $G$ be the polar function of $g$. Using obvious notations, for
any $\mu \in {\rm Dom}(g),$ ${\rm
\>polar\>}g_{0^+}(\mu)=\lim_{\gamma\downarrow 0}(\gamma G(\mu))
=0$. So, ${\rm \>polar\>}g_{0^+}=l_{{\rm dom\> G}}$. By the
conjugacy relationship applied to closed functions, $g_{0^+}$ is
the support function of ${\rm Dom}(G)$; so $g_{0^+}$ is finite
everywhere iff ${\rm Dom}(G)$ is bounded, or iff $g$ is uniformly
Lipschitz, or finally
iff $g$ has linear growth. \\
The recession function of the quadratic function $q_k(z)=c+k|z|^2$
is infinite except in $z=0$, and its polar function is the null
function. More generally, convex functions such that
$g_{0^+}=l_{\{0\}}$ admit finite polar function $G$ and this
condition is sufficient.
\subsubsection{Perspective Function \label{Perspective function}}
Let us consider a closed convex
function $g$ such that $g(0)=0$. The perspective function
associated with $g$ is the function $\tilde g$ defined on
${\mathbb R}^+\times{\mathbb R}^n $ as: 
 $$\tilde{g}(\gamma,z)=
\left\{\begin{array}{ccc} &\gamma g(\frac{z}{\gamma}) \quad &{\rm if}\> \gamma>0\\
&\lim_{\gamma \rightarrow 0} \gamma g(\frac{z}{\gamma})=g_{0^+}(z)
\quad &{\rm if}\> \gamma =0
\end{array}
\right.
$$
Note first that the perspective function of $g$ corresponds to the
$\gamma $-dilated of $g$, seen as a function of both variables $z$
and $\gamma$, when $\gamma>0$. It is prorogated for $\gamma=0$ by
the recession function $g_{0^+}$. Note that the risk tolerance
coefficient is considered as a risk factor itself. $\tilde{g}$ is
a positive homogeneous convex function (for more details, please
refer to Part B \cite{Hiriart-Urruty-Lemarchal}).\\
The dual function of $\tilde{g}$, defined on $\mathbb{R} \times
\mathbb{R}^{n}$, is given by:
$$\tilde{G}(\theta,\mu)= 0 \quad {\rm if}\> G(\mu)\leq -\theta
\>, \quad {\rm and} \quad +\infty \quad {\rm otherwise}.$$ If
$g(0)<\infty$, note that $G(\mu)$ is bounded.
\subsection{Infimal Convolution of Convex Functions and Minimization Programs \label{Subsection inf-convolution convex functions}}
Addition and inf-convolution of closed convex functions are two
dual operations with respect to the conjugacy relation.\\
Let $g^A$ and $g^B$ be two closed convex functions from $\mathbb
R^n\cup \{+\infty\}$. By definition, the infimal convolution of
$g^A$ and $g^B$ is the function $g^{A}\square g^{B}$ defined as:
\begin{equation}
\label{eq:infconvolution} \big(g^{A}\square g^{B}\big)(z)=
\inf_{y^A+y^B=z}(g^A(y^A)+g^B(y^B))=\inf_{y}(g^A(z-y)+g^B(y)).
\end{equation}
If $g^{A}\square g^{B} \not\equiv \infty$, then its polar
function, denoted by $G^{AB}$, is simply the sum of the polar
functions of $g^A$ and $g^B$:
$$G^{AB}(\mu)=G^A(\mu)+G^B(\mu)$$
\subsubsection{Inf-Convolution as a Proper Convex Function}
The function $g^{A}\square g^{B}$ may take the value $-\infty$,
which is contrary to the assumption made in Subsection
\ref{Subsection results convex analysis}. To avoid this
difficulty, we assume that both functions $g^{A}$ and $g^{B}$ have
a common affine minorant $\langle s,.\rangle-b$. This assumption
may be expressed in terms of their recession functions, both of
them being also bounded from below by $\langle s,.\rangle$.
Therefore, $g^{A}_{0^+}(z)+g^{B}_{0^+}(-z)\geq 0$ for any $z$ and
consequently $\big(g^{A}_{0^+}\square g^{B}_{0^+}\big)(0) \geq 0$.
Note that this condition can also be expressed in terms of the
polar functions of $g^A$ and $g^B$ as ${\rm dom}(G^A)\cap {\rm
dom}(G^B)\neq \emptyset$.
\subsubsection{Existence of \index[sub]{exact inf-convolution}Exact Inf-Convolution}
We are interested in the existence of a solution to the
inf-convolution problem (\ref{eq:infconvolution}). When a solution
exists, the infimal convolution is said to be {\it exact}.\\
The previous conditions are almost sufficient, as proved in
\index[aut]{Rockafellar} Rockafellar \cite{Rockafellar} since, if
we assume
\begin{equation}\label{condition inf-convolution atteint}
g^A_{0^+}(z)+g^B_{0^+}(-z)>0,\quad \forall z\neq 0
\end{equation}
then $g^A \square g^B$ is a closed convex function, and for any z,
the infimum is attained by some $x^*$: $$ g^A\square
g^B(z)=\inf_x\{g^A(z-x)+g^B(x)\}=g^A(z-x^*)+g^B(x^*).$$ The
condition (\ref{condition inf-convolution atteint}) is satisfied
if ${\rm intdom}(G^A)\cap {\rm intdom}(G^B)\neq \emptyset$ (in
fact, the true interior corresponds to the relative interior
defined in Section A by \index[aut]{Hiriart-Urruty}Hiriart-Urruty
and \index[aut]{Lemar\'{e}chal} Lemar\'{e}chal
\cite{Hiriart-Urruty-Lemarchal}).\vspace{2mm}\\
{\bf Examples of exact inf-convolution:} We now mention different
cases where the inf-convolution has a solution.
\begin{itemize}
\item First, when both convex functions $g^A$ and $g^B$ are
dilated, then their inf-convolution is exact without having to
impose any particular assumption, as we have already noticed when
working with static risk measures in the first part (see
Proposition \ref{Theorem dilated risk measure}). More precisely,
assume that $g^A$ and $g^B$ are dilated from a given convex
function $g$ such that $g^A=g_{\gamma_A}$ and $g^B=g_{\gamma_B}$,
then $g^A\,\square\,g^B=g_{\gamma_A+\gamma_B}$ and for any $z$, an
optimal solution $x^*$ to the inf-convolution problem is given by
$x^*=\frac{\gamma_B}{\gamma_A+\gamma_B}$. \item More generally, if
$g^A$ is bounded from below and if $g^B$ satisfies the
qualification constraint ensuring that $\inf_z g^B(z)$ is reached
for some $z$ (in other words, $g^B$ has a strictly positive
recession function $g^B_{0^+}$), then the condition
(\ref{condition inf-convolution atteint}) is satisfied and the
inf-convolution $g^A \square g^B$ has a non-empty compact set of
solutions.
\end{itemize}
\subsubsection{Characterization of Optima}
We are now interested on the characterization of optima in the
case of exact inf-convolution. This can be done in terms of the
subdifferentials of the different convex functions involved. More
precisely, let us consider $z^A$ and $z^B$ respectively in ${\rm
dom}(g^A)$ and in ${\rm dom}(g^B)$ and $z=z^A+z^B$ in ${\rm
dom}(g^A \square g^B)$. Then, $\partial g^A(z^A) \cap \partial
g^B(z^B) \subset \partial (g^A \square
g^B)(z)$.\\
Moreover, if $\partial g^A(z^A) \cap \partial g^B(z^B) \neq
\emptyset$, then the inf-convolution $g^A \square g^B$ is exact at
$z=z^A+z^B$ and $\partial g^A(z^A) \cap \partial g^B(z^B) =
\partial (g^A \square g^B)(z)$. (For more details, please refer to
\cite{Hiriart-Urruty-Lemarchal}).\vspace{1mm}\\
In particular, as $0$ belongs to the domain of $g^A$ and $g^B$, if
$\partial g^A(0) \cap
\partial g^B(0) \neq \emptyset$, then $\partial (g^A \square
g^B)(0)=\partial g^A(0)
\cap \partial g^B(0)$ and the inf-convolution is exact at $0$.\\
Moreover, if both functions are centered, i.e. $g^A(0)=g^B(0)=0$,
then the inf-convolution is also centered as $(g^A \square
g^B)(0)=g^A(0)+g^B(0)=0$.
\subsubsection{Regularization by Inf-Convolution}
As convolution, the infimal convolution is used in regularization
procedures. The most famous regularizations are certainly, on the
one hand, the {\it Lipschitz regularization} $g_{(k)}$ of $g$
using the inf-convolution with the kernel $b_k(z)= k|z|$ and on
the other hand, the {\it differentiable regularization}, also
called {\it Moreau-Yosida regularization}, $g_{[k]}$ of $g$ using
the inf-convolution with the kernel $q_k(z)= \frac{k}{2}|z|^2$.
Both regularizations do not have however the same
"efficiency".\vspace{-2mm}
\paragraph{\sc Lipschitz Regularization}
We first consider the inf-convolution $g_{(k)}$ of $g$ using the
kernel $b_k(z)= k|z|$ or more generally using functions whose
polar's domain is bounded (or equivalently with a finite recession
function).\\
The function $g_{(k)}$ is finite, convex, non decreasing w.r. to
$k$. Moreover,  its inf-convolution $g_{(k)}$ is
Lipschitz-continuous, with Lipschitz constant $k$. More generally,
{\it the inf-convolution of two convex functions, one of them
satisfying (H1), also satisfies (H1) without any condition on the
other function.}\\
If $z_0\in {\rm int\>dom}g$, then $g_{(k)}(z_0)=g(z_0)$ for $k$
large enough.
When $g=l_C$ is the indicator function of a closed convex set $C$, $g_k=k {\rm dist}(.,C)$.\\
This regularization is used in the book's chapter dedicated to
BSDEs to show the existence of BSDE with continuous coefficient.
\vspace{-2mm}
\paragraph{\sc Moreau-Yosida Regularization}
We now consider the inf-convolution $g_{[k]}$ of $g$ using the
kernel $q_k(z)= \frac{k}{2}|z|^2$. The function $g_{[k]}$ is
finite, convex, non decreasing w.r. to $k$. Moreover, $g_{[k]}$ is
differentiable and its gradient is Lipschitz-continuous with
Lipschitz constant $k$. In other words, the polar function of
$g_{[k]}$ is \index[sub]{strongly convex}strongly convex with
module $k$, equivalently $G_{[k]}(.)-\frac{k}{2}|.|^2$ is still a
convex function (for more details,
please refer to \index[aut]{Cohen} Cohen \cite{Cohen}).\\
There exists a point $J_k(z)$ that attains the minimum in the
inf-convolution problem with $q_k$. The maps $z\rightarrow J_k(z)$
are Lipschitz continuous with a constant $1$, independent of $k$
and monotonic in the following sense
$(J_k(z)-J_k(y))\,\,(z-y)^*\geq ||J_k(z)-J_k(y)||^2$.
Moreover,$\nabla g_{[k]}=k(z-J_k(z))$.\\
More generally, {\it the inf-convolution of two convex functions,
one of them being strongly convex, satisfies (H3) without any
condition on the other function.}

\end{document}